\documentclass[11pt]{amsart}

\title{Higher Spherical Scissors Congruence I: Hopf Algebra}

\author[Klang]{Inbar Klang}
\address[Klang]{
Department of Mathematics,
Vrije Universiteit Amsterdam - Faculty of Science,
De Boelelaan 1111,
1081 HV Amsterdam,
The Netherlands}
\email{i.klang@vu.nl}

\author[Kuijper]{Josefien Kuijper}
\address[Kuijper]{
Department of Mathematics, University of Toronto
Bahen Centre, Room 6290
40 St. George St.,
Toronto, Ontario M5S 2E4,
Canada}
\email{josefien.kuijper@utoronto.ca}

\author[Malkiewich]{Cary Malkiewich}
\address[Malkiewich]{
Department of Mathematics, Binghamton University,
PO Box 6000, 
Binghamton, New York 13902,
USA}
\email{cmalkiew@binghamton.edu}

\author[Mehrle]{David Mehrle}
\address[Mehrle]{
Department of Mathematics, 
University of Kentucky, 
719 Patterson Office Tower, 
Lexington, Kentucky 40506,
USA}
\email{davidm@uky.edu}

\author[Wittich]{Thor Wittich}
\address[Wittich]{
Department of Mathematics, Universit\"{a}t Osnabr\"{u}ck, 
Albrechtstrasse 28a, 49076 Osnabr\"{u}ck,
Germany}
\email{thor.wittich@uni-osnabrueck.de}

\usepackage{adjustbox}
\usepackage{amssymb,amsthm,amsmath} % standard stuff
\usepackage{bbm} % blackboard bold 1
\usepackage{enumerate} % enumerate tools
\usepackage[margin=1.25in]{geometry} % margins and stuff
\usepackage[pagebackref,
			colorlinks,
			citecolor=OI3,
			linkcolor=OI3,
			urlcolor=OI5]{hyperref} % clickable hyperlinks
\usepackage{mathrsfs} % for \mathscr
\usepackage{thmtools} % this fixes \cref for arXiv. No clue why.
\usepackage{tikz} % tikz pictures
\usepackage{tikz-cd} % commutative diagrams 
\usepackage{todonotes} % to do notes
\usepackage{xcolor} % colors

\usepackage[nameinlink,capitalize]{cleveref} % smart inline references

\definecolor{OI1}{RGB}{230,159,0}	% orange
\definecolor{OI2}{RGB}{86,180,233}	% light blue
\definecolor{OI3}{RGB}{0,158,115}	% green
\definecolor{OI4}{RGB}{240,228,66}	% yellow
\definecolor{OI5}{RGB}{0,114,178}	% dark blue
\definecolor{OI6}{RGB}{213,94,0}	% terracotta/reddish brown
\definecolor{OI7}{RGB}{204,121,167}	% pink

\newcommand{\Q}{\mathbb{Q}}
\newcommand{\R}{\mathbb{R}}
\newcommand{\Sph}{\mathbb{S}}
\newcommand{\Z}{\mathbb{Z}}

\DeclareMathAlphabet\mathbfcal{OMS}{cmsy}{b}{n} 

\newcommand{\sma}{\wedge} % smash product
\newcommand{\sfrac}[2]{{}^{#1}\!/_{\!#2}} % slanted fraction

\makeatletter
\newcommand*{\boxwedge}{\mathbin{\mathpalette\@boxwedge{}}}
\newcommand*{\@boxwedge}[2]{%
  \sbox0{$#1\boxtimes\m@th$}%
  \dimen2=.5\dimexpr\wd0-\ht0-\dp0\relax % side bearing
  \dimen@=\dimexpr\ht0+\dp0\relax
  \def\lw{.07} % line width as factor for height of \boxtimes
  \kern\dimen2 % side bearing
  \tikz[
    line width=\lw\dimen@,
    line join=round,
    x=\dimen@,
    y=\dimen@,
    baseline=0
  ]
  \draw
    (\lw/2,0) rectangle (1-\lw,1-\lw)
    (\lw,0) -- (.5,1-\lw-\lw/2) -- (1-\lw-\lw/2,0)
  ;%
  \kern\dimen2 % side bearing
}
\makeatother

\newcommand{\apt}{\mathrm{apt}} % apartment
\DeclareMathOperator{\CT}{CT} % cone of tits complex 
\DeclareMathOperator{\C2T}{C^2T} % double cone of tits complex 
\DeclareMathOperator{\Ls}{LS} % lee-szczarba group
\DeclareMathOperator{\Pt}{Pt} % polytope group
\newcommand{\Sah}{\mathcal{S}} % spectral Sah algebra
\newcommand{\funSah}{\mathfrak{S}} % spectral Sah functor
\DeclareMathOperator{\St}{St} % steinberg module
\DeclareMathOperator{\ST}{ST} % suspension of tits complex 
\DeclareMathOperator{\tPt}{\widetilde{Pt}} % reduced polytope 
\newcommand{\tcP}{\widetilde{\mathcal{P}}} % mathcal p with tilde
\DeclareMathOperator{\T}{T} % tits complex
\DeclareMathOperator{\GL}{GL} % general linear group

\DeclareMathOperator{\colim}{colim} % colimit
\DeclareMathOperator{\ev}{ev} % evaluation
\DeclareMathOperator{\hocolim}{hocolim} % homotopy colimit
\newcommand{\id}{\mathrm{id}} % identity
\DeclareMathOperator{\pr}{pr} % projection
\newcommand{\sd}{\operatorname{sd}} % subdivision
\DeclareMathOperator{\sgn}{sgn} % sign
\DeclareMathOperator{\sh}{sh} % shear map
\newcommand{\Step}{\textup{Step}} % step functions
\newcommand{\Sub}{\textup{Sub}} % subspaces
\DeclareMathOperator{\Sym}{Sym} % symmetric power

\DeclareMathOperator{\Ab}{Ab} % abelian groups 
\DeclareMathOperator{\CMon}{CMon} % commutative monoids
\DeclareMathOperator{\Comm}{Comm} % commutative operad
\DeclareMathOperator{\CHopf}{CHopf} % commutative hopf algebras
\DeclareMathOperator{\Dip}{Dip} % finite dim inner product spaces
\newcommand{\Fun}{\mathrm{Fun}} % functors
\newcommand{\Pol}[2]{\mathcal{P}^{#1}_{#2}} % polytopes in geometry #1 with isometry group #2
\DeclareMathOperator{\Sp}{Sp} % spectra
\DeclareMathOperator{\Spp}{\mathbf{Sp}} % bold spectra
\newcommand{\Top}{\mathrm{Top}} % unpointed spaces
\newcommand{\topp}{\mathrm{Top}_*} % pointed spaces

\numberwithin{equation}{section}
\numberwithin{figure}{section}

\newtheorem{corollary}[equation]{Corollary}
\newtheorem{lemma}[equation]{Lemma}
\newtheorem{proposition}[equation]{Proposition}
\newtheorem{theorem}[equation]{Theorem}

\newtheorem{lettertheorem}{Theorem} 
\theoremstyle{definition}
\newtheorem{conjecture}[equation]{Conjecture}
\newtheorem{construction}[equation]{Construction}
\newtheorem{definition}[equation]{Definition}
\newtheorem{example}[equation]{Example}
\newtheorem{notation}[equation]{Notation}
\newtheorem{remark}[equation]{Remark}
\newtheorem*{question}{Question}

\crefname{equation}{}{} % make \cref behave like eqref for equations
\crefname{figure}{figure}{figures} % make \cref write `figure` instead of `fig`

\begin{document}

%%%%%%%%%%%%%%%%
%%% ABSTRACT %%% 
%%%%%%%%%%%%%%%%
\begin{abstract}
	In the study of the generalization of Hilbert's Third Problem to spherical geometry, Sah constructed a Hopf algebra of spherical polytopes with product given by join and coproduct given by a generalized Dehn invariant. 
	Using Zakharevich's reinterpretation of scissors congruence via algebraic K-theory, we lift the Sah algebra to an $(E_\infty, E_1)$-Hopf algebra spectrum whose $\pi_0$ is the classical Sah algebra. 
	As an application, we show that the reduced spherical scissors congruence $K$-theory groups $\widetilde K_{2n}\big(\mathcal{P}^{S^{2k+1}}_{O(2k+2)}\big)$ are nonzero for all nonnegative integers $n$ and $k$. 
\end{abstract}

\maketitle

\begingroup%
\setlength{\parskip}{0em} % no paragraph skips in TOC
\setcounter{tocdepth}{1}
\tableofcontents
\endgroup%

%%%%%%%%%%%%%%%%%%%%
%%% INTRODUCTION %%%
%%%%%%%%%%%%%%%%%%%%
\section{Introduction}

Two polytopes $P, Q$ in $n$-dimensional Euclidean space are said to be \textit{scissors congruent} if it is possible to decompose $P$ into finitely many polytopes and assemble $Q$ out of the pieces. 
More precisely, $P$ and $Q$ are scissors congruent if one can write $P = \cup_{i=1} ^k P_i$ and $Q = \cup_{i=1} ^k Q_i$, where all the intersections $P_i \cap P_j$ and $Q_i \cap Q_j$ have measure zero, and there is an isometry $P_i \cong Q_i$ for every $i$. 
Such polytopes have the same volume.

In dimension 2, $P$ and $Q$ are scissors congruent if and only if they have the same area. 
Hilbert's third problem asks about dimension 3: if $P$ and $Q$ have the same volume, are they necessarily scissors congruent?
The answer is no; an additional scissors congruence invariant is given by the Dehn invariant
\[
	D(P) = \sum_e l(e) \otimes \theta(e) \in \R \otimes \R / \pi \Z,
\]
where the sum ranges over the edges of $P$, $l(e)$ denotes the length of $e$, and $\theta(e)$ denotes the dihedral angle of $P$ at $e$. 
Dehn invariants also exist in higher dimensions, and the Dehn--Sydler theorem \cite{jessen68, jessen72,Sydler} states that in dimensions 3 and 4, $P$ and $Q$ are scissors congruent if and only if their volumes and Dehn invariants agree. 

In addition to Euclidean space of all dimensions, one can also consider other geometries, such as spherical and hyperbolic geometry. 
In fact, the dihedral angle in the Dehn invariant is most naturally considered as an element of a (reduced) spherical scissors congruence group. 
Given an $n$-dimensional geometry $X$, the \textit{scissors congruence group} $\mathcal{P}(X)$ is defined to be the free abelian group on ($n$-dimensional) polytopes $P \subseteq X$, modulo the relations $[P] = [g(P)]$ for every isometry $g$ of $X$, and $[P] = \sum_{i=1} ^k [P_i]$ if $P = \cup_{i=1} ^k P_i$ with $P_i \cap P_j$ of measure zero for all $i \neq j$. In the case of the sphere, the \textit{reduced scissors congruence group} $\tcP(S^{n-1})$ is the quotient of $\mathcal{P}(S^{n-1})$ by suspensions of polytopes. 
(The suspension of a polytope $P \subseteq S(V)$ for a codimension-one subspace $V \subseteq \R^n$ is formed by taking its join with the two points in $S^{n-1}$ that are perpendicular to $V$.)

The volume is then a group homomorphism $\mathcal{P}(X) \to \mathbb{R}$, and the Dehn invariants are homomorphisms of the form 
\(
	\mathcal{P}(X^n) \to \mathcal{P}(X^{n-c}) \otimes \tcP(S^{c-1}),
\) 
where $X^{n-c}$ is the same type of geometry as $X = X^n$, but of dimension $n-c$.
Dupont and Sah \cite{ds1} classified the kernel of the Dehn invariant in 3-dimensional spherical geometry and 3-dimensional hyperbolic geometry, and posed a generalization of Hilbert's third problem:

\begin{question}[Generalized Hilbert's Third Problem] 
	Do volume and the Dehn invariants completely classify polytopes up to scissors congruence in Euclidean, spherical, and hyperbolic geometries of all dimensions?
\end{question}

This question is often studied via algebraic structures induced by the Dehn invariant and other polytope operations. 
Sah \cite{sah_79} gave a Hopf algebra structure on the direct sum $\bigoplus_{n\geq0}  \widetilde{\mathcal{P}}(S^{n-1})$, which was revisited through the lens of shuffle algebras by Cathelineau in \cite{cathelineau_sc}. 
We call this Hopf algebra the \textit{Sah algebra}. 
The coproduct in the Sah algebra is given by the Dehn invariants. 
In addition, the direct sum of scissors congruence groups in all dimensions (whether in Euclidean, spherical, or hyperbolic geometry) is a comodule over the Sah algebra. 
The Sah algebra also shows up in the following formulation of a conjecture due to Goncharov.
\begin{conjecture}[Goncharov]
	The Hopf algebra cohomology 
	\(
		H^i\big(
			\bigoplus_{n\geq0}\widetilde{\mathcal{P}}(S^{n-1})
		\big)_n
	\) 
	is isomorphic to the positive eigenspace of the complex conjugation action on $K_{2n-i}(\mathbb{C})^{(n)}$, where the latter is the weight $n$ part of the Adams decomposition of algebraic K-theory.
\end{conjecture}

In this paper, we lift the Hopf algebra structure of 
\(
	\bigoplus_{n\geq0}  \widetilde{\mathcal{P}}(S^{n-1})
\) 
to reduced scissors congruence $K$-theory spectra. 
First defined by Zakharevich in \cite{inna-assem}, higher scissors congruence $K$-theory groups are the homotopy groups of a space built out of polytopes and scissors congruences between them. 
For example, $\pi_0$ of this space is precisely $\mathcal{P}(X)$, while $\pi_1$ of this space describes scissors congruence automorphisms of polytopes. 
The space of polytopes and scissors congruences between them is in fact an infinite loop space, and therefore gives rise to a spectrum $K(\Pol{X}{I(X)})$, where $I(X)$ denotes the group of isometries of $X$. 
Higher scissors congruence $K$-theory groups of $X$ are the homotopy groups of this spectrum.

As in many $K$-theoretic endeavors, considering $K$-theory as a spectrum (rather than a sequence of groups or a space) is  beneficial. 
For example, the third named author proved in \cite{scissors_thom} that scissors congruence $K$-theory spectra are Thom spectra of certain bundles over recognizable base spaces. 
The descriptions of the base spaces, and the fact that Thom spectra are twisted suspensions, allow for calculations of scissors congruence $K$-groups that were previously out of reach.

Let us denote by $\widetilde K(\Pol{S^{n-1}}{O(n)})$ the reduced spherical scissors congruence $K$-theory spectrum in dimension $n$. 
In this paper, we prove:

\begin{lettertheorem}[\Cref{cor:spectral_Sah_alg_is_Hopf}]\label{thm-hopf-intro}
    The spectrum 
    \(
    	\Sah 
			:= 
		\bigvee_{n \geq 0} 
			\widetilde K(\Pol{S^{n-1}}{O(n)})
	\) 
	is a commutative Hopf algebra spectrum.
\end{lettertheorem}

\vspace*{-1em}

\noindent In addition, we prove that this spectrum $\Sah$ is indeed a lift of the Sah algebra:

\begin{lettertheorem}[\Cref{cor:hopf_alg_struct_agree}]\label{thm-pi0-intro}
    As Hopf algebras,
    \(
    	\pi_0(\Sah) 
    		\cong 
		\bigoplus\limits_{n\geq0} 
			\widetilde{\mathcal{P}}(S^{n-1})
	\).
\end{lettertheorem}

\vspace*{-1em}

In light of \cref{thm-pi0-intro}, we call $\Sah$ the \emph{spectral Sah algebra}. We also show that the $E_1$-coalgebra structure of the spectral Sah algebra cannot be improved to an $E_2$-coalgebra structure (\cref{not_cocomm}).
 
\begin{remark}
	Note that related Hopf algebra structures on the level of group homology were recently obtained. 
	Indeed, the third author in \cite{scissors_thom} shows that if we are working rationally, $\pi_m(\widetilde K(\Pol{S^{n-1}}{O(n)}))$ agrees with the twisted group homology $H_m(O(n),\St(\R^n)^t)$ where $O(n)$ is considered as a discrete group. 
	Therefore Cathelineau's Hopf algebra structure from \cite{cathelineau_sc} gives a Hopf algebra structure on the $K$-groups $\pi_m(\widetilde K(\Pol{S^{n-1}}{O(n)}))$ after rationalization. 
	Furthermore, both Ash, Miller and Patzt \cite{AshMillerPatzt}, and Brown, Chan, Galatius and Payne \cite{BrownChanGalatiusPayne} produce Hopf algebra structures on 	\(
		H_*(\GL(\mathbb{Z}),\St(\Q)) 
			= 
		\bigoplus_n H_*(\GL_n(\mathbb{Z}), \St_n(\Q)).
	\) 
	Actually, Ash, Miller and Patzt produce a Hopf algebra structure on $H_*(\GL(\mathbb{Z}),\St(\Q) \otimes \chi)$ for any character $\chi$. 
	In particular, their result also allows $\St(\Q)^t$ as coefficients. 
	Although the comparison to our work still needs to be explored further, their product map agrees with ours.
\end{remark}
 
In \cref{sec:application} we use the bigraded Hopf algebra structure on the rational homotopy groups of $\Sah$ to construct a large nonzero subalgebra.

\begin{lettertheorem}[\Cref{cor:nonzero_subalg}]
	The rational homotopy groups
	\begin{equation*}
		\bigoplus_{m,n \geq 0}
			\widetilde K_m(\mathcal{P}^{S^{n-1}}_{O(n)}) \otimes \Q
	\end{equation*}  
contain as a subalgebra the free commutative algebra on 
\(
	\bigoplus_{k \geq 0} \Lambda^{2k+1}(\R/\Q)[2k,2].
\)
\end{lettertheorem}

In future work, we hope to prove that scissors congruence $K$-theory spectra for Euclidean, spherical, and hyperbolic geometries give rise to comodule spectra over this Hopf algebra. 
We aim to use these additional algebraic structures to make new computations of scissors congruence $K$-theory. 

To resolve the coherence issues that arise in the proof of  \cref{thm-hopf-intro}, we lift the spectral Sah algebra 
\(
	\Sah 
		= 
	\bigvee_{n \geq 0} 
		\widetilde K(\Pol{S^{n-1}}{O(n)})
\) 
to a diagram $\funSah: \Dip \to \Sp^O$, where $\Sp^O$ denotes the category of orthogonal spectra and $\Dip$ denotes the category of finite-dimensional inner product spaces with inner product preserving linear isomorphisms between them. 
%By work of the third author \cite{scissors_thom}, we know that 
%\(
%	\widetilde K(\Pol{S^{n-1}}{O(n)}) ) 
%		\simeq 
%	(\mathbb{S}^{-\R ^n} \wedge \Sigma\!\ST(\R^n))_{hO(n)}.
%\) 
%In words, it is the homotopy orbits of a desuspension of the doubly-suspended Tits complex $\Sigma\!\ST(\R^n)$. 
%(For a definition of $\Sigma\!\ST$, see \cref{df:SigmaST(V)}.) 
A crucial observation (\cref{prop-hocolimS}) is that the spectral Sah algebra $\Sah$ is equivalent to the homotopy colimit of 
\(
	\funSah: \Dip \to \Sp^O,
\) 
where
\[
	\funSah(V) = \mathbb{S}^{-V} \wedge \Sigma\!\ST(V).
\]
In light of this observation, we call $\funSah$ the \emph{spectral Sah functor}. 

The category of Dip-diagrams, $\Fun(\Dip, \Sp^O)$, has a symmetric monoidal product given by Day convolution. 
In \cref{sec:Sah_functor_is_bialg}, we prove that direct sum of inner product spaces makes $\funSah$ into a commutative monoid under Day convolution, and the Dehn invariant makes $\funSah$ into an $E_1$-coalgebra under Day convolution.
Furthermore, these structures are compatible, so that $\funSah$ is an $(E_\infty, E_1)$-bialgebra. 

Showing that $\funSah$ is a Hopf algebra takes considerably more work. 
The key observation here is that being a Hopf algebra is a property of a bialgebra, not structure on it. 
A bialgebra is Hopf if and only if its \emph{shear map} 
%\(
%	(\mu \otimes \id) \circ (\id \otimes \delta)
%\) 
is an isomorphism \cite[Proposition 5.6]{KKMMW0}. 
To show that the shear map for $\funSah$ is an isomorphism,   we 
%use polytope groups and Lee--Szczarba groups in \cref{sec:apartments} to 
build a $\Dip$-diagram $\tPt \colon \Dip \to \Ab$ of abelian groups whose colimit is the (ordinary, non-spectral) Sah algebra.
In \cref{lem:polytope_functor_hopf}, we show that this diagram is a Hopf algebra in $\Fun(\Dip, \Ab)$. 
Since $\tPt$ is a Hopf algebra, its shear map is an isomorphism. 
Keeping careful track of the isomorphism $H_0(\funSah) \cong \tPt$, we prove in that the shear map for $\funSah$ is also an isomorphism, and therefore $\funSah$ is a Hopf algebra (\cref{Hopfupourlife}). 

Using tools from \cite{KKMMW0}, a colimit argument shows that the spectral Sah algebra $\Sah = \hocolim\funSah$ is a commutative Hopf algebra spectrum in the $\infty$-category of spectra. We must use the $\infty$-category of spectra because coalgebras in any symmetric monoidal model category of spectra are necessarily cocommutative by \cite[Theorem 1.1]{PS_coalgebras}; we show that the spectral Sah algebra is not cocommutative in \cref{not_cocomm}. 
Finally, we show that the induced Hopf algebra structure on $\pi_0(\Sah)$ is the classical one.

\begin{remark}
	In \cite{cz-hilbert} Campbell and Zakharevich define a ``derived Dehn invariant'' at the level of pointed simplicial sets, which becomes the classical Dehn invariant after taking coinvariants and homology. 
	In \cref{subsec:CZ_comparison}, we also explain how their procedure produces a map of spectra, and demonstrate that our coproduct contains their derived Dehn invariant as a summand.
\end{remark}

\textbf{Acknowledgments.}
	The authors would like to thank the organizers of the Collaborative Research Workshop on K-theory and Scissors Congruence at Vanderbilt University in July 2024, during which this paper was initiated, and the conference Scissors congruence and K-theory at the University of Pennsylvania in July 2025, during which this paper was almost finished. 
	Both were sponsored by the NSF grant FRG: Collaborative Research: Trace Methods and Applications for Cut-and-Paste K-Theory. 
	The authors also thank Alexander Kupers for an important insight that one should prove the Hopf algebra structure at the level of diagrams, and Maximilien Peroux for extensive guidance on model categories of coalgebras, and how to ultimately avoid using them.

	JK was partially supported by the Knut and Alice Wallenberg grant KAW 2023.0416. 
	CM was partially supported by the National Science Foundation (NSF) grants DMS-2052923 and DMS-2506430 and by a Simons Fellowship.
	DM was partially supported by NSF grant DMS-2135884.
	TW was partially supported by the research training group 2240: Algebro-geometric Methods in Algebra, Arithmetic and Topology. 

%%%%%%%%%%%%%%%%%%%%%%%%%%%%%%%%
%%% THE SPECTRAL SAH ALGEBRA %%%
%%%%%%%%%%%%%%%%%%%%%%%%%%%%%%%%
\section{The spectral Sah algebra}

In this section we recall the classical definition of the classical Sah algebra $\bigoplus_{n\geq0}  \widetilde{\mathcal{P}}(S^{n-1})$. 
We also develop a spectral model $\Sah$ alongside it, such that $\pi_0(\Sah)= \bigoplus_{n\geq0}  \widetilde{\mathcal{P}}(S^{n-1})$ as abelian groups. 
In the subsequent sections we construct the algebraic structure on $\Sah$ that descends to the Hopf algebra structure on $\bigoplus_{n\geq0}  \widetilde{\mathcal{P}}(S^{n-1})$, which justifies us calling $\Sah$ the ``spectral Sah algebra.''

%%% TITS COMPLEX AND SUSPENSIONS
\subsection{The Tits complex and its suspensions}

For any real vector space $V$ of positive finite dimension, the \textit{Tits complex} $\T(V)$ is the realization of the poset of proper nonzero subspaces of $V$. We write this as:
\[ 
	\T(V) = | 0 \subsetneq U \subsetneq V |. 
\]
Note that the realization of any poset can be constructed without using degeneracy maps, as the realization of a $\Delta$-complex or a semisimplicial set. 
It has a $k$-simplex for every \emph{strict} flag $U_0 \subsetneq U_1 \subsetneq \cdots \subsetneq U_k$, and each face of this $k$-simplex is another strict flag.

Recall that the cone $\CT(V)$ is formed by taking the product $\T(V) \times I$ of $T(V)$ with an interval $I$ and quotienting $\T(V) \times \{1\}$ to a single point. 
This is homeomorphic to the realization of the larger poset in which we have formally added an initial or terminal object, corresponding to the cone point. 
As a result, we can think of the cone $\CT(V)$ as the realization of the poset $| 0 \subsetneq U \subseteq V |$, and the double cone $\C2T(V) = \mathrm{C}(\CT(V))$ as the realization of the poset $| 0 \subseteq U \subseteq V |$.

\begin{definition}\label{df:SigmaST(V)}
    For any finite-dimensional real vector space $V$, we define the \emph{doubly-suspended Tits complex} $\Sigma\!\ST(V)$ to be the realization of the poset of all subspaces of $V$, quotiented by the union of the realizations of the posets that do not involve 0 or $V$:
\[ 
	\Sigma\!\ST(V) 
		= 
	\frac{| 0 \subseteq U \subseteq V |}%
		{| 0 \subseteq U \subsetneq V | \cup | 0 \subsetneq U \subseteq V |}. 
\]
\end{definition}

\begin{lemma}
    When $\dim V > 0$, this space $\Sigma\!\ST(V)$ is in fact homeomorphic to the reduced suspension of the unreduced suspension of $\T(V)$.
\end{lemma}

\begin{proof}
    The unreduced suspension $\ST(V)$ is the quotient of the cone $\CT(V)$ by $\T(V)$, so it is homeomorphic to the quotient of realizations of posets
	\[ 
		\ST(V) 
			= 
		\frac{| 0 \subsetneq U \subseteq V |}%
			{| 0 \subsetneq U \subsetneq V |}. 
	\]
    
    When we take the reduced suspension of this, we use the general fact that for any space $X$ and subspace $A$, we have an identification
    \[ 
    	\Sigma\left( \frac{X}{A} \right) 
			\cong 
		\frac{CX}{CA \cup_A X}. 
	\]
    Taking $X = | 0 \subsetneq U \subseteq V |$, $A = | 0 \subsetneq U \subsetneq V |$, and forming cones by adding the initial object 0 to the posets gives the result.
\end{proof}

The case where $V = 0$ is treated separately. 
Recall that for any space $X$, the quotient of $X$ by the empty set $\varnothing$ does not give $X$, but rather $X$ with a disjoint basepoint added:
\[ 
	X/\varnothing 
		= 
	\textrm{cofiber}(\varnothing \to X) \cong X_+. 
\]
As a result, we have
\[ 
	\Sigma\!\ST(0) 
		= 
	\frac{| 0 \subseteq U \subseteq 0 |}%
		 {
		  | 0 \subseteq U \subsetneq 0 | 
			\cup 
		  | 0 \subsetneq U \subseteq 0 |
		 } 
		\cong 
	\frac{*}{\varnothing} \cong S^0. 
\]
In general, the space $\Sigma\!\ST(V)$ is homotopy equivalent to a wedge of spheres whose dimension matches $V$. 
We will recall this in more detail in \cref{solomon_tits_w_apts} below.

\begin{definition}\label{classical_sah_algebra}
    Let $n = \dim V$. 
    The \emph{Steinberg module} $\St(V)$ is defined to be the $n$-th homology of $\Sigma\!\ST(V)$:
    \[ 
    	\St(V) = H_n(\Sigma\!\ST(V),*). 
	\]
\end{definition}    
When $n \geq 2$, this is also isomorphic to $H_{n-2}(\T(V))$.
    
Note that $\St(\R^n)$ carries an action by the orthogonal group $O(n)$, regarded as a discrete group acting on the poset of subspaces of $\R^n$. 
We take this action and ``twist'' it by tensoring it with the determinant representation $\det$, which is $\Z$ with $O(n)$ acting by multiplication by the determinant. 
Note that the tensor $\St(\R^n) \otimes \det$ is also sometimes written as $\St(\R^n)^t$, e.g. \cite{cathelineau_sc,scissors_thom}.

\begin{definition}\label{sah_algebra}
	Our first definition of the Sah algebra is as the direct sum of the groups of $O(n)$-coinvariants:
	\[ 
		\bigoplus_{n \geq 0} (\St(\R^n) \otimes \det)_{O(n)}. 
	\]
\end{definition}

This might be more appropriately called the \emph{Cathelineau algebra}, since Cathelineau used this definition to define the Hopf algebra structure algebraically in \cite[Section 8]{cathelineau_sc}, while Sah's original definition in \cite[Section 6]{sah_79} was more geometric. 
We recall in \cref{sec:apartments} why these agree with other definitions in the literature, such as \cite{dupont_book} and \cite{rudenko}.

%% PASSING TO SPECTRA 
\subsection{Passing to spectra}
\label{sec:intro_to_spo}

% \begin{definition}\label{spo}
Let $\Sp^O$ denote the category of orthogonal spectra \cite{mandell_may_shipley_schwede}. 
This is a closed symmetric monoidal category, whose tensor is the smash product $\sma$ and whose unit is the sphere spectrum $\Sph = \Sigma^\infty S^0 = F_0 S^0$. 
By standard abuse of notation, we also use $\sma$ to denote the operation that takes a spectrum $X$ and a based space $K$, and smashes each level of $X$ with $K$.
% \end{definition}

We regard orthogonal spectra in the coordinate-free way, so that they have a level for each finite-dimensional real inner product space $V$. 
For any such $V$, there is an evaluation functor
\[ 
	\ev_V\colon \Sp^O \to \topp 
\]
that takes each spectrum $X$ to the based space $X(V)$.

% \begin{definition}\label{free_spectrum}
The left adjoint of this functor is the \emph{free spectrum} $F_V A$ on a based space $A$ at the vector space $V$. 
Concretely, this is an orthogonal spectrum whose $W$-th level is $\mathscr J(V,W) \sma A$, for the category $\mathscr J$ defined in \cite[Definition II.4.1]{mandell_may}.
% \end{definition}

Let $O(V)$ be the group of linear isometries of $V$, with left action $O(V) \times V \to V$. 
It turns out that $\mathscr J(V,V) = O(V)_+$, and therefore the free spectrum $F_V A$ carries a natural right $O(V)$-action by precomposition
\[ 
	\mathscr J(V,W) \sma \mathscr J(V,V) \sma A \to \mathscr J(V,W) \sma A. 
\]
It also turns out that $\mathscr J(0,V) = S^V$, the one-point compactification of $V$, and that the left $O(V)$-action on $S^V$ can be described as the composition
\[ 
	\mathscr J(V,V) \sma \mathscr J(0,V) \to \mathscr J(0,V). 
\]

\begin{definition}\label{positive_and_negative_spheres}
    Let $V$ be any finite-dimensional real inner product space. We define two spectra with left action of $O(V)$:
    \vspace*{-1em}
    \begin{itemize}
        \item Let $\Sph^V = \Sigma^\infty S^V = F_0 S^V$. 
        	  We let $O(V)$ act on the left through its action on $V$.
        \item Let $\Sph^{-V} = F_V S^0$. 
        	  We let $O(V)$ act on the left by applying the inverse and then applying the right action on $F_V S^0$.
    \end{itemize}
\end{definition}

\begin{lemma}
    As spectra with left $O(V)$-action, $\Sph^V$ and $\Sph^{-V}$ are inverses under the smash product. 
    Accordingly, there is an $O(V)$-equivariant map of spectra
    \[ 
    	\Sph^V \sma \Sph^{-V} \to \Sph 
	\]
    that is a stable equivalence after forgetting the $O(V)$-action.
\end{lemma}

\begin{proof}
    The map is
    \[ 
    	F_0 S^V \sma F_V S^0 \cong F_V S^V \to F_0S^0, 
	\]
    where the first isomorphism is an instance of \cite[Lemma II.4.8]{mandell_may} and the second map is the canonical stable equivalence from \cite[Definition III.4.3]{mandell_may}, given by the composition
    \[ 
    	\mathscr J(V,W) \sma \mathscr J(0,V) \sma S^0 
			\to 
		\mathscr J(0,W) \sma S^0. 
	\]
    If we act on the source of this map by an element $\rho \in O(V)$, it has the effect of adding in two copies of $\mathscr J(V,V)$, selecting the points $\rho^{-1}$ and $\rho$ inside, then composing everything together. 
    But the composition cancels the $\rho^{-1}$ and $\rho$ out, so this agrees with the trivial action in the target, which proves that the map is $O(V)$-equivariant.
\end{proof}

In fact, these actions are continuous if we give $O(V)$ its usual topology, but in this paper we will only be interested in the discrete topology on $O(V)$. 

Now take $V = \R^n$ with the usual action of $O(n)$.

\begin{definition}\label{spectral_sah_algebra}
    We define an orthogonal spectrum called the \emph{spectral Sah algebra} by
    \[ 
    	\Sah 
			= 
    	\bigvee_{n \geq 0} 
			EO(n)_+ 
				\sma_{O(n)} 
			(\Sph^{-\R^n} \sma \Sigma\!\ST(\R^n)). 
	\]
    where $EO(n)$ is any contractible free right $O(n)$-CW complex, and $\Sph^{-\R^n}$ is the de-suspension by the $O(n)$-representation $\R^n$ given in \cref{positive_and_negative_spheres}. 
    This is exactly the notion of homotopy orbits in the category of spectra with $O(n)$-action, so we abbreviate the above definition as
    \[ 
    	\Sah 
			= 
		\bigvee_{n \geq 0} \ST(\R^n)^{1-\R^n}_{hO(n)}. 
	\]
    In the $n = 0$ term, \cref{df:SigmaST(V)} gives $\Sigma\!\ST(0) \cong S^0$, so that we get the sphere spectrum
    \[ 
    	(\Sigma\!\ST(\R^0))_{hO(0)} \cong \Sph^0 = \Sph. 
	\]
\end{definition}
%We will see that the spectrum $\Sph^{-\R^n} \sma \Sigma\!\ST(\R^n)$ is equivalent to a wedge sum of sphere spectra $\Sph^0$ by the Solomon-Tits theorem (\cref{solomon_tits_w_apts}). 
It is a consequence of \cref{solomon_tits_w_apts} below that $\pi_0(\Sah)$ is isomorphic to the classical Sah algebra:
\[ 
	\pi_0(\Sah) = \bigoplus_{n \geq 0} (\St(\R^n) \otimes \det)_{O(n)}, 
\] 
as abelian groups.
Our goal is to prove that $\Sah$ is an $(E_\infty,E_1)$-Hopf algebra in spectra, that on $\pi_0$ recovers the Hopf algebra structure on the classical Sah algebra from \cite{sah_79,cathelineau_sc,dupont_book,rudenko}.

\begin{remark}
This definition of the spectral Sah algebra has not appeared in the literature yet, but it is motivated by \cite{scissors_thom}, which shows that it is the wedge sum of the natural spectral versions of the groups $\widetilde{\mathcal P}(S^{n-1})$, called the reduced spherical scissors congruence $K$-theory:
\[ 
	\widetilde K(\Pol{S^{n-1}}{O(n)}) 
		\simeq 
	\ST(\R^n)^{1-\R^n}_{hO(n)}. 
\]
Therefore, we may also write 
\(
	\Sah 
		\simeq 
	\bigvee_{n \geq 0}
		\widetilde{K}(\mathcal{P}^{S^{n-1}}_{O(n)}).
\)
\end{remark}

\begin{remark}
The homotopy of the spectral Sah algebra is not known, even rationally on $\pi_0$. 
We do know that all of the odd summands are rationally trivial both on $\pi_0$ and above (see \cref{odd_vanishing}). 
However the rational homotopy type of the rest is not completely understood, not even the $n = 4$ summand on $\pi_0$, which has to do with three-dimensional spherical geometry. 
There is a four-term exact sequence for it, see \cite[(7.15)]{dupont_book}:
\[
\begin{tikzcd}
    0 
    	\ar{r} 
		& 
	H_3(SU(2);\Q) 
		\ar{r} 
		& 
	\tcP(S^3) \hspace{-1.2pt}\otimes\hspace{-1pt} \Q 
		\ar{r} 
		& 
	\R/\Q \hspace{-1pt}\otimes\hspace{-1pt} \R/\Q 
		\ar{r} 
		& 
	H_2(SU(2);\Q) 
		\ar{r} 
		& 
	0
\end{tikzcd}
\]
However, $H_3(SU(2);\Q)$, where $SU(2)$ has the discrete topology, is difficult to completely understand.
\end{remark}

%The Hopf algebra structure does give us a large subalgebra of $\Sah$ rationally: the classes in $\widetilde K_{2m}(S^1) \cong \Lambda^{2m+1}(\R/\Q)$ are all primitive for degree reasons, and therefore the free commutative algebra on these classes is a subalgebra of the spectral Sah algebra (see \cref{cor:nonzero_subalg} below). 
%The resulting classes are given geometrically by iterated joins of arcs lying in perpendicular copies of $S^1$ inside $S^{2d-1}$, and higher scissors congruence classes that come from applying interval exchange transformations to these arcs.

%%% A DIFFERENT MODEL OF SIGMA ST V
\subsection{A different model of \texorpdfstring{$\Sigma\!\ST(V)$}{∑ST(V)}}\label{Section 2.3}

Let $V$ be a real vector space of positive finite dimension. We denote by $\Sub(V)$ the poset of linear subspaces of $V$ ordered by inclusion (including $0$ and $V$).
 
\begin{definition}
    Let $\textup{Poset}\big([0,1], \Sub(V)\big)$ be the set of order-preserving functions from $[0,1]$, with the usual ordering, to $\textup{Sub(V)}$. 
    Let $\Step\big([0,1],\Sub(V)\big)$ be the quotient
    \[
    	\textup{Poset}\big([0,1], \Sub(V)\big)/\sim 
	\]
	where $f\sim g$ if the set of points where $f$ and $g$ differ is finite.  
\end{definition}

Note that an order-preserving function $[0,1]\to \Sub(V)$ can only have a finite number of values: $f(1)$ has to be a subspace of $V$, and the values that $f$ takes on $[0,1)$ form a non-degenerate flag of subspaces of $f(1)$, which has to have finite length since $V$ is finite-dimensional. 
We can think an element $\phi$ of $\Step\big([0,1],\Sub(V)\big)$ as a subdivision of $[0,1]$ into a finite number of intervals $[0,s_1]$, $[s_1,s_2],\dots, [s_{k},1]$, and a flag of subspaces $V_0\subseteq \dots \subseteq V_k$, where  the interior of each interval $[s_i,s_{i+1}]$ is sent to $V_i$ (where we set $s_0=0$ and $s_{k+1} = 1$), and we are agnostic about what happens to the cut points $s_0,\dots, s_{k+1}$. 
Note that giving the inner cut points $0 < s_1 <  \dots < s_{k}<1$ is equivalent to giving interval lengths $t_0,\dots,t_k>0$ such that $t_0 = s_1, t_0+t_1 = s_2, \dots, t_0 + \dots +t_k = 1$. 

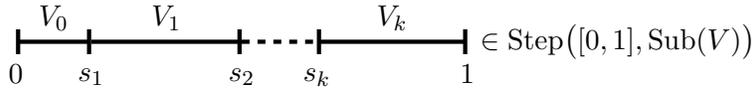
\begin{figure}[h]
\begin{center}
	\begin{tikzpicture}
		\node at (8,0) {$ \in \Step\big([0,1],\Sub(V)\big)$};
		\draw[ultra thick,|-|] 
			(0,0) node[below=1.5mm]{$0$} 
				--node[above]{$V_0$} 
			(1,0);
        \draw[ultra thick,-|] 
        	(1,0) node[below=2mm]{$s_1$} 
				--node[above]{$V_1$} 
			(3,0) node[below=2mm]{$s_2$};
		\draw[ultra thick,dashed]
			(3,0) -- (4,0);
		\draw[ultra thick,|-|]
			(4,0) node[below=2mm]{$s_k$}
				--node[above]{$V_k$} 
			(6,0) node[below=1.5mm]{$1$};
	\end{tikzpicture}
\end{center}
\caption{A schematic illustration of an element of $\Step\big([0,1],\Sub(V)\big)$.} 
\end{figure}
 
\begin{definition}
    We define a function 
    \[
    	c:{\Step}\big([0,1],\Sub(V)\big) 
			\to 
		|0 \subseteq U \subseteq V|
	\]
    that sends 
    \(
    	\phi \in {\Step}\big([0,1],\Sub(V)\big),
	\) 
	determined by the flag $U_0\subseteq \dots \subseteq U_k$, and interval lengths $t_0, \dots, t_k$ with $t_0+\dots + t_k = 1$, to the point in the simplex $\langle U_0,\dots, U_k\rangle$ with coordinates $t_0,\dots, t_k$. 
	In particular, the constant step function with value $U$ is sent to the 0-simplex in $|0 \subseteq U \subseteq V|$ associated to the space $U$.
\end{definition}

\begin{proposition}\label{prop:set_model}
    The function $c$ induces a bijection on the quotients
	\[
		\tilde c:{\Step}\big([0,1],\Sub(V)\big)/{\sim} 
			\to  
		\Sigma\!\ST(V)
	\] 
	where on the left we identify a function $\phi$ with the base point if $\phi(t) \neq 0$ for all $t>0$, or if $\phi(t) \neq V$ for all $t<1$.
\end{proposition}

The step functions $\phi$ identified with the basepoint are the functions whose associated flags of subspaces do not start with $0$ or end with $V$. 

\begin{proof}
	We first check that the original function $c$ is a bijection. 
	Indeed, every point of $|0 \subseteq U \subseteq V|$ is in the interior of a $k$-simplex for some $k\geq 0$, and hence in the image of $c$. 
	On the other hand, if $c(\phi) = c(\psi)$ is in the interior of the simplex $\langle U_0,\dots, U_k\rangle$ with coordinates $t_0,\dots, t_k$, then $\phi,\psi$ can both be represented by poset functions with inner cut points $t_0, t_0+t_1,\dots, t_0+\dots+t_{k-1}$ and values $U_0, \dots, U_k$. 
	Therefore they are the same as elements of $\Step\big([0,1],\Sub(V)\big)$.

	It is automatic that $c$ passes to a bijection on the quotient: if $\phi$ is a step function such that  $\phi(t) \neq 0$ for all $t>0$, or such that $\phi(t) \neq V$ for all $t<1$, then this corresponds to a point in $|0 \subseteq U \subseteq V|$ in which $U_0\neq 0$ or $U_k \neq V$.
\end{proof}
        
%%%%%%%%%%%%%%%%%%%%%%%%%%%%%%%%%%%%%%%%
%%% DIP DIAGRAMS AND DAY CONVOLUTION %%%
%%%%%%%%%%%%%%%%%%%%%%%%%%%%%%%%%%%%%%%% 
\section{\texorpdfstring{$\Dip$}{Dip}-diagrams and Day convolution}
\label{sec-dip-diagrams}

In order to define the product and coproduct on the spectral Sah algebra, we will present the spectral Sah algebra as a homotopy colimit of a diagram of spectra, indexed by a certain small category. 
We will define a product and coproduct on this diagram, where the symmetric monoidal structure on the category of such diagrams is given by Day convolution.
Our small category will be called $\Dip$.

\begin{definition}
    The category $\Dip$ has an object for each finite-dimensional real inner product space $V$, and a morphism for each linear isomorphism $V \cong W$ preserving the inner product.
\end{definition}

This is considered as a discrete category, so the morphisms do not have a topology. 
A skeleton is given by the set of vector spaces $\R^n$ for $n \geq 0$, with only automorphisms, given by $O(n)$ as a discrete group. 
Moreover, the direct sum (of inner product spaces on objects, and of linear maps on morphisms) defines a symmetric monoidal structure on $\Dip$. 
The unit for this symmetric monoidal structure is the inner product space $0$.

Recall from \cref{sec:intro_to_spo} that we let $\Sp^O$ denote the category of orthogonal spectra. 
Our category of interest is the category of functors from $\Dip$ to orthogonal spectra, denoted $\Fun(\Dip, \Sp^O)$.

\begin{remark}
	The category $\Fun(\Dip, \Sp^O)$ is similar to the functor categories used for orthogonal calculus, which hints at intriguing connections. 
	It is also reminiscent of the category of FI-modules of \cite{ChurchEllenbergFarb}, as well as the category of VB-modules of \cite{JS}. 
	The category of FI-modules is the category of functors from FI (the category of finite sets and injections) to a module category. 
	The category of VB-modules is the category of functors from VB (the category of finite dimensional vector spaces and linear bijections) to a module category. 
	We might therefore call $\Fun(\Dip, \Sp^O)$ the category of ``VO-modules.'' 
	These are functors from VO (the category of finite dimensional vector spaces and orthogonal isomorphisms) to a module category. 
	In our case, the target category is the category of modules over the (orthogonal) sphere spectrum $\mathbb{S}$.
\end{remark}

The category $\Sp^O$ has a stable model structure that was first defined in \cite[Theorem 9.2]{mandell_may_shipley_schwede}, which makes $\Sp^O$ into a cofibrantly generated model category. 
A standard result in model category theory says that $\Fun(\Dip, \Sp^O)$ is therefore a cofibrantly generated model category as well, with the projective model structure, in which weak equivalences and fibrations are determined objectwise.

Let $\Fun(\Dip, \Sp^O)^c$ denote the cofibrant diagrams in the projective model structure, and let $\Fun(\Dip, \Sp^O)^{\operatorname{pc}}$ denote the \emph{pointwise cofibrant} diagrams, for which $F(V)$ is a cofibrant spectrum for each $V \in \Dip$ separately. 
Then we have the strict inclusion
\[ 
	\Fun(\Dip, \Sp^O)^{\operatorname{c}}
		\subsetneq 
	\Fun(\Dip, \Sp^O)^{\operatorname{pc}}. 
\]

As we will need to consider commutative monoids in this category, we also recall the positive stable model structure from \cite[Theorem 14.2]{mandell_may_shipley_schwede}. 
This has the same weak equivalences as the stable model structure, but fewer cofibrations, and is therefore Quillen equivalent to the stable model structure. 
We let $\Sp^O_+$ denote the category of orthogonal spectra with the positive stable model structure, and $\Fun(\Dip, \Sp^O_+)$ denote the model category of functors $\Dip \to \Sp^O_+$, with the corresponding projective model structure. 
So we have an additional strict inclusion of model categories
\[ 
	\Fun(\Dip, \Sp^O_+)^{\operatorname{c}}
		\subsetneq 
	\Fun(\Dip, \Sp^O)^{\operatorname{c}}
		\subsetneq 
	\Fun(\Dip, \Sp^O)^{\operatorname{pc}}. 
\]

Since both $\Dip$ and $\Sp^O$ are symmetric monoidal, the category $\Fun(\Dip, \Sp^O)$ also has a symmetric monoidal product given by Day convolution.

\begin{definition}\label{day_convolution}
    For any $F, F': \Dip \to \Sp^O$, define their \emph{Day convolution} 
    \[
    	F \boxwedge F': \Dip \to \Sp^O
	\] 
	as the left Kan extension 
	\(
		\mathrm{Lan}_{\oplus}(F \wedge F'). 
	\)
    This is the left Kan extension of the following diagram:
	\[
	\begin{tikzcd}[column sep = 3.5em]
    	\Dip \times \Dip 
	        \ar[r,"F \times F'"]
    	    \ar[d, "\oplus"] 
        	& 
	    \Sp^O \times \Sp^O 
    	    \ar[r, "\wedge"] 
        	& 
	    \Sp^O
    	    \\
	    \Dip 
	    	\ar[urr,dashed, 
				out=0, 
				in=-150,
				swap,
				"\mathrm{Lan}_{\oplus}(F \wedge F')"
				]
	\end{tikzcd}
	\]
\end{definition}

We can alternatively describe $F \boxwedge F'$ explicitly as follows:
\[
	(F \boxwedge F')(V) 
		= 
	\bigvee_{U \subseteq V} F(U) \wedge F'(V\ominus U),
\]
where $V \ominus U = U^\perp$ denotes the orthogonal complement of $U$ in $V$. 
More generally, for $F_1, ..., F_n \in \Fun(\Dip, \Sp^O)$, 
\[
	(F_1 \boxwedge \ldots \boxwedge F_n)(V) 
		= 
	\bigvee_{0 = U_0 \subseteq U_1 \subseteq \ldots \subseteq U_{n-1} \subseteq U_n = V} 
		F_1(U_1 \ominus U_0) \wedge \ldots \wedge F_n(U_n \ominus U_{n-1}).
\]
The unit for $\boxwedge$ is given by the functor $\mathbbm{1}: \Dip \to \Sp^O$ with $\mathbbm{1}(0) = \mathbb{S}$ and $\mathbbm{1}(V) = {*}$ for $V \neq 0$.

\begin{remark}\label{rem-symmetry-explicit}
    Explicitly, the symmetry isomorphism $(F \boxwedge F')(V) \to (F' \boxwedge F)(V)$ sends the $U$-summand $F(U) \sma F'(U^\perp)$ to the $U^{\perp}$-summand $F'(U^\perp) \sma F((U^\perp)^\perp) = F'(U^\perp) \sma F(U)$ by the swap map.
\end{remark}

\begin{proposition}\label{prop-day-sm}
   The model categories $\Fun(\Dip, \Sp^O_+)$ and $\Fun(\Dip, \Sp^O)$ are symmetric monoidal topological model categories under Day convolution.
\end{proposition}
% That is, fibrations and weak equivalences are determined objectwise. We note that the same is true for $\Fun(\Dip, \topp)$.

\begin{proof}
    The fact that the model structure is symmetric monoidal follows from \cite[Theorem 4.1]{BB}; see also \cite[Theorem 5.18]{Whi}, or \cite[Theorem 6.5 and Lemma 6.6]{mandell_may_shipley_schwede}. 
    To see that it is also topological, we construct a strong symmetric monoidal Quillen left adjoint from topological spaces, by taking each topological space $K$ to the diagram $\Sigma^\infty_+ K \sma \mathbbm{1}$, which at the vector space 0 is the spectrum $\Sigma^\infty_+ K$, and at all other vector spaces is the zero spectrum. 
    Note that taking Day convolution of $X$ with this diagram has the effect of smashing all of the spectra in $X$ with the based space $K_+$. 
    One checks that this functor is both strong symmetric monoidal and left Quillen, and the result follows.
\end{proof}

In particular, both of the above model categories are simplicial model categories. 
The simplicial set of maps between two objects is formed by taking singular simplices of the topological space of maps. 
This means that all of the simplicial sets are Kan complexes, which is a hypothesis in \cite[Theorem 5.7]{KKMMW0} which we use in \cref{sec:hop_functor}.

Day convolution also preserves weak equivalences between pointwise cofibrant objects. 
Recall that a diagram in $\Fun(\Dip,\Sp^O)$ is pointwise cofibrant if each spectrum $F(V)$ is cofibrant in the stable model structure on orthogonal spectra. 

\begin{proposition}\label{prop-ptwise-cof}
    If $F, F', G, G'\in \Fun(\Dip, \Sp^O)$ are pointwise cofibrant and $F \to G$ and $F' \to G'$ are pointwise equivalences of diagrams, then $F \boxwedge F' \to G \boxwedge G'$ is a pointwise equivalence as well.
\end{proposition}

\begin{proof}
    Since $(F \boxwedge F')(V) = \bigvee_{U \subseteq V} F(U) \wedge F'(U^\perp)$, this follows from the fact that the wedge sum and smash product of spectra preserve weak equivalences between cofibrant spectra.
\end{proof}

\begin{proposition}\label{prop-colim-left}
   The functor $\colim \colon \Fun(\Dip, \Sp^O) \to \Sp^O$ is a left Quillen functor.
\end{proposition}

\begin{proof}
    Recall that the fibrations and weak equivalences in $\Fun(\Dip, \Sp^O)$ are determined objectwise. 
    Therefore the constant diagram functor
    \[ 
    	\Sp^O \to \Fun(\Dip, \Sp^O), 
	\]
    preserves both fibrations and acyclic fibrations, so it is right Quillen. 
    Since the colimit is the left adjoint of this, the colimit is therefore left Quillen.
\end{proof}

Our main functor of interest is the following. 
Let $\Sph^{-V} = F_V S^0$ be the orthogonal spectrum that is the shift desuspension of $S^0$ by the vector space $V$. 
As an orthogonal spectrum, it is the functor (free diagram) represented by $V$.

\begin{definition}\label{sahfunctor}
    The \emph{spectral Sah functor} is the functor $\funSah$, whose value on an inner product space $V$ is
    \[ 
    	\funSah(V) 
			= 
		\mathbb{S}^{-V} \wedge \Sigma\!\ST (V) \cong F_V \Sigma\!\ST (V), 
	\]
    and whose value on a linear isomorphism preserving the inner product $V \to W$, is the natural map of spectra 
    \(
    	\mathbb{S}^{-V} \wedge \Sigma\!\ST (V) 
			\to 
		\mathbb{S}^{-W} \wedge \Sigma\!\ST (W).
	\)
\end{definition}

Here we take as our model for $\Sigma\!\ST (V)$ the quotient of geometric realizations of posets
\[ 
	\Sigma\!\ST(V) 
		= 
	\frac{| 0 \subseteq U \subseteq V |}{| 0 \subseteq U \subsetneq V | \cup | 0 \subsetneq U \subseteq V |} 
\]
as in \cref{df:SigmaST(V)}. 
Note that $\funSah$ is a pointwise cofibrant $\Dip$-diagram. 
The spectral Sah functor $\funSah$ models the spectral Sah algebra in the following sense:

\begin{proposition}\label{prop-hocolimS}
    The homotopy colimit of $\funSah$ is equivalent to the spectral Sah algebra.
\end{proposition}

\begin{proof}
    The category $\Dip$ is equivalent to the groupoid whose objects are $\mathbb{N}$ and in which morphisms from $n$ to itself are $O(n)$ with the discrete topology. 
    All other morphism sets are empty.
    This groupoid is the disjoint union $\amalg _n \mathcal{B}O(n)$, where $\mathcal{B} G$ denotes the category associated to the group $G$; it has one object with morphisms given by $G$.
    Thus
    \[
    	\underset{\Dip}\hocolim\, 
		\funSah 
			\simeq 
		\bigvee_n 
		\underset{\mathcal{B}O(n)}\hocolim\,  
		(\mathbb{S}^{-\R^n} \wedge \Sigma\!\ST (\R^n)) 
			\simeq 
		\bigvee_n (\mathbb{S}^{-\R^n} 
		\wedge \Sigma\!\ST (\R^n))_{hO(n)}
		= \Sah
	\]
    as required.
\end{proof}

%%%%%%%%%%%%%%%%%%%%%%%%%%%
%%% BIALGEBRA STRUCTURE %%%
%%%%%%%%%%%%%%%%%%%%%%%%%%% 
\section{Bialgebra structure on the spectral Sah functor}\label{sec:Sah_functor_is_bialg}

In this section, we define a product and a coproduct on the spectral Sah functor $\funSah$. 
More precisely, we show that $\funSah$ is both a commutative monoid and a coalgebra over the little intervals operad $(D_1)_+$ in the symmetric monoidal category $\textup{Fun}(\Dip,\Sp^O)$. 
At the end of \cref{sec:hop_functor} we will use infinity-categorical techniques to show that, after taking the homotopy colimit over $\Dip$, this indeed gives a bialgebra structure on $\Sah$. 

We first reduce the problem to defining the product and coproduct structure on the following functor.

\begin{definition}
	\label{topological-Sah-functor}
	Let $G \colon \Dip \to \topp$ denote the functor $V\mapsto \Sigma\!\ST(V)$. 
\end{definition}

\begin{definition}\label{def-desusp}
    Let 
    \[
    	D \colon \Fun(\Dip, \topp ) \to \Fun(\Dip, \Sp^O)
	\]
    be the functor that sends $F \colon \Dip \to \topp$ to the functor given by 
    \[
    	V \mapsto \Sph^{-V} \sma F(V)
	\]
    where $\Sph^{-V} = F_VS^0$ is the negative $V$-sphere from \cref{positive_and_negative_spheres}.
\end{definition}

\begin{lemma}\label{negative_sphere_coherence}
    The assignment $V \mapsto \Sph^{-V}$ extends to a strong symmetric monoidal functor $\Dip \to \Sp^O$.
\end{lemma}

\begin{proof}
    This is a standard fact about orthogonal spectra, but we will still explain the proof. 
    We have an isomorphism $i_{U,V}\colon \Sph^{-U} \sma \Sph^{-V} \cong \Sph^{-(U \oplus V)}$ that arises from the maps $\mathscr J(U,W_1) \sma \mathscr J(V,W_2) \to \mathscr J(U \oplus V,W_1 \oplus W_2)$ that add the embeddings and the vectors in the orthogonal complement together, see e.g. \cite[Lemma II.4.8]{mandell_may}. 
    We also have a unique isomorphism of orthogonal spectra $\Sph^{-0} \cong \Sph$. 
    The interested reader may check the following associativity, unitality, and symmetry axioms, which verify that we have a strong symmetric monoidal functor:
    \[ 
    	\begin{tikzcd}[column sep = 4em]
	        \Sph^{-U} \sma \Sph^{-V} \sma \Sph^{-W} 
    	    	\rar{i_{U,V} \sma \id} 
				\dar[swap]{\id \sma i_{V,W}} 
				& 
			\Sph^{-(U \oplus V)} \sma \Sph^{-W} 
				\dar{i_{U \oplus V,W}} 
				\\
	        \Sph^{-U} \sma \Sph^{-(V \oplus W)} 
	        	\rar{i_{U,V \oplus W}} 
				& 
			\Sph^{-(U \oplus V \oplus W)}
	    \end{tikzcd}
    \]
    \vspace{1em}
    \[
   		\begin{tikzcd}
	        \Sph \sma \Sph^{-V} 
   		     	\rar{\cong} 
				\dar[swap]{\cong} 
				& 
			\Sph^{-(0)} \sma \Sph^{-V} 
				\dar{i_{0,V}} 
				\\
	        \Sph^{-V} 
    	    	\rar{\cong} 
				& 
			\Sph^{-(0 \oplus V)}
    	\end{tikzcd}
	    \hspace{2em}
    	\begin{tikzcd}
	        \Sph^{-U} \sma \Sph^{-V} 
	        	\dar[swap]{i_{U,V}} 
				\rar{\cong} 
				& 
			\Sph^{-V} \sma \Sph^{-U} 
				\dar{i_{V,U}} 
				\\
	        \Sph^{-(U \oplus V)} 
	        	\rar{\cong} 
				& 
			\Sph^{-(V \oplus U)}
	    \end{tikzcd}
	\]
\end{proof}

Note that the functor category $\Fun(\Dip, \topp)$ has a symmetric monoidal structure coming from the symmetric monoidal structure $(\topp, \wedge, S^0)$ on $\topp$ and Day convolution. 

\begin{lemma}\label{lem:D_symm_monoidal}
    The functor $D$ is strong symmetric monoidal for the Day convolution product.
\end{lemma}

\begin{proof}
    Given $F,F'$ in $\textup{Fun}(\Dip, \topp )$ and $V$ in $\Dip$, there is a natural isomorphism 
    \[
    	(DF\boxwedge DF')(V) \to D(F\boxwedge F')(V).
	\]
    Indeed, there is a natural isomorphism
    \[
    	\bigvee_{U\subseteq V}\Sph^{-U} \sma F(U)\wedge\Sph^{-(V\ominus U)} \sma  F'(V\ominus U) 
			\to 
		\Sph^{-V} \sma \left(\bigvee_{U\subseteq V} F(U)\wedge F'(V\ominus U)\right)
	\]
    since $\Sph^{-U} \sma \Sph^{-(V\ominus U)} \simeq \Sph^{-V}$ and the smash product commutes with the wedge sum. 
    Similarly there is a natural isomorphism
    \[
    	\mathbbm{1}_{\Sp}(V) \cong D(\mathbbm{1}_{\topp})(V)
	\]
    given by the unique isomorphism of zero spectra $* \cong *$ when $V \neq 0$, and the unique isomorphism of orthogonal spectra $\Sph \cong \Sph^{-0} \wedge S^0$ when $V = 0$.
    
    The associativity, unitality, and symmetry axioms for this strong symmetric monoidal functor then follow from those same conditions for the functor in \cref{negative_sphere_coherence}.
\end{proof}

Recall the functor $G \colon \Dip \to \topp$ from \cref{topological-Sah-functor}. Observe that $\funSah = D(G)$. 
Therefore, to build commutative monoid and $(D_1)_+$-coalgebra structures on $\funSah$, it suffices to show that $G$ is a commutative monoid and a $(D_1)_+$-coalgebra in $\textup{Fun}(\Dip, \topp)$. 

%%% THE PRODUCT 
\subsection{The product}\label{sec:product}

In order to show that $G$ is a commutative monoid for Day convolution, we need to give a map 
\[
	G \boxwedge G \to G,
\]
or equivalently, maps 
\[ 
	(G \boxwedge G)(V) 
		= 
	\bigvee_{U\subseteq V} G(U) \wedge G(V \ominus U) \to G(V).
\]
Therefore it suffices to define continuous maps 
\[
	G(V) \wedge G(V') \to G(V \oplus V')
\]
for all $G,V$ in $\Dip$. 

\begin{construction}\label{cons:product}
	We recall that 
	\[  
		G(V) = \Sigma\!\ST(V) 
			= 
		\frac{| 0 \subseteq U \subseteq V |}
			 {| 0 \subseteq U \subsetneq V | 
			 	\cup 
			  | 0 \subsetneq U \subseteq V |} 
	\]
	We will use the shorthand notation $P(V)$ for the poset of spaces $0 \subseteq U \subseteq V$ and denote by $R(V)$ its realization. 
	Similarly, we denote by $P_<(V)$ the sub-poset of proper subspaces, and denote by $R_<(V)$ its realization. 
	We denote by $P_>(V)$ the sub-post of non-zero subspaces, and denote by $R_>(V)$ its realization. 

	Using the identity
	\[ 
		\frac{X \times Y}{(A \times Y) \cup (X \times B)} 
			\cong 
		\frac{X}{A} \sma \frac{Y}{B}, 
	\]
	we see that we need to give a continuous map 
	\[
		\frac{R(V)\times R(V')}%
			 {((R_>(V)\hspace{-1pt}
			 	  \cup\hspace{-1pt} 
				R_<(V))\times R(V')) \hspace{-1pt}
				  \cup\hspace{-1pt} 
			   (R(V) \times ( R_>(V') \hspace{-1pt}
			      \cup\hspace{-1pt} R_<(V')))
			  } 
			\to 
		\frac{R(V\oplus V')}%
			 {R_>(V \oplus V') 
			 	\hspace{-1pt}\cup\hspace{-1pt} 
			  R_<(V\oplus V')}.
	\]
	Now, since geometric realization preserves finite limits, we get a map 
	\[
		R(V)\times R(V') \to R(V\oplus V') 
	\]
	from the poset map 
	\begin{align*}
		P(V) \times P(V') 
			&
			\to 
		P(V\oplus V') 
			\\
	    (U,U') 
	    	& 
			\mapsto 
		U\oplus U'.
	\end{align*}
	It is clear that the sub-posets $P_>(V) \times P(V')$ and $P(V) \times P_>(V')$ land in $P_>(V\oplus V')$ and $P_<(V)\times P(V')$ and $P(V) \times P_<(V')$ land in $P_<(V \oplus V')$, so this gives the desired map $G(V)\wedge G(V')\to G(V \oplus V')$.
\end{construction}

\begin{proposition}\label{prop-comm-step}
    The functor $G$ is a commutative monoid object in $\Fun(\Dip, \topp)$.
\end{proposition}

\begin{proof}
    We need to provide a unit for the multiplication as well.
    Recall that the monoidal unit in  $\Fun(\Dip, \topp)$ is the functor $\mathbbm{1}$ such that $\mathbbm{1}(0) = S^0$ and $\mathbbm{1}(V) = *$ for $V \neq 0$. 
    The unit map $\eta \colon \mathbbm{1} \to G$ is a homeomorphism in degree 0 (recall that $\Sigma ST(0)=S^0$), and inclusion of the basepoint in all other degrees. 

    Associativity and commutativity follow since this is true for the direct sum of vector spaces. 
    That $\eta$ is indeed the unit follows since $V\oplus 0 = V$. 
 \end{proof}
 
\begin{corollary}\label{cor-sah-comm}
    The functor $\funSah = D(G)$ is a commutative monoid in $\Fun(\Dip, \Sp^O)$.
\end{corollary}

\begin{figure}[!ht]
\begin{center}
    \begin{tikzpicture}[xscale=2]
        % pale vertical lines 
        \draw[gray,dotted] (0.015,0) -- (0.015,-3);
        \draw[gray,dotted] (0.985,-1) -- (0.985,-3);
        \draw[gray,dotted] (1.985,0) -- (1.985,-3);
        \draw[gray,dotted] (2.985,-1) -- (2.985,-3);
        \draw[gray,dotted] (3.985,0) -- (3.985,-3);        
    
        % interval in F(V)
        \draw[ultra thick,OI6,|-|] (0,0) --node[above]{$V_0$} (2,0);
        \draw[ultra thick,OI6,-|] (2,0) --node[above]{$V_1$} (4,0);
        \node at (4.5,0) {$\in \textcolor{OI6}{G(V)}$};
        
        % inverval in F(U)
        \draw[ultra thick,OI2,|-|] (0,-1) --node[above]{$U_0$} (1,-1);
        \draw[ultra thick,OI2,-|] (1,-1) --node[above,fill=white]{$U_1$} (3,-1);
        \draw[ultra thick,OI2,-|] (3,-1) --node[above]{$U_2$} (4,-1); 
        \node at (4.5,-1) {$\in \textcolor{OI2}{G(U)}$};

        % interval in F(V \oplus U)
        \draw[ultra thick,|-|] (0,-3) --node[above]{$V_0 \oplus U_0$} (1,-3);
        \draw[ultra thick,-|] (1,-3) --node[above]{$V_0 \oplus U_1$} (2,-3);
        \draw[ultra thick,-|] (2,-3) --node[above]{$V_1 \oplus U_1$} (3,-3);
        \draw[ultra thick,-|] (3,-3) --node[above]{$V_1 \oplus U_2$} (4,-3); 
        \node at (4.7,-3) {$\in G(V \oplus U)$};
    \end{tikzpicture}
\end{center}
\caption{An illustration of the product $G(V) \wedge G(U) \to G(V \oplus U)$.}
\label{ProductFigure}
\end{figure}
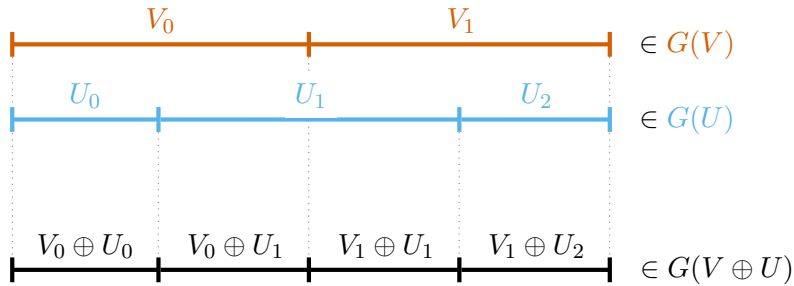

In order to be able to compare this algebra structure with the $(D_1)_+$-coalgebra structure that we will give in the next section, we now give a description of this using the step function model of $\Sigma ST(V)$ from \cref{Section 2.3}. 
We recall that by \cref{prop:set_model}, $G(V)$ can be identified with a quotient of $\Step\big([0,1],\Sub(V)\big)$.
% for a suitable equivalence relation $\sim$. 

In that model, the commutative monoid structure on $G$ is given by the coordinatewise direct sum. 
That is,
\[
	\bigvee _{U \subseteq V} G(U) \sma G(V \ominus U) \to G(V)
\]
sends 
\(
	\phi \in \Step\big([0,1],\Sub(U)\big)/\sim
\) 
and 
\(
	\psi \in \Step\big([0,1],\Sub(V \ominus U)\big)/\sim
\) 
to the equivalence class of $\phi \oplus \psi$ defined as 
\(
	(\phi \oplus \psi)(t) = \phi(t) \oplus \psi(t).
\)

%%% THE DEHN INVARIANT 
\subsection{The Dehn invariant as coproduct}

We next show that the spectral Sah functor $\funSah$ admits a coassociative coproduct. 
More precisely, we show that $\funSah$ is a coalgebra over the little intervals operad $D_1$. 
We can view $D_1$ as an $E_1$-operad in the category of pointed spaces, via the strong symmetric monoidal functor $(-)_+ \colon \Top \to \topp$ that sends an unpointed space $X$ to $X$ with a disjoint basepoint. 
We use the fact that $\Sp^O$, and therefore $\Fun(\Dip,\Sp^O)$, is naturally tensored over $\topp$. 
We denote this by $X \boxwedge F$ for a pointed space $X$ and functor $F \colon \Dip \to \Sp^O$.

\begin{definition}\label{df:little_intervals_operad}
	The \emph{little intervals operad} $D_1$ in $\Top$ is the operad where $D_1(0) = *$, and for $n\geq 1$, the space of $n$-ary operations $D_1(n)$ is the subspace of $\textup{Map}(I^{\sqcup n},I)$ given by maps
	\[
  		I \sqcup \dots \sqcup I \xrightarrow{e} I
	\]
	such that for the restriction to the $i$-th component of the disjoint union
	\[
		e_i \colon I \hookrightarrow I \sqcup \dots \sqcup I \xrightarrow{e} I
	\]
	is an embedding of the form $x\mapsto a+bx$ with $0 < b\leq 1$, and such that the interiors $e_i((0,1))$ are disjoint.
	Hence, we can identify a point in $D_1(n)$ with a tuple of a rectilinear embeddings $(e_1,\dots, e_n)$ of $I$ into $I$, with almost-disjoint images. 

	The unit 
	\(
		e \colon *\longrightarrow D_1(1)
	\)
	picks out the element $\textup{id}_I$, and the composition map 
	\[
		D_1(n)\times D_1(k_1)\times \dots \times D_1(k_n) \longrightarrow D_1(k_1+\dots + k_n) 
	\]
	sends tuples $(e_1,\dots, e_n)$ and $(e^i_1,\dots, e^i_{k_i})$ for $1\leq i\leq n$, to the tuple 
  	\[
		(e_1 \circ e^1_1, \dots, e_1 \circ e^1_{k_1}, 
		 e_2 \circ e^2_1, \dots, e_2 \circ e^2_{k_2},
		 \dots,
		 e_n \circ e^n_1, \dots, e_n \circ e^n_{k_n}).
	\]

	The $\Sigma_n$-action on $D_1(n)$ is given by permuting the tuples $e=(e_1,\dots, e_n)$. 
	For $\sigma \in \Sigma_n$, we have $\sigma^* e = (e_{\sigma(1)},\dots,e_{\sigma(n)})$.
\end{definition}

\begin{remark}\label{rmk:D_1(n)_as_set_topology}
    As a set, $D_1(n)$ can be identified with the subset of $\mathbb{R}^{2n}$ consisting of tuples $(a_1,b_1,\dots, a_n,b_n)$ such that $a_i<b_i$ for  $1 \leq i \leq n$, and such that the open intervals $(a_i,b_i)$ do not overlap for varying $i$. 
    The topology on $D_1(n)$ is the subspace topology obtained from $\mathbb{R}^{2n}$.
\end{remark}

As noted below \cref{lem:D_symm_monoidal}, it suffices to show that the functor $G \colon V\mapsto \Sigma ST(V)$ is a $(D_1)_+$-coalgebra in $\textup{Fun}(\Dip,\topp)$. We use the equivalent definition of \cref{prop:set_model}:
\[
	G \colon V \mapsto \Step\big([0,1],\Sub(V)\big)/\sim.
\] 

\begin{construction}\label{const:E_1_coalg}
	For $n\geq 0$, we define maps
	\[
		\theta_n \colon D_1(n)_+ \boxwedge G\to  G^{\hspace{1pt}\boxwedge n}
	\] 
	as follows. 
	The map 
	\[
		\theta_0 \colon D_1(0)_+ \boxwedge G \to \mathbbm{1}
	\]
	is the obvious one: evaluated on $0$ in $\Dip$, this is the identity map
	\(
		D_1(0)_+ \wedge G(0) \cong S^0 \to S^0
	\)
	and evaluated on any non-zero $V$ in $\Dip$, this is the trivial map 
	\(
		D_1(0)_+ \wedge G(V) \cong * \to S^0.
	\)

	If $n=1$, we can set $\theta_1^V(\phi) = \phi\circ e_1$, where $\theta_1^V$ is the natural transformation $\theta_1$ at $V \in \Dip$.  
	Note $\phi\circ e_1$ is equivalent to the base point of $G(V)$ if any of the cut points of $\phi$ are not in the interior of the image of $e_1$. 

	For $n> 1$, we define a map
	\(
		\theta_n \colon D_1(n)_+ \boxwedge G \to G^{\hspace{1pt}\boxwedge n}
	\)
	as follows. 
	Evaluated on $V$ in $\Dip$, it suffices to define a map
	\begin{equation}\label{eq:coproduct}
    	\theta_n^V \colon D_1(n) \times G(V) 
			\to 
		\bigvee_{0=V_0 \subseteq \dots \subseteq V_n=V} 
			G(V_1 \ominus V_0)\wedge 
				\dots 
			\wedge G(V_n \ominus V_{n-1} )
	\end{equation}
	that sends $D_1(n) \times \{*\}$ to the base point. 
	Let $e=(e_1,\dots, e_n)$ and $\phi \in G(V)$. 
	We assume that the images of the $e_i$ are in ascending order $e_1(1)<e_2(1)< \dots<e_n(1)$. 
	Indeed, if this is not the case, then $e = \sigma^*e'$ for some $\sigma \in \Sigma_n$ and a suitable tuple $e'$ for which the embeddings are in order; by equivariance we then define $\theta_n^V(e,\phi) = \sigma^*\theta^V_n(e',\phi)$.

	To determine in which component of the wedge sum $\theta^V_n(e,\phi)$ lands, we define a flag of subspaces
	\[ 
		0=U_0 \subseteq \dots \subseteq U_n = V
	\]
	by setting $U_1=\phi(e_1(1)), \dots, U_{n-1} = \phi(e_{n-1})(1)$.
	If $e_i(1)$ coincides with a cut point of $\phi$ for any $i$, then we set $\theta^V_n(e,\phi) = *$. 
	   
	Now for $1 \leq i \leq n$, we modify $\phi\circ e_i$ to obtain an element of $G(U_i\ominus U_{i-1})$. 
	Note that for all elements $t$ in the image of $e_i$, we have $\phi(t)\subseteq U_{i}$. 
	On the other hand, $U_{i - 1} = \phi(e_{i-1}(1)) \subseteq \phi(t)$ for elements $t$ in the image of $e_i$ since we assume that the images of the $e_i$ are in ascending order.  
	Therefore we can consider $\phi\circ e_i$ as an element of $\Step\big([0,1],\Sub(U_{i})\big)/\sim$, and by taking the orthogonal compliment of $U_{i-1}$ inside $U_i$, as an element of
	\[
		\Step\big(
			[0,1],\Sub(U_{i} \ominus U_{i-1})
		\big)/\sim\ \  
			\cong 
		\ G(U_{i}\ominus U_{i-1}).
	\]

	We denote this element of $ G(U_{i}\ominus U_{i-1})$ by $\phi_i$. 
	Together, these $\phi_i$ form an element $(\phi_1,\dots, \phi_n)$ in 
	\(
		\bigvee_{0=V_0 \subseteq \dots \subseteq V_n=V} 
			G(V_1 \ominus V_0)\wedge 
				\dots 
			\wedge G(V_n \ominus V_{n-1} )
	\) 
	in the component indexed by the flag $U_1 \subseteq \dots \subseteq U_n $. 
\end{construction}

\begin{figure}[ht]
\begin{center}
    \begin{tikzpicture}
        % pale lines
        \draw[dotted,gray] (0.025,0) -- (0.025,-3) -- (-1,-4);
        \draw[dotted,gray] (0.975,-1) -- (0.975,-3) -- (0,-4);
        \draw[dotted,gray] (1.975,0) -- (1.975,-3) -- (1,-4);
        \draw[dotted,gray] (1.975,-3) -- (2.025,-4);
        \draw[dotted,gray] (2.975,-1) -- (2.975,-3) -- (2.975,-4);
        \draw[dotted,gray] (3.975,-1) -- (3.975,-3) -- (3.975,-4);
        \draw[dotted,gray] (4.975,0) -- (4.975,-3) -- (4.975,-4);
%        \draw[dotted,gray] (4.975,-1) -- (4.975,-3);        
        \draw[dotted,gray] (4.975,-3) -- (6,-4);
        \draw[dotted,gray] (5.975,-1) -- (5.975,-3) -- (7,-4);        
        \draw[dotted,gray] (7.975,-1) -- (7.975,-3) -- (9,-4);        
        \draw[dotted,gray] (9.975,0) -- (9.975,-3) -- (11,-4);
    
        % interval in little disks
        \draw[ultra thick,OI6,|-|] (0,0) -- (2,0);
        \draw[ultra thick,OI6,-|] (2,0) -- (5,0);
        \draw[ultra thick,OI6,-|] (5,0) -- (10,0);
        \node at (11.15,0) {$\in \textcolor{OI6}{D_1(3)_+}$};

        % element of F(V)
        \draw[ultra thick,OI2,|-|] (0,-1) --node[above]{$0\!=\!U_0$} (1,-1);
        \draw[ultra thick,OI2,-|] (1,-1) --node[above,fill=white]{$U_1$} (3,-1);
        \draw[ultra thick,OI2,-|] (3,-1) --node[above]{$U_2$} (4,-1); 
        \draw[ultra thick,OI2,-|] (4,-1) --node[above,fill=white]{$U_3$} (6,-1); 
        \draw[ultra thick,OI2,-|] (6,-1) --node[above]{$U_4$} (8,-1); 
        \draw[ultra thick,OI2,-|] (8,-1) --node[above]{$U_5=V$} (10,-1); 
        \node at (11,-1) {$\in \textcolor{OI2}{G(V)}$};

        % combined interval
        \draw[ultra thick,OI2,|-|] (0,-3) --node[above]{$U_0$} (1,-3);
        \draw[ultra thick,OI2,-|] (1,-3) --node[above,fill=white]{$U_1$} (3,-3);
        \draw[ultra thick,OI6] (1.975,-3.15) -- (1.975,-2.85); 
        \draw[ultra thick,OI2,-|] (3,-3) --node[above]{$U_2$} (4,-3);
        \draw[ultra thick,OI2,-|] (4,-3) --node[above]{$U_3$} (6,-3);
        \draw[ultra thick,OI6] (4.975,-3.15) -- (4.975,-2.85); 
        \draw[ultra thick,OI2,-|] (6,-3) --node[above]{$U_4$} (8,-3);
        \draw[ultra thick,OI2,-|] (8,-3) --node[above]{$U_5$} (10,-3); 

        % element in F(V_1)
        \draw[ultra thick,|-|] (-1,-4) --node[below]{\scriptsize $0$} (0,-4);
        \draw[ultra thick,-|] (0,-4) --node[below]{\scriptsize $U_1$} (1,-4);
        \node at (0,-5) {$\in G(V_1)$};

        % element in F(V_3/V_1)
        \draw[ultra thick,|-|] (2,-4) --node[below]{\scriptsize $U_1/U_1$} (3,-4);
        \draw[ultra thick,-|] (3,-4) --node[below]{\scriptsize $U_2/U_1$} (4,-4);
        \draw[ultra thick,-|] (4,-4) --node[below]{\scriptsize $U_3/U_1$} (5,-4);
        \node at (3.5,-5) {$\in G(U_3\ominus U_1)$};

        % element in F(V_5/V_3)
        \draw[ultra thick,|-|] (6,-4) --node[below]{\scriptsize $U_3/U_3$} (7,-4);
        \draw[ultra thick,-|] (7,-4) --node[below]{\scriptsize $U_4/U_3$} (9,-4);
        \draw[ultra thick,-|] (9,-4) --node[below]{\scriptsize $U_5/U_3$} (11,-4);
        \node at (8.5,-5) {$\in G(U_5\ominus U_3)$};
    \end{tikzpicture}
\end{center}
\caption{An illustration of the coproduct $D(3)_+ \wedge G(V) \to G^{\hspace{1pt}\protect \boxwedge 3}(V)$. 
This particular element lands in the product 
\(
	G(U_1) \wedge G(U_3/U_1) \wedge G(U_5/U_3),
\)
which includes into $G^{\hspace{1pt}\protect \boxwedge 3}(V)$.}
\label{CoproductFigure}
\end{figure}
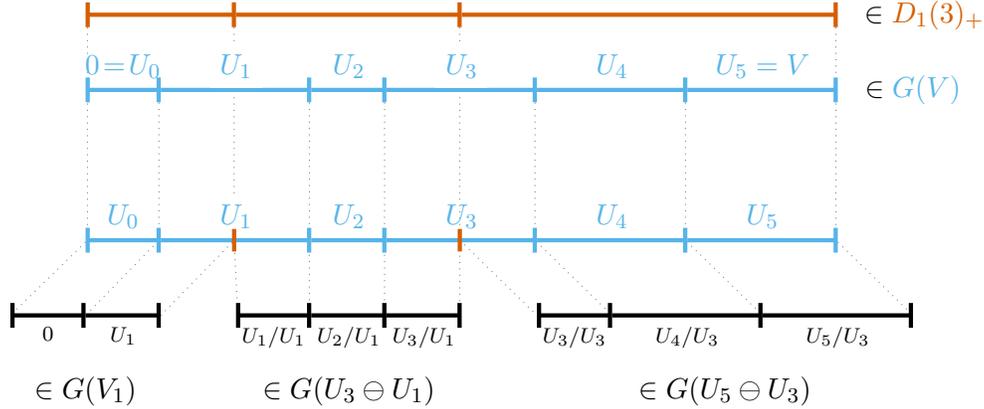

\begin{lemma}
    The structure maps as defined in \cref{const:E_1_coalg} are continuous.
\end{lemma}

\begin{proof}
    It suffices to show that the maps $\theta_n^V$ in  \cref{eq:coproduct} are continuous for all $V$. 
    Recall that 
    \[ 
    	G(V) = \Sigma\!\ST(V) 
			= 
		\frac{| 0 \subseteq U \subseteq V |}
			 {| 0 \subseteq U \subsetneq V | 
			 	\cup 
			  | 0 \subsetneq U \subseteq V |}. 
	\]
    For $k\geq 0$, let $(\Delta^k_{\textup{top}})_{U_0 \subsetneq \dots \subsetneq U_k}$ in $G(V)$ be the topological simplex associated to the strict flag $U_0 \subsetneq \dots \subsetneq U_k$ in the geometric realization $| 0 \subseteq U \subseteq V |$. 
    Since the topology on $G(V)$ is generated by these simplices, it suffices to show that the restriction of $\theta^V_n$ to each of these simplices
	\begin{equation}\label{eq:coprod_cont}
		f\colon D_1(n) \times (\Delta^k_{\textup{top}})_{U_0 \subsetneq \dots \subsetneq U_k} 
			\to 
		\bigvee_{0=V_0 \subseteq \dots \subseteq V_k=V} G(V_1)\wedge \dots \wedge G(V\ominus V_{n-1} )
    \end{equation}
    is continuous.
    
	As in \cref{rmk:D_1(n)_as_set_topology}, we identify $D_1(n)$ with the subset of $\mathbb{R}^{2n}$ of tuples $(a_1,b_1,\dots,a_n,b_n)$, such that $a_i<b_i$ for all $1 \leq i \leq n$, and the open intervals $(a_i,b_i)$ do not overlap for varying $i$. 
	We identify $(\Delta^k_{\textup{top}})_{U_0\subsetneq \dots \subsetneq U_k}$ with the standard topological $k$-simplex $\Delta^k$.
	In order to check that the map $f$ of \cref{eq:coprod_cont} is continuous, we divide $D_1(n)\times \Delta^k$ into pieces on which $f$ is given by a concrete formula. Recall that a point of $\Delta^k = (\Delta^k_{\textup{top}})_{U_0\subsetneq \dots \subsetneq U_k}$ can be thought of as a step function $\phi \colon [0,1] \to \Sub(V)$, which is determined by the flag of subspaces $U_0\subsetneq \dots \subsetneq U_k$ and cut points
	\[ 
		0 \leq s_1 \leq \dots \leq s_k \leq 1, 
	\]
	or equivalently, ``step lengths'' $t_0=s_1, t_1 = s_2-s_1, \dots, t_k = 1-s_k$, which give a point in the standard topological $k$-simplex. 
	Now what the coproduct map does for a given $e=(a_1,b_1,\dots,a_n,b_n)$ in $D_1(n)$ depends on where the cut points of $\phi$ end up relative to the intervals $[a_1,b_1],\dots, [a_n,b_n]$. 
	For simplicity of notation, we assume $n=2$. Then there are $k+1$ subsets of $D_1(2)\times \Delta^k$ where the coproduct does something non-trivial: assuming that $b_1\leq a_2,$ we need to consider the sets
	\[
		S_i = \Big\{
			\big((a_1,b_1,a_2,b_2), (s_1,\dots, s_k)\big) 
				\ \Big \vert\  
			a_1\leq s_1 \leq \dots\leq s_{i} \leq b_1 
				\text{ and } 
			a_2\leq s_{i+1} \leq \dots \leq s_k \leq b_2  
		\Big\} 
	\]
	for $0 \leq i \leq k$ (in the case that $b_1\leq a_1$, the roles of $a_1$ and $a_2$ and of $b_1$ and $b_2$ get swapped).
	These are the sets where the cut points all end up somewhere in $[a_1,b_1]$ or $[a_2,b_2]$. 
	Note that these sets are closed in $D_1(2)\times \Delta^k$.
	The closure of the complement of $S_0 \cup \dots \cup S_{k}$ on the other hand, is the closed set of configurations where at least one of the cut points is outside of, or on the boundary of one of the intervals $[a_1,b_1]$ or $[a_2,b_2]$; we denote this set by $S_\infty$. 
	This gives a covering of $D_1(2)\times\Delta^k$ by finitely many closed subsets, and we will give a formula for $f$ on each one. 

	On each of the sets $S_i$ for $0 \leq i \leq k$, we check that $f_i := f|_{S_i}$ has a linear formula and is therefore continuous. 
	For example, a point of $S_3$ looks like 
	\(
		((a_1,b_1,a_1,b_2), (s_1,\dots, s_k)),
	\) 
	where $a_1\leq s_1\leq s_2 \leq s_3 \leq b_1$ and $a_2 \leq s_4\leq \dots\leq s_k \leq b_2$. In this case $f_3$ defines a map
	\[
		f_3\colon S_3 \to G(U_3)\wedge G(V \ominus U_3).
	\]
	It rescales $(s_1,s_2,s_3)$ and $(s_4,\dots, s_k)$ so that they give the cut points of $\phi\circ e_1$ and $\phi \circ e_2$, where $e_1$ and $e_2$ are the rectilinear embeddings of $[0,1]$ in $[0,1]$ with endpoints $a_1,b_1$ and $a_2,b_2$ respectively. 
	Concretely, the formula to do this is given by
	\[ 
		\left(\frac{s_1-a_1}{b_1-a_1}, \frac{s_2-a_1}{b_1-a_1} ,\frac{s_3-a_1}{b_1-a_1}\right)	 
		\quad \textup{and} \quad 
		\left( \frac{s_4-a_2}{b_2-a_2},\dots, \frac{s_k-a_2}{b_2-a_2} \right) 
	\]
	respectively. 
	Here the relevant flags are $0=U_0\subseteq U_1 \subseteq U_2 \subseteq U_3$ and $ 0=U_3 \ominus U_3 \subseteq \dots \subseteq U_k \ominus U_3$. 
	Note that these are independent of the point in $S_3$. Since the formula is continuous in each variable (recall from  \cref{rmk:D_1(n)_as_set_topology} that $a_i < b_i$), we have established that    
	\[
		f_3\colon S_3 \to G(U_3)\wedge G(V \ominus U_3) \subseteq \bigvee_{0 \subseteq  V_1 \subseteq V} G(V_1) \wedge G(V\ominus V_1)
	\]
	is continuous. 
	Likewise, on $S_0, \dots, S_k$, the maps $f_0,\dots, f_k$ have similar formulas and are therefore continuous. 
	Lastly, on $S_\infty$, $f_\infty = f|_{S_\infty}$ takes everything to the base point of $\bigvee_{0 \subseteq  V_1 \subseteq V} G(V_1) \wedge G(V\ominus V_1)$. 

	All together this shows that the map $f$ of \cref{eq:coprod_cont} is continuous, and therefore that $\theta_n^V$ is continuous, when $n = 2$. 
	The proof for higher $n$ proceeds in the same way, with a larger collection of sets $S_i$, one for each choice of how many of the cut points are in each of the intervals $[a_j,b_j]$.
\end{proof}

\begin{proposition}\label{prop-coalginFun}
    The maps in \cref{const:E_1_coalg} make $G$ into a coalgebra over $(D_1)_+$.
\end{proposition}

\begin{proof}
	One may check that the for the identity element $\textup{id}_I$ in $D_1(1)$, the induced map $G \to G$ is the identity. 
	To see that the square 
    \begin{center}
        \begin{tikzcd}
            D_1(n)_+\boxwedge D_1(k_1)_+ \boxwedge \dots \boxwedge D_1(k_n)_+ \boxwedge G \ar[d] 
            	\ar[r] 
            	& 
            D_1(k_1)_+\boxwedge\dots \boxwedge D_1(k_n)_+\boxwedge G^{\hspace{1pt}\boxwedge n} 
            	\ar[d]
            	\\
            D_1(k_1+\dots+k_n)_+ \boxwedge G 
            	\ar[r] 
				& 
			G^{\hspace{1pt}\boxwedge k_1+\dots+k_n}
        \end{tikzcd}
    \end{center}
    commutes, all we need are the facts that composition of maps is associative, and that for real finite-dimensional vector spaces $A\subseteq B \subseteq C$, we have $C\ominus B = (C \ominus A)\ominus (B \ominus A)$.

	Lastly, we need to see that for $\sigma \in \Sigma_n$, the square 
	\begin{center}
    	\begin{tikzcd}
        	D_1(n)_+ \boxwedge G
				\arrow[r, "\theta_n"] 
				\arrow[d, "\sigma^*\boxwedge \id"]  
				& 
			G^{\hspace{1pt}\boxwedge n} 
				\arrow[d, "\sigma_*"] 
				\\
			D_1(n)_+\boxwedge G 
				\arrow[r, "\theta_n"] 
				& 
			G^{\hspace{1pt}\boxwedge n}
	    \end{tikzcd}
	\end{center}
	commutes. 
	This follows from \cref{const:E_1_coalg}. 
	Indeed, for $e=(e_1,\dots, e_n)$ in $D_1(n)_+$, we have 
	\[
		\sigma_*(\theta_n^V(e,\phi)) 
			= 
		\sigma_*(\phi_1,\dots, \phi_n) 
			= 
		(\phi_{\sigma(1)},\dots, \phi_{\sigma(n)})
			= 
		(\phi\circ e_{\sigma(1)},\ldots,\phi\circ e_{\sigma(n)}) 
	\]
	and this is equal to $\theta_n^V(\sigma^* e,\phi)$ since $\sigma^*e = (e_{\sigma(1)},\dots e_{\sigma(n)})$.    
\end{proof}

\begin{corollary}\label{cor:sah_functor_coalg}
    The spectral Sah functor $\funSah$ is a $(D_1)_+$-coalgebra in $\textup{Fun}(\Dip,\Sp^O)$.
\end{corollary}

We next show that the coalgebra structure is compatible with the product defined in the previous section. 
In particular, we show $\funSah$ is a $(D_1)_+$-coalgebra in $\textup{CMon}(\textup{Fun}(\Dip,\Sp^O))$. 
This will allow us to show that the spectral Sah functor $\funSah$, and therefore the spectral Sah algebra $\Sah$, is an $(E_\infty, E_1)$-bialgebra (\cref{Hopfupourlife}).

Again, by \cref{lem:D_symm_monoidal}, it suffices to show that $G$ is a $(D_1)_+$-coalgebra in the category $\textup{CMon}(\textup{Fun}(\Dip, \topp))$. 
We use that $\textup{CMon}(\textup{Fun}(\Dip, \topp))$ is enriched over $\topp$. 

\begin{lemma}\label{lem-bialginFun}
    The functor $G$ is a $(D_1)_+$-coalgebra in 
    \(
    	\CMon(\textup{Fun}(\Dip, \topp)).
	\)
\end{lemma}

\begin{proof}
    We have shown already that $G$ is in $\CMon(\textup{Fun}(\Dip, \topp))$. Moreover, we have maps 
    \[
    	\gamma_n \colon D_1(n)_+ \to \textup{Map}(G,G^{\hspace{1pt}\boxwedge n})
	\]
    exhibiting $G$ as a coalgebra over $(D_1)_+$, where the mapping space on the right is understood to be in $\textup{Fun}(\Dip, \topp)$. 
    These are the transpose of the maps defined in \cref{const:E_1_coalg}. 
    %$\textup{Map}_{\textup{CMod}}(G,G^{\hspace{1pt}\boxwedge n})$
    Now we show that $\gamma_n$ in fact takes values in the mapping space in $\textup{CMod}(\textup{Fun}(\Dip, \topp))$.
    Let $e = (e_1,\dots, e_n)$ be an element of $D_1(n)_+$, where the $e_i\colon[0,1]\to [0,1]$ are rectilinear embeddings and set $a_i=e_i(0)$ and $b_i = e_i(1)$. 
    We show that the associated co-operation 
    \(
    	\gamma_n(e) \colon G \to G^{\hspace{1pt}\boxwedge n}
	\)
	is a map of commutative monoids. 
	We consider the diagram
	\begin{center}
    	\begin{tikzcd}[column sep = huge]
        	G\boxwedge G 
				\arrow[r, 
					"\gamma_n(e)\,\boxwedge\, \gamma_n(e)"]
        		\arrow[d, "\mu"]  
				& 
			G^{\hspace{1pt}\boxwedge n} 
				\boxwedge 
			G^{\hspace{1pt}\boxwedge n} 
				\arrow[d, "\mu^{\hspace{1pt}\boxwedge n}"] 
				\\
			G 
				\arrow[r, "\gamma_n(e)"] 
				& 
			G^{\boxwedge n}
	    \end{tikzcd}
	\end{center}
	which, on a $V \in \Dip$, evaluates to 
	\begin{equation}\label{eq:prod_coprod}
    	\begin{tikzcd}
			\bigvee_{U\subseteq V} 
				G(U)\wedge G(V \ominus U) 
				\arrow[r]
				\arrow[d]  
				&  
			(G^{\hspace{1pt} \boxwedge n} 
				\boxwedge 
			G^{\hspace{1pt} \boxwedge n})(V) 
				\arrow[d]  
				\\
			G(V) 
				\arrow[r] 
				& 
			\displaystyle
			\bigvee_{0=W_0 \subseteq \dots \subseteq W_n=V}
				G(W_1)\wedge \dots \wedge G(W_n\ominus W_{n-1})
		\end{tikzcd}   
	\end{equation}
	where 
	\(
		(G^{\hspace{1pt}\boxwedge n} 
			\boxwedge 
		G^{\hspace{1pt}\boxwedge n})(V)
	\) 
	can be expanded to 

	\begin{adjustbox}{width=\textwidth}
	\(
		\bigvee_{U \subseteq V}
		\left(
			\bigvee_{0=U_0 \subseteq \dots \subseteq U_n=U} 
				G(U_1)\wedge \dots \wedge G(U_n\ominus U_{n-1})
		\right) 
			\wedge  
		\left(
			\bigvee_{0=V_0\subseteq\dots\subseteq V_n=V\ominus U} 
				G(V_1)
					\wedge \dots \wedge 
				G(V_n\ominus V_{n-1})\right).
	\)
	\end{adjustbox}

	Let $\phi\colon[0,1] \to \Sub(U)$ and $\psi\colon[0,1] \to \Sub(V\ominus U)$ be poset maps that represent step functions in $G(U)$ and $G(V \ominus U)$ respectively, and form an element in the top left corner of the diagram. 
	Following the left arrow of the diagram, we add these point-wise to obtain 
	\[
		\phi\oplus \psi\colon[0,1]\to \Sub(V)
	\]
	and pre-compose with $e_i$ to get 
	\[
		(\phi\oplus \psi)\circ e_i\colon [0,1]
			\to 
		\Sub(V).
	\]
	Tracing this along the bottom arrow of the diagram, we get
	\[
		(\phi\oplus \psi)_i\colon[0,1]
			\to 
		\Sub(W_i\ominus W_{i-1}), 
			\quad 
		t \mapsto (\phi\oplus \psi)(e_i(t)) \ominus W_{i-1}
	\]
	in $G(W_i \ominus W_{i-1})$,  where $W_i = (\phi \oplus \psi)(b_i) = \phi(b_i)\oplus \psi(b_i)$ for $1 \leq i \leq n-1$ and $W_0 = 0, W_n = V$. 
	For the alternative route, we first pre-compose with $e_i$ to get
	\(
		\phi \circ e_i\colon[0,1] \to \Sub(U) 
	\)
	and $\psi \circ e_i\colon[0,1] \to \Sub(V \ominus U),$ from which we obtain
	\[
		\phi_i\colon[0,1]\to \Sub(U_i\ominus U_{i-1}), \quad t \mapsto \phi(e_i(t)) \ominus U_{i-1}
	\]
	and
	\[
		\psi_i\colon[0,1]\to \Sub(V_i\ominus V_{i-1}), \quad t\mapsto \psi(e_i(t))\ominus V_{i-1}
	\]
	where $U_i = \phi(b_i)$ and $V_i = \psi(b_i)$. 
	Adding $\phi_i$ and $\psi_i$ point-wise then gives a map $\phi_i\oplus \psi_i$ that takes values in subspaces of 
	\begin{multline*}
		(\phi(b_i) \ominus \phi(b_{i-1})) 
			\oplus 
		(\psi(b_i) \ominus \psi(b_{i-1}))
			= 
		(\phi(b_i) \oplus \psi(b_i)) \ominus (\phi(b_{i-1})
			\oplus 
		\psi(b_i)) 
			= 
		W_i \ominus W_{i-1}   
	\end{multline*}
	and sends $t\in [0,1]$ to 
	\[
		(\phi(e_i(t)) \ominus U_{i-1}) 
			\oplus 
		(\psi(e_i(t)) \ominus  V_{i-1}) 
			= 
		(\phi(e_i(t)) \oplus \psi(e_i(t))) \ominus W_{i-1} 
			= 
		(\phi\oplus \psi)(e_i(t))\ominus W_{i-1}.
	\]
	This shows that $(\phi\oplus \psi)_i = \phi_i\oplus \psi_i$ for $1 \leq i \leq n$, and in particular that these maps represent the same step function in $G(W_i \ominus W_{i-1})$.
	Hence \cref{eq:prod_coprod} commutes, as desired.
\end{proof}

\begin{corollary}\label{cor:sah_functor_bialg}
	The spectral Sah functor $\funSah$ is a $(D_1)_+$-coalgebra in $\CMon(\textup{Fun}(\Dip,\Sp^O))$.
\end{corollary}

%%% COMPARISON TO CAMPBELL-ZAKHAREVICH 
\subsection{Comparison to Campbell-Zakharevich's derived Dehn invariant} 
\label{subsec:CZ_comparison}

In \cite[Definition 2.12]{cz-hilbert}, Campbell and Zakharevich define a topological version of the Dehn invariant that reduces to the classical one on homology groups. 
Their ``derived Dehn invariant'' is given as a map of pointed simplicial sets
\[
	D_i \colon F_\bullet^V 
		\to 
	\bigvee_{\substack{U \subseteq V \\ \dim(U) = i}}
		F_\bullet^U \,\tilde{\star}\, F_\bullet^{V\ominus U}
\]
for any inner product space $V$ and any $0<i< \dim(V)$, where $\tilde{\star}$ is the \emph{reduced join} of simplicial sets (recalled below). 
Here $U$ ranges over all $i$-dimensional subspaces of $V$, and $F_\bullet^V$ is a pointed simplicial set called the \textit{RT-building} associated to $V$ (\cite[Definition 1.8]{cz-hilbert}).
It is defined as the quotient of nerves
\begin{equation}\label{fv_defn}
	F^V_\bullet 
		= 
	\frac{N_\bullet(0 \subsetneq U \subseteq V)}%
		 {N_\bullet(0 \subsetneq U \subsetneq V)}
\end{equation}
and so it has an $n$-simplex for each chain of non-zero linear subspaces $0 \subsetneq U_0\subseteq \dots \subseteq U_n = V$, and those faces in which $U_n \neq V$ are collapsed to the basepoint. 
Upon realization, this gives the unreduced suspension $ST(V)$. 

Note that we have changed the conventions slightly from \cite{cz-hilbert}, letting $i$ denote the dimension of the vector space $U$, rather than the dimension of its underlying geometry $S(U)$, and denoting by $F_\bullet^V$ the simplicial set that is called $F_\bullet^{S(V)}$ in \cite{cz-hilbert}. 
As a result, what we call $D_i$ in the above display is in fact called $D_{i-1}$ in \cite{cz-hilbert}.

Recall that the simplicial join and its reduced version have $n$-simplicies given by  
\[
	(F_\bullet^U \star F_\bullet^{V\ominus U})_n 
		= 
	\coprod_{j=-1}^{n} F_j^U \times F_{n-j-1}^{V\ominus U}, 
		\qquad 
	(F_\bullet^U \,\tilde{\star}\, F_\bullet^{V\ominus U})_n 
		= 
	\bigvee_{j=0}^{n-1} F_j^U \wedge F_{n-j-1}^{V\ominus U}.
\]
On realization, this gives the usual topological join $\ST(U) \star \ST(V \ominus U)$ and its quotient by the subspace that is the join of each space with the basepoint of the other:
\[ 
	|F_\bullet^U \,\tilde{\star}\, F_\bullet^{V\ominus U}| 
		\cong 
	\frac{\ST(U) \star \ST(V \ominus U)}%
		 {(\ST(U) \star \{*\}) \cup (\{*\} \star \ST(V \ominus U))} 
		\cong 
	\ST(U) \sma S^1 \sma \ST(V \ominus U), 
\]
where $S^1$ is the interval $I$ modulo its endpoints. 
Geometrically, this space has an $(n-1)$-simplex for every $(j-1)$-simplex in $F_j^U$ and $(n-j-1)$-simplex in $F_{n-j-1}^{V \ominus U}$.
A point in this simplex with barycentric coordinates $(t_1,\ldots,t_n)$ corresponds to the points in the three spaces $\ST(U)$, $I$, and $\ST(V \ominus U)$ with coordinates
\[ 
	\left( 
		\frac{t_1}{t_1 + \cdots + t_j}, 
		\cdots, 
		\frac{t_j}{t_1 + \cdots + t_j}
	\right), 
		(t_1 + \cdots + t_j), 
	\left( 
		\frac{t_{j+1}}{t_{j+1} + \cdots + t_n}, 
		\cdots, 
		\frac{t_n}{t_{j+1} + \cdots + t_n}
	\right). 
\]
Along the most obvious homeomorphism between $S^1 \sma \ST(V \ominus U)$ and $\Sigma\!\ST(V \ominus U)$, this gives the two points in $\ST(U)$ and $\Sigma\!\ST(V \ominus U)$ with coordinates
\[ 
	\left( 
		\frac{t_1}{t_1 + \cdots + t_j}, 
		\cdots, 
		\frac{t_j}{t_1 + \cdots + t_j}
	\right), 
	\Big((t_1 + \cdots + t_j), t_{j+1}, \cdots, t_n\Big). \]
If we suspend one more time, we get a homeomorphism
\[ 
	\Sigma|F_\bullet^U \,\tilde{\star}\, F_\bullet^{V\ominus U}| \cong \Sigma\!\ST(U) \sma \Sigma\!\ST(V \ominus U). 
\]
Let's describe this more explicitly. 
An $n$-simplex on the left corresponds to an $(n-1)$-simplex in $F_\bullet^U \,\tilde{\star}\, F_\bullet^{V\ominus U}$, which comes from a $(j-1)$-simplex in $F_\bullet^U$ and an $(n-j-1)$-simplex in $F_\bullet^{V\ominus U}$. 
In other words, a flag of $j$ subspaces $0 \subsetneq U_1 \subsetneq \dots \subsetneq U_j = U$, which we can extend to a flag of $(j+1)$ subspaces to represent the cone point of $\Sigma\!\ST(U)$
\[ 
	0 = U_0 \subsetneq U_1 \subsetneq \dots \subsetneq U_j = U, \]
and a flag of $(n-j)$ subspaces $0 \subsetneq V_1 \subsetneq \dots \subsetneq V_{n-j} = V \ominus U$, which we can extend to a flag of $(n-j+1)$ subspaces to represent the cone point of $\Sigma\!\ST(V \ominus U)$
\[ 
	0 = V_0 \subsetneq V_1 \subsetneq 
		\dots 
	\subsetneq V_{n-j} = V \ominus U. 
\]
A point in this simplex with barycentric coordinates $(t_0,\ldots,t_n)$ is sent to the point in the corresponding product $\Delta^j \times \Delta^{n-j}$ with coordinates
\begin{equation}\label{homeo1}
	\left( 
		t_0, t_1\frac{t_1 + \cdots + t_n}{t_1 + \cdots + t_j},
		\cdots, 
		t_j\frac{t_1 + \cdots + t_n}{t_1 + \cdots + t_j}
	\right), 
	\left(
		\frac{t_1 + \cdots + t_j}{t_1 + \cdots + t_n},
		\frac{t_{j+1}}{t_1 + \cdots + t_n},
		\cdots,
		\frac{t_n}{t_1 + \cdots + t_n}
	\right).
\end{equation}
Everything goes to the basepoint if any of the three quantities $t_0$, $t_1 + \cdots + t_j$, or $t_{j+1} + \cdots + t_n$ is equal to zero, so the formula only has to make sense in the case that all these quantities are positive.

For the purpose of our proof, we will pick a different identification of 
\(
	\Sigma|F_\bullet^U \,\tilde{\star}\, F_\bullet^{V\ominus U}|
\)
with $ \Sigma\!\ST(U) \sma \Sigma\!\ST(V \ominus U)$, where the point with barycentric coordinates $(t_0,\ldots,t_n)$ is sent to the point in the corresponding product $\Delta^j \times \Delta^{n-j}$ with coordinates
\begin{equation}\label{homeo2}
	\left(
		\left( 
			\frac{t_0}{s_a},
			\dots, 
			\frac{t_{j-1}}{s_a}, 
			\frac{\frac{1}{2}t_j}{s_a} 
		\right),
		\left( 
			\frac{\frac{1}{2}t_j}{s_b},
			\frac{t_{j+1}}{s_b},
			\dots, 
			\frac{t_{n}}{s_b}  
		\right) 
	\right),
\end{equation}
where $s_a = \Sigma_{i=0}^{j-1} t_i + \frac{1}{2}t_j$ and $s_b = \frac{1}{2}t_j + \Sigma_{i=j+1}^{n} t_i$.

\begin{lemma}\label{better_homeo}
    The formula \cref{homeo2} defines a continuous map 
    \[
    	\Sigma|F_\bullet^U \,\tilde{\star}\, F_\bullet^{V\ominus U}| 
			\to 
		\Sigma\!\ST(U) \sma \Sigma\!\ST(V \ominus U)
	\] 
	that is homotopic to the one defined by \cref{homeo1}, and therefore a weak homotopy equivalence.
\end{lemma}

\begin{proof}
	The rule gives a bijection between the $n$-simplices on the left and the pairs of a $j$-simplex and an $(n-j)$-simplex on the right. 
	The above formula is continuous on the $n$-simplex, so we check it respects the face maps. 
	The first face of the $n$-simplex in which $t_0 = 0$ is glued to the basepoint by the definition of the suspension $\Sigma|\dots|$, and on the right it is also the basepoint by the definition of the suspension for $\Sigma\!\ST(U)$. 
	The first through $(j-1)$st faces agree along this map since the formula is unchanged and both sides delete the corresponding subspace $U_i$. 
	The $j$th face agrees because on both sides it is the basepoint; on the left this is because it makes $U_j \neq U$, while on the right we pass to the last face of $\Delta^j$ and the first face of $\Delta^{n-j}$, both of which go to the basepoint. 
	(One of them would have been enough.) 
	For the $(j+1)$st through $(n-1)$st faces we again get the same formula on both sides, and for the $n$th face we again get the basepoint on both sides.

	We therefore have a continuous map 
	\(
		\Sigma|F_\bullet^U \,\tilde{\star}\, F_\bullet^{V\ominus U}| 
			\to 
		\Sigma\!\ST(U) \sma \Sigma\!\ST(V \ominus U).
	\) 
	In fact, by the same argument we can give a continuous homotopy between this map and the map defined by \cref{homeo1}, because for each $\lambda \in [0,1]$ the weighted sum of $\lambda$ times \cref{homeo1} and $(1-\lambda)$ times \cref{homeo2} gives a formula that respects the face maps in the same way.
\end{proof}

Now we are ready to define Campbell and Zakharevich's topological Dehn invariant $D_i$ for $0<i<\dim(V)$ from \cite[Definition 2.9]{cz-hilbert}.
Let $0 \subsetneq U_0\subseteq \dots \subseteq U_n  = V$ be an $n$-simplex in $F^V_\bullet$. If none of the $U_j$ have dimension $i$, then $D_i$ sends such a simplex to the basepoint.
If there is a maximal $j$ such that $\dim U_j = i$, then the simplex is sent to the $n$-simplex in $F_\bullet^U \,\tilde{\star}\, F_\bullet^{V\ominus U}$ defined by the $j$-simplex in $F^{U_j}_\bullet$ given by the $(j+1)$ subspaces
\[ 
	U_0\subseteq \dots \subseteq U_j  = U_j 
\]
and the $(n-j-1)$-simplex in $F^{V\ominus U_j}_\bullet$ given by the $(n-j)$ subspaces
\[ 
	(U_{j+1}\ominus U_j)\subseteq\dots\subseteq(U_n\ominus U_j) 
		= 
	(V \ominus U_j). 
\]
After geometric realization, suspension, and the weak homotopy equivalence from \cref{better_homeo}, this derived Dehn invariant gives a map
% $$D_i:\ST(V) \to \bigvee_{\substack{U \subseteq X \\ \textup{dim}\ U = i}} \ST(U) \wedge \Sigma\!\ST(V\ominus U).$$
% Its suspension
\[
	\Sigma D_i\colon\Sigma\!\ST(V) 
		\to 
	\bigvee_{\substack{U \subseteq V \\ \textup{dim}\ U = i}}
		\Sigma\!\ST(U) \wedge \Sigma\!\ST(V\ominus U).
\]
It takes the $n$-simplex for the flag $0 = U_0 \subseteq \dots \subseteq U_n = V$ to the $j$-simplex in $|0 \subseteq U \subseteq U_j|$ given by the $(j+1)$ subspaces
\[ 
	0 = U_0\subseteq \dots \subseteq U_j = U_j, 
\]
and the $(n-j)$-simplex in $|0 \subseteq U \subseteq V \ominus U_j|$ given by the $(n-j+1)$ subspaces
\[ 
	0 = (U_j\ominus U_j) 
		\subseteq (U_{j+1}\ominus U_j) 
		\subseteq \dots 
		\subseteq (U_n\ominus U_j) 
		= (V \ominus U_j), 
\]
where again $j$ is the maximal index such that $\dim U_j = i$.
On the simplex itself, the map $\Delta^n \to \Delta^j \times \Delta^{n-j}$ is given by the formula \cref{homeo2}.

One might wonder if the suspension of the Campbell-Zakharevich derived Dehn invariant is anything like our Dehn invariant given by the $(D_1)_+$-action defined in \cref{const:E_1_coalg}. Recall that part of the $(D_1)_+$-coalgebra structure is given by a map 
\[
	\theta_2^V\colon\Sigma\!\ST(V) \times D_1(2) 
		\to 
	\bigvee_{U \subseteq V} 
		\Sigma\!\ST(U) \wedge \Sigma\!\ST(V \ominus U),
\]
where for $(\phi,e)\in \Sigma\!\ST(V) \wedge D_1(2),$ the element $e\in D_1(2)$ determines in which summand $\Sigma\!\ST(U) \wedge \Sigma\!\ST(V \ominus U)$ the image lies, and in particular, what the dimension of $U$ is.
On the other hand, for the suspended derived Dehn invariant $\Sigma D_i$ above, the dimension of $U$ is fixed. 
However, from $\theta^V_n$ we can obtain ``fixed-dimension Dehn invariants'' as follows.

\begin{definition}
    Let $\Sigma\!\ST(V)_i$ be the subspace of $\Sigma\!\ST(V) \sma D_1(2)_+$ consisting of those pairs $(\phi,e)$ such that $e = (e_1,e_2)$ is the element of $D_1(2)$ with $e_1 \colon I \to I$ the rectilinear embedding with image $[0,t]$ and $e_2 \colon I \to I$ the rectilinear embedding with image $[t,1]$. 
    Here $t$ is defined in terms of $\phi$ as follows: when we represent $\phi$ by a step function with images $0 = U_0 \subsetneq \dots \subsetneq U_n = V$, we let $t$ be the midpoint of the interval for the unique subspace $U_j$ whose dimension is equal to $i$. 
    If no such subspace appears, we instead let $t$ be the coordinate of the cut point where $\phi$ jumps from a subspace of dimension $< i$ to a subspace of dimension $> i$.
\end{definition}

Note that the forgetful map $\Sigma\!\ST(V)_i \to \Sigma\!\ST(V)$ is a homeomorphism: it has an inverse defined by choosing the unique set of embeddings $e$ that satisfy the above condition, and the formula defining these embeddings is continuous on each $n$-simplex in $\Sigma\!\ST(V)$.
% The homeomorphism is  containing only $k$-simplices associated to flags $0=U_0 \subseteq \dots \subseteq U_k = V$ where one of the $U_j$ has dimension $i$.
The restriction of $\theta_2^V$ to $\Sigma\!\ST(V)_i$ gives a map
\[
	(\theta_2^V)_i\colon\Sigma\!\ST(V)_i
		\to 
	\bigvee_{\substack{U \subseteq V\\\textup{dim}\ U = i}}
		\Sigma\!\ST(U) \wedge \Sigma\!\ST(V \ominus U).
\]
% Note that by definition, $\Sigma D_i$ is trivial outside of $\Sigma\!\ST(V)_i$. 

\begin{proposition}
    The following diagram commutes.
    \begin{center}
        \begin{tikzcd}
         	\Sigma\!\ST(V)_i 
             	\arrow[d, "\cong"] 
				\arrow[r, "{(\theta^V_2)_i}"] 
				& 
			\bigvee_{\dim U = i} \Sigma\!\ST(U) \wedge \Sigma\!\ST(V \ominus U)  
				\\
            \Sigma|F_\bullet^V| 
            	\arrow[r, "\Sigma D_i"] 
				& 
			\bigvee_{\dim U = i} \Sigma| F_\bullet^U \,\tilde{\star}\, F_\bullet^{V\ominus U}| 
				\arrow[u, "\simeq"]
        \end{tikzcd}
    \end{center}
\end{proposition}

\begin{proof}
    Here the left-hand map is the canonical identfication and the right-hand map is the weak homotopy equivalence of \cref{better_homeo}.
    To see that the diagram commutes, we only have to check what happens at a single point in $\Sigma\!\ST(V)$, which lies in the interior of an $n$-simplex indexed by a flag 
    \[
    	0 = U_0 \subsetneq \dots \subsetneq U_{n} = V
	\]
    with barycentric coordinates $(t_0,\dots,t_n)$ such that all $t_j > 0$. 
    We have spent most of this subsection describing what happens when we take the bottom route through the diagram: we go to the basepoint if none of the $U_j$ has dimension $i$, and otherwise if $U_j$ is the unique subspace of dimension $i$, we go to the point in $\Sigma\!\ST(U) \wedge \Sigma\!\ST(V \ominus U)$ lying in the product of simplices $\Delta^j \times \Delta^{n-j}$ corresponding to the two flags of subspaces 
    \[ 
    	0 = U_0\subsetneq \dots \subsetneq U_j = U_j, 
	\]
    \[ 
    	0 = (U_j\ominus U_j) 
			\subsetneq(U_{j+1}\ominus U_j) 
			\subsetneq \dots 
			\subsetneq (U_n\ominus U_j) = (V \ominus U_j), 
	\]
    and with barycentric coordinates 
    \[ 
    	\left(
			\left( 
				\frac{t_0}{s_a},
				\dots, 
				\frac{t_{j-1}}{s_a}, 
				\frac{\frac{1}{2}t_j}{s_a} 
			\right), 
			\left( 
				\frac{\frac{1}{2}t_j}{s_b}, 
				\frac{t_{j+1}}{s_b},
				\dots, 
				\frac{t_{n}}{s_b}   
			\right) 
		\right), 
	\]
    where $s_a = \Sigma_{i=0}^{j-1} t_i + \frac{1}{2}t_j$ and $s_b = \frac{1}{2}t_j + \Sigma_{i=j+1}^{n} t_i$.

    On the other hand, for the corresponding point $(\phi,e)$ in $ \Sigma\!\ST(V)_i$, $\phi$ is the step function taking the values $U_0,\dots, U_n$ with step lengths $t_0,\dots, t_n$. If there is no subspace here of dimension $i$, then the endpoint of the embedding $e_1$ lines up with a cut point of $\phi$, and therefore $(\theta_2^V)_i$ takes this to the basepoint. Otherwise
    % in the open $(n+1)$-simplex $$0 \subseteq U_0 \subseteq \dots \subseteq U_{n} = V$$  such that $\dim U_j = i$. Recall that we can see $\phi$ as a point in the topological simplex with non-zero coordinates $t_0,\dots, t_n$ adding up to 1, or equivalently, as a step function with values $0,U_0,\dots, U_n$ and step lengths $t_0,\dots, t_{n+1}.$ Then
    $(\theta_2^V)_i$ acts by precomposing the step function $\phi$ with $e_1$ and $e_2$. In terms of simplicial coordinates, this is exactly the continuous map $\Delta^n \to \Delta^{j+1} \times \Delta^{n-j}$ above, since the images of $e_1$ and $e_2$ cover all of $I$ and cut it in the middle between $t_0+\dots + t_{j}$ and $t_0+\dots +t_{j+1}$, i.e., half-way through the step where $\phi$ has value $U_j$. 
    This shows that the diagram commutes, as desired.
\end{proof}

%%%%%%%%%%%%%%%%%%%%%%%%%%%%%%%%%%%%%%%%%%%%%%%
%%% APARTMENTS AND THE SOLOMON-TITS THEOREM %%%
%%%%%%%%%%%%%%%%%%%%%%%%%%%%%%%%%%%%%%%%%%%%%%%
\section{Apartments and the Solomon-Tits theorem}\label{sec:apartments}

One of our goals is to prove that $\pi_0$ of the spectral Sah algebra is the classical Sah algebra: a direct sum of $O(n)$-coinvariants of the Steinberg modules $\St(\R^n)$.
In this section, we describe several equivalent ways to understand these Steinberg modules
\[ 
	\St(V) = H_{\dim(V)}(\Sigma\!\ST(V)), 
\]
along with the action of $O(V)$. 
Our goal is to describe them with enough data to make them into a diagram of abelian groups over the category $\Dip$.

These descriptions include a classical description by Lee and Szczarba and a more geometric description that appears in Sah's work. 
We clarify the relationship between these different models and use the concept of an ``apartment'' from Steinberg module theory to relate these models back to the Steinberg module. 
In the next sections we will describe the Hopf algebra structure on these groups and then use it to confirm that our spectral Sah algebra is indeed a Hopf algebra, and not just a bialgebra, in spectra.

%%% POLYTOPE GROUP AND LEE-SZCZARBA GROUp
\subsection{The polytope group and the Lee-Szczarba group}\label{subseq:polytope_group}

We begin by introducing two groups, presented in different ways and with different conventions, that are nonetheless isomorphic to each other, and also to the Steinberg module 
\(
	\St(V) = H_n(\Sigma\!\ST(V),*).
\) 
Working with these groups gives us a more combinatorial description of the spectral Sah algebra, and we will need this description to verify the Hopf algebra structure.

\begin{definition}\label{pt}
    For any real inner product space $V$ of finite dimension $n$, the \emph{(spherical) polytope group} $\Pt(V)$ is the group with one generator for each nondegenerate geodesic simplex $[v_1,\ldots,v_n]$ in the unit sphere $S(V)$, modulo the relation that when one simplex $P$ is covered by finitely many simplices $P_i$ with disjoint interiors, we have $[P] = \sum_i [P_i]$.
	When $n = 0$, the polytope group has a single generator, the empty simplex.
   
	We also allow the vectors $v_i$ to have non-unit length, in which case $[v_1,\ldots,v_n]$ refers to the simplex we obtain after rescaling the vectors to have unit length.
	Equivalently, instead of thinking of $[v_1,\ldots,v_n]$ as a simplex in $S(V)$, we can think of it as a conical region in $V$ obtained by all nonnegative linear combinations of the vectors,
	\[ 
		\left\{ \ 
			\sum_{i=1}^n a_i v_i \ \Big\vert \ a_i \geq 0 \ 
		\right\}. 
	\]
    
	For each codimension-one subspace $U \subseteq V$ with orthogonal complement spanned by $v \in S(V)$, we define a suspension map $\Pt(U) \to \Pt(V)$ by taking each polytope in $S(U)$ to its join with the vectors $v$ and $-v$. 
	In other words, it takes each simplex or conical region $[v_1,\ldots,v_{n-1}]$ to the sum
	\begin{equation}\label{pt_relation}
		[v_1,\ldots,v_{n-1},v] + [v_1,\ldots,v_{n-1},-v].
	\end{equation}
	The \emph{reduced polytope group} $\tPt(V) = \Pt(V)/\Sigma$ is the quotient of $\Pt(V)$ by the image of all of these suspension maps.
\end{definition}

Note that the order of the points $[v_1,\ldots,v_n]$ does not affect the resulting element of $\tPt(V)$, so rearranging them does not introduce a sign. 
On the other hand, by the relation \cref{pt_relation} in $\tPt(V)$, a sign is introduced each time we negate a vector:
\begin{align*}
    [v_2,v_1,v_3,\ldots,v_n] &= [v_1,v_2,v_3,\ldots,v_n], \\
    [-v_1,v_2,\ldots,v_n] &= -[v_1,v_2,\ldots,v_n].
\end{align*}

We let the orthogonal group $O(V)$ act on $\Pt(V)$ by acting on the vertices of each simplex, and again we do not introduce a sign, even when the matrix reverses orientation. 
This makes the groups $\Pt(-)$ and $\tPt(-)$ functors $\Dip\to\Ab$. 
We call the latter functor the \emph{polytope functor}.

\begin{definition}\label{ls}
    The \emph{Lee-Szczarba group} $\Ls(V)$ is the free abelian group on the $n$-tuples of vectors $(v_1,\ldots,v_n)$, sometimes called \emph{basic sharblies} \cite{AshMillerPatzt}, modulo the usual simplicial relation for each $(n+1)$-tuple
    \[ 
    	\sum_{i=1}^{n+1} 
			(-1)^i (v_1,\ldots,\hat{v_i},\ldots,v_{n+1}) 
		= 0, 
	\]
    and also modulo all tuples with all vectors lying in a proper subspace of $V$. 
    When $n = 0$, we do not impose the relations, so that $\Ls(0) = \Z$ with the empty tuple $()$ as the generator.
\end{definition}

Note that in this definition the generators are again given by $n$-tuples of linearly independent vectors $(v_1,\ldots,v_n)$ as in $\Pt(V)$, but now the order matters. 
If we swap any two of the elements it introduces a sign.
Similarly, it follows from the relations that if we replace $v_i$ with $\lambda v_i$ for $\lambda \neq 0$, it does not change the resulting element of $\Ls(V)$. 
So in a sense this follows the opposite sign rules as $\tPt(V)$:
\begin{align*}
    (v_2,v_1,v_3,\ldots,v_n) &= -(v_1,v_2,v_3,\ldots,v_n), \\
    (-v_1,v_2,\ldots,v_n) &= (v_1,v_2,\ldots,v_n).
\end{align*}
Without loss of generality the vectors $v_i$ are on the unit sphere $S(V)$, and we can negate any one of them without changing the element of $\Ls(V)$.

We let $O(V)$ act on $\Ls(V)$ by applying the matrix to each of the vectors $v_i$. 
We let $\Ls(V) \otimes \det$ denote $\Ls(V)$ with the same $O(V)$-action except that there is an additional factor of $-1$ each time the matrix reverses orientation. 
This makes the groups $\Ls(-) \otimes \det$ into a functor $\Dip\to\Ab$.

We call 
\[
	\begin{tikzcd}[cramped,row sep = 0ex,
			/tikz/column 1/.append style={anchor=base east},
    		/tikz/column 2/.append style={anchor=base west}]
		\Ls \otimes \det \colon \Dip \ar[r] & \Ab\\
		V \ar[r, mapsto] & \Ls(V) \otimes \det
	\end{tikzcd}
\]
the \emph{Lee-Szczarba functor}.

\begin{proposition}\label{pt_is_ls}
    Each orientation of $V$ gives an $O(V)$-equivariant isomorphism
    \[ 
    	\Ls(V) \otimes \det \cong \tPt(V). 
	\]
    Picking an orientation for each $V$ therefore gives an isomorphism as functors $\Dip\to\Ab$.
\end{proposition}

\begin{proof}
    We send $(v_1,\ldots,v_n)$ to $\pm [v_1,\ldots,v_n]$, the sign depending on whether the orientation on $V$ given by the basis $(v_1,\ldots,v_n)$ agrees with the chosen one. 
    The fact that this is an isomorphism follows from \cite[Corollary 2.11]{dupont_book} or \cite{sah_79}; an explicit proof appears in \cite[Theorem 6.4]{KLMMS-2}.
\end{proof}

%%%% APARTMENTS 
\subsection{Apartments}\label{subsec:apartments}

Now we relate the above two groups $\tPt(V)$ and $\Ls(V)$ to the Steinberg module $\St(V)$. 
Recall that $\St(V)$ is the homology of $\Sigma\!\ST(V)$. 
The idea is to give an explicit description of each sphere in $\Sigma\!\ST(V)$, one for each generator of $\Ls(V)$, which all together give an isomorphism $\Ls(V) \cong \St(V)$. %This amounts to defining in $\St(V)$ an ``apartment'' class for each generator of $\Ls(V)$. This can be done explicitly, as is usually the case in the literature, or homotopy-coherently, as in \cite{scissors_thom,KLMMS-2}.

Let $V$ be a finite-dimensional vector space of dimension $n$.
For the first two definitions, we require $n > 0$.

\begin{definition}\label{apartment_1}
    Let $n > 0$. 
    For each $n$-tuple of linearly independent vectors $(v_1,\ldots,v_n)$ from $V$, the corresponding \emph{apartment} is defined as the map
    \[ 
    	\partial \Delta^{n-1} \to \T(V) 
			= 
		|0 \subsetneq U \subsetneq V| 
	\]
    that takes the barycentric subdivision of $\partial \Delta^{n-1}$, forming a simplicial complex $\sd(\partial\Delta^{n-1})$, and then takes the map of simplicial complexes specified on vertices by
    \[ 
    	(\varnothing \subsetneq S \subsetneq \{v_1,\ldots,v_n\}) 
			\mapsto 
		\langle S \rangle. 
	\]
    This sends inclusions of subsets to inclusions of spaces, and always gives a proper nonzero subspace of $V$ when $\varnothing \neq S \neq \{v_1,\ldots,v_n\}$, and so defines the desired map $\sd(\partial\Delta^{n-1}) \to \T(V)$.
\end{definition}

\begin{figure}[ht]
	\centerline{\includegraphics[scale=1]{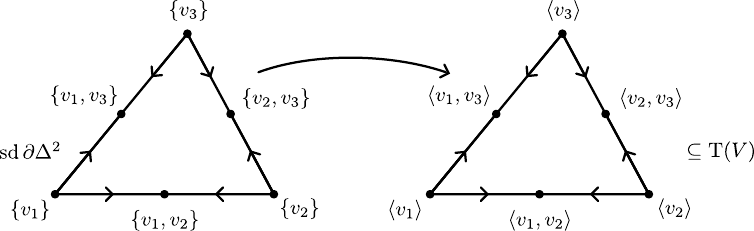}}
	\caption{The apartment map in dimension 3}
	\label{fig:apt}
\end{figure}

Note that the definition still makes sense if the vectors $(v_1,\ldots,v_n)$ are nonzero and linearly dependent. 
However, in that case the resulting map will be zero on homology.

Next we describe how this map suspends.
\begin{definition}\label{apartment_2}
	Let $n > 0$. 
    The \emph{suspended apartment} for $(v_1,\ldots,v_n)$ is the map of quotients
    \[ 
    	\frac{\Delta^{n-1}}{\partial \Delta^{n-1}} 
			\to 
		\frac{\CT(V)}{\T(V)} 
			\cong 
		\frac{|0 \subsetneq U \subseteq V|}%
			 {|0 \subsetneq U \subsetneq V|} 
	\]
    induced by the map of simplicial complexes 
    \(
    	\sd(\Delta^{n-1}) \to |0 \subsetneq U \subseteq V|
	\) 
	determined by the map of vertices
    \[ 
    	(\varnothing \subsetneq S \subseteq \{v_1,\ldots,v_n\}) 
			\mapsto 
		\langle S \rangle. 
	\]
\end{definition}

\begin{figure}[h]
	\centerline{\includegraphics[scale=1]{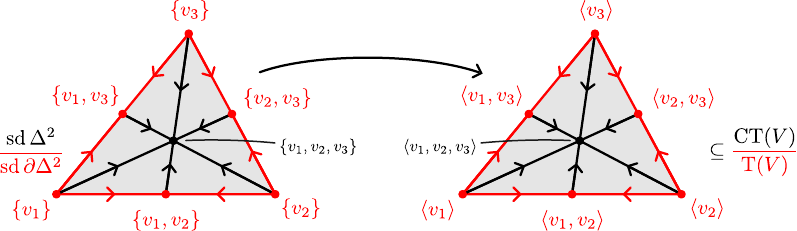}}
	\caption{The suspended apartment map in dimension 3}
	\label{fig:apt_susp}
\end{figure}

\begin{definition}\label{apartment_3}
    The \emph{doubly-suspended apartment} for $(v_1,\ldots,v_n)$ is the map of quotients
    \[ 
    	\frac{I^n}{\partial I^n} 
			\to 
		\frac{|0 \subseteq U \subseteq V|}
			 {|0 \subseteq U \subsetneq V| 
			 	\cup 
			  |0 \subsetneq U \subseteq V|} 
	\]
    induced by the map of posets $\underline{2}^n \to \{0 \subseteq U \subseteq V\}$ that sends
    \[ 
    	(\varnothing \subseteq S \subseteq \{v_1,\ldots,v_n\}) 
			\mapsto 
		\langle S \rangle. 
	\]
\end{definition}

\begin{figure}[ht]
	\centerline{\includegraphics[scale=1]{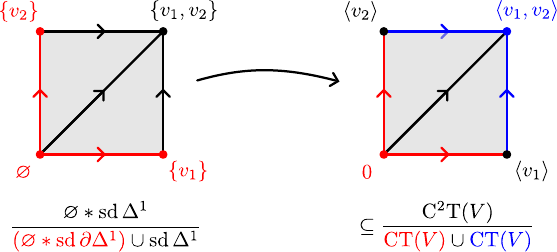}}
	\caption{The doubly-suspended apartment map in dimension 2}
	\label{fig:apt_2susp_1}
\end{figure}
\begin{figure}[ht]
	\centerline{\includegraphics[scale=1]{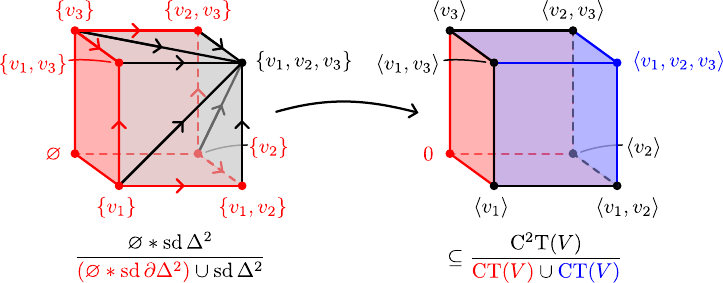}}
	\caption{The doubly-suspended apartment map in dimension 3}
	\label{fig:apt_2susp_2}
\end{figure}

Notice that the formula for the doubly-suspended apartment makes sense when $n = 0$ as well, where the above map of quotients becomes the identity map $S^0 \to S^0$.

\begin{lemma}
    These two types of apartment maps are identified with the one-fold and two-fold suspensions of the apartment map for $\T(V)$, respectively.
\end{lemma}

\begin{proof}
    The first one follows because both $\Delta^{n-1}$ and $\CT(V)$ are contractible, and the given map $\sd(\Delta^{n-1}) \to \CT(V)$ agrees on the boundary with the map $\sd(\partial \Delta^{n-1}) \to \T(V)$ of \cref{apartment_1}.

    For the second one, we write $I^n$ as a simplicial complex in the standard way, arising from the cube-shaped poset $\underline{2}^n$. 
    Note that this is isomorphic to the simplicial complex $0 * \sd(\Delta^{n-1})$, a model for $\Delta^n$ in which instead of taking the barycentric subdivision of all of $\Delta^{n}$, we only subdivide $\Delta^{n-1}$ and then join the result with 0.
    
    Taking the join of the previous map with $0$, which is represented by $S = \varnothing$, gives the above formula. Quotienting out by the subspaces described above gives the suspension of the map from \cref{apartment_2}.
\end{proof}

The following form of the Solomon-Tits theorem can be found in \cite[Corollary 2.11]{dupont_book} or \cite[Theorem A]{KLMMS-2}.

\begin{proposition}\label{solomon_tits_w_apts}
    If $n = \dim(V)$, then the space $\Sigma\!\ST(V)$ is homotopy equivalent to a wedge sum of $n$-spheres, and the apartment maps for each generator of $\Ls(V)$ induce an $O(V)$-equivariant isomorphism
    \[ 
    	[\apt]\colon \Ls(V) \cong H_n(\Sigma\!\ST(V),*) 
			\cong 
		\St(V). 
	\]
\end{proposition}

\begin{remark}\label{rmk:orientation}
    Once we start considering a map to homology, we should be clear about the choice of orientation for the simplicial complex $I^n \cong 0 * \sd(\Delta^{n-1})$. 
    It is easiest to think of ``orientation'' as a property that simplices have, and that is specified by the choice of ordering on the vertices, up to any even permutation. 
    Under this convention, we orient $\Delta^n$ by the ordering $0,1,2,\ldots,n$. 
    When we subdivide the face of $\Delta^n$ spanned by $\{1,\ldots,n\}$ and join to 0, we give each simplex 
    \(
    	\{0\} 
		\subseteq \{0,x_1\} 
		\subseteq \cdots 
		\subseteq \{0,x_1,\ldots,x_n\}
	\) 
	the orientation given by whichever permutation to $\{x_1,\ldots,x_n\}$ is required to rearrange it into $\{1,\ldots,n\}$. 
	This turns out to agree with the orientation on $I^n$ that results from ordering its coordinate directions as $\{1,\ldots,n\}$.
\end{remark}

\begin{remark}
    The isomorphism of \cref{solomon_tits_w_apts} does not depend on a choice of orientation of $V$. 
    Each basis $(v_1,\ldots,v_n)$ gives the cube $I^n$ an orientation that allows us to make a well-defined homology class.
\end{remark}

\begin{notation}
    We use $\apt(v_1,\ldots,v_n)$ to refer to the map of spaces $\frac{I^n}{\partial I^n} \to \Sigma\!\ST(V)$ of \cref{apartment_3}, and $[\apt(v_1,\ldots,v_n)]$ to refer to the corresponding class in $H_n(\Sigma\!\ST(V),*)$, using the orientation conventions discussed in \cref{rmk:orientation}. 
    As in \cref{solomon_tits_w_apts}, we also use $[\apt]$ to refer to the homomorphism $\Ls(V) \to H_n(\Sigma\!\ST(V),*)$ that takes each tuple to its apartment class.
\end{notation}

We can now give a precise relationship between the two groups $\Ls(V)$ and $\tPt(V)$ and the homology of $\Sigma\!\ST(V)$.
Consider the following diagram of abelian groups.
\begin{equation}\label{un-confuse}
	\begin{tikzcd}[arrows=rightarrow]
			&
		\Ls(V) \otimes \det 
			\ar{d}{[\apt]}[swap]{\cong} 
			\ar{r}{\textup{\cref{pt_is_ls}}}[swap]{\cong} 
			& 
		\tPt(V) 
			\ar{d}{[\apt]}[swap]{\cong} 
			\\
		\St(V) \otimes \det 
			\arrow[r, equals] 
			& 
		H_n(\Sigma\!\ST(V),*) \otimes \det 
			\ar{r}[swap]{\cong} 
			& 
		H_0(\Sigma^{-V}\Sigma\!\ST(V))
	\end{tikzcd}
\end{equation}
The vertical map on the left is the apartment map of \cref{apartment_3}, and is an isomorphism by \cref{solomon_tits_w_apts}. 
The vertical map on the right is defined in the same way, except that we think of the unit cube $I^n$ as embedded in $V$, and define a map $S^V \to \Sigma\!\ST(V)$ by collapsing $V$ onto this cube and then applying the apartment map. 
The horizontal map along the top is defined as in \cref{pt_is_ls} by sending $(v_1,\ldots,v_n)$ to $\pm[v_1,\ldots,v_n]$, and the bottom horizontal map is the suspension isomorphism on homology.
    
\begin{proposition}\label{clarifying_diagram}
    In the diagram \cref{un-confuse}, all maps are $O(V)$-equivariant, and the horizontal maps depend on a choice of orientation of $V$, while the vertical maps do not. 
    The diagram commutes provided we use the same orientation on $V$ throughout.
\end{proposition}

\begin{proof}
    We leave the $O(V)$-equivariance and commutativity of the diagram to the reader.
    We observe that the horizontal maps each change by a factor of $-1$ when the orientation of $V$ is reversed, while the vertical maps do not use the orientation of $V$ in their definition.
\end{proof}

%%%%%%%%%%%%%%%%%%%%%%%%%%%%%%%%%%%%%%%%%%%%%%
%%% HOPF ALGEBRA STRUCTURE ON THE HOMOLOGY %%%
%%%%%%%%%%%%%%%%%%%%%%%%%%%%%%%%%%%%%%%%%%%%%% 
\section{Hopf algebra structure on the homology of \texorpdfstring{$\Sigma\!\ST(V)$}{∑ST(V)}}\label{sec:Hopf_alg_diagram}

In this section we define a version of the classical Sah algebra in which we don't take orbits, and check that it is still a Hopf algebra. 
The most natural way to describe this is via the polytope functor $\widetilde \Pt$, which is a diagram of abelian groups on the category $\Dip$ that we studied in \cref{sec:apartments}.
We will see that it is a Hopf algebra with respect to the Day convolution product on such diagrams. 
We use this to recover the result that the Sah algebra (after taking orbits) is a Hopf algebra in abelian groups. 
We will also use this in the next section to verify that the spectral Sah functor $\funSah$ is a Hopf algebra in $\Fun(\Dip, \Sp^O)$.

%%% POLYTOPE HOPF ALGEBRA
\subsection{The polytope Hopf algebra}

We next recall Sah's original description of the Hopf algebra structure on the groups $\tPt(\R^n)_{O(n)}$,
% \[ \bigoplus_{n \geq 0} (\Pt(\R^n)/\Sigma)_{O(n)} \]
see \cite[Sections 6.1 and 6.3]{sah_79}. 
It arises from taking coinvariants of the following product, coproduct, and involution on the groups $\tPt(\R^n)$.

\begin{definition}\label{pt_hopf}
	The product
	\begin{align*}
    	\mu\colon \tPt(U) \otimes \tPt(V) &\to \tPt(U \oplus V)
	\end{align*}
	takes simplices $P \subseteq S(U)$ and $Q \subseteq S(V)$ to their join $P * Q \subseteq S(U) * S(V) \cong S(U \oplus V)$. 
	Alternatively, the product is given by the formula
	\[ 
		\mu \colon 
		[u_1,\ldots,u_m],[v_1,\ldots,v_n] 
			\mapsto 
		[u_1,\ldots,u_m,v_1,\ldots,v_n], 
	\]
	where each vector in $U$ or $V$ is regarded as a vector in $U \oplus V$ on the right-hand side. 
	The coproduct
	\[
    	\delta\colon \tPt(V) 
			\to 
		\bigoplus_{0 \subseteq U \subseteq V} 
			\tPt(U) \otimes \tPt(V \ominus U)
	\]
	takes each simplex $[v_1,\ldots,v_n]$ to the sum
	\[ 
		[v_1,\ldots,v_n] 
			\mapsto 
		\sum_{S\subseteq \{v_1,\dots, v_n\}} 
			[S] \otimes [\pr_{V \ominus \langle S\rangle} S^c]. 
	\]
	Here $\pr_{V \ominus \langle S\rangle} S^c$ means that we take the set $S^c = \{v_1,\ldots,v_n\} \setminus S$ and then project onto the orthogonal complement of the span of $S$.
	Geometrically, each term of this sum is one face of the simplex (whose span is some subspace $U \subseteq V$), and the link of that face as a subset of the unit sphere of the complement, $S(U^\perp)$.

	The involution or antipode is given by
	\[
		\begin{tikzcd}[cramped, row sep = 0ex]
			\alpha \colon \tPt(V) 
				\ar[r] 
				& 
			\tPt(V)
				\\
			\phantom{a \colon}{[v_1,\ldots,v_n]}
				\ar[r, mapsto]
				&
			{[v_1^\vee, \ldots, v_n^\vee]}
		\end{tikzcd}
	\]
	where 
	\(
		v_i^\vee 
		= \pr_{V \ominus \langle \{v_1,\ldots,\hat{v_i},\ldots,v_n\}\rangle} (-v_i).
	\)
	So $v_i^\vee$ is the projection of $-v_i$ onto the orthogonal complement of the span of the other vectors. 
	In other words, $v_i^\vee$ is the unique vector (up to scaling by $\lambda \in \R_{>0}$) such that
	\[ 
		v_i^\vee \cdot v_j = 0 \quad \textup{for all} \quad j \neq i, \quad \textup{and} \quad v_i^\vee \cdot v_i < 0. 
	\]
% that is orthogonal to $v_j$ for every $j \neq i$ and such that $v_i \cdot v_i^\vee < 0$.
	This condition is symmetric enough in $v_i$ and $v_i^\vee$ that we can observe right away that $(v_i^\vee)^\vee$ is a positive multiple of $v_i$ for every $i$. 
	This shows that $\alpha^2 = \id$.
\end{definition}

More generally, for any subset $S \subseteq \{v_1,\ldots,v_n\}$, we let $S^\vee$ denote the set of dual vectors
\[ 
	S^\vee = \{ v_i^\vee \mid i \in S \} 
		= \{ \pr_{V \ominus \langle \{v_1,\ldots,\hat{v_i},\ldots,v_n\}\rangle} (-v_i) \mid i \in S \}. 
\]
Note that this depends on the entire basis $\{v_1,\ldots,v_n\}$ and not just the subset $S$.

Before proving that this is a Hopf algebra, we need two results from \cite{sah_79}. We include the proofs since the book \cite{sah_79} is difficult to obtain and uses less modern language.

\begin{lemma}[{\cite[page 116]{sah_79}}]\label{orthogonal_lemma}
    We have an equality of linear transformations
    \[ 
    	\pr_{V \ominus \langle \{v_1,\ldots,\hat{v_j},\ldots,v_n\}\rangle} 
			= 
		\pr_{V \ominus \langle \pr_{V \ominus \langle S \rangle} (S^c \setminus \{v_j\}) \rangle} 
		\circ 
		\pr_{V \ominus \langle S \rangle} 
	\]
    for each value of $j \in S^c$.
\end{lemma}

\begin{proof}
	Choose an orthonormal basis $\{u_1,\dotsc,u_r\}$ of $\langle S \rangle$, then extend it to an orthonormal basis $\{u_1,\dotsc,u_{n-1}\}$ of $\langle v_1,\ldots,\hat{v_j},\ldots,v_n \rangle$, and finally to an orthonormal basis $\{u_1,\ldots,u_n\}$ of $V$. 
	Then
	\begin{align*}
    	\pr_{V \ominus \langle \{v_1,\ldots,\hat{v_j},\ldots,v_n\}\rangle} 
			&= 
		\pr_{\langle u_n \rangle}, 
			\\
    	\pr_{V \ominus \langle S \rangle} 
			&= 
		\pr_{\langle u_{r+1},\dotsc,u_n \rangle}, 
		\\
    	\langle 
			\pr_{V \ominus \langle S \rangle} 
				(S^c \setminus \{v_j\}) 
		\rangle 
			&= 
		\langle u_{r+1},\dotsc,u_{n-1} \rangle,
	\end{align*}
	and so 
	\[
	    \pr_{V \ominus \langle \pr_{V \ominus \langle S \rangle} (S^c \setminus \{v_j\}) \rangle} 
	    	= 
		\pr_{\langle u_1,\dotsc,u_r,u_n \rangle}.
	\]
	Therefore if we write each vector in $V$ as $\sum_{i=1}^n \lambda_i u_i \in V$, then
	\[
		\pr_{V \ominus \langle \{v_1,\ldots,\hat{v_j},\ldots,v_n\}\rangle} 
		\left(\sum_{i=1}^n \lambda_i u_i \right) 
		= \lambda_n u_n, 
	\]
	while
	\[
		\pr_{V \ominus \langle \pr_{V \ominus \langle S \rangle} (S^c \setminus \{v_j\}) \rangle} 
			\circ 
		\pr_{V \ominus \langle S \rangle}\left( \sum_{i=1}^n \lambda_i u_i \right) 
			= 
		\pr_{V \ominus \langle \pr_{V \ominus \langle S \rangle} (S^c \setminus \{v_j\}) \rangle}\left(\sum_{i= r+1}^n \lambda_i u_i\right) 
			= 
		\lambda_n u_n.
	\] 
These agree, as claimed.
\end{proof}

\begin{proposition}[{\cite[Chapter 6, Proposition 3.4]{sah_79}}]\label{sphere_decomposition}
    The unit sphere $S(V)$ decomposes into the union of the simplices $[S \amalg (S^c)^\vee]$ for $S \subseteq \{v_1,\ldots,v_n\}$, and therefore
    \[ 
    	\sum_{S \subseteq \{v_1,\ldots,v_n\}} [S \amalg (S^c)^\vee] 
		= 0 
		\quad \textup{ in the group } \tPt(V). 
	\]
\end{proposition}

\begin{figure}[h]
    \begin{minipage}[t]{0.5\textwidth}\vspace{0pt}
        \[
	\begin{tikzpicture}[scale=1.6]
		\draw[thick] (0,0) circle (1cm);
		\draw[thick,->] (0,0) -- (0:1.5) node[right]{$v_1$};
		\draw[thick,->] (0,0) -- (30:1.5) node[right]{$v_2$};
		\draw[thick,dashed,->] 
			(0,0) -- (120:1.5) node[above]{$v_1^\vee$};
		\draw[thick,dashed,->] 
			(0,0) -- (270:1.5) node[below]{$v_2^\vee$};
		\draw[ultra thick, OI2] 
			(1,0) arc (0:30:1) node[midway,right]{$[v_1,v_2]$} ;
		\draw[ultra thick, OI7] 
			(30:1) arc (30:120:1)
				node[midway,above]{$[v_2,v_1^\vee]$};
		\draw[ultra thick,OI3] 
			(120:1) arc (120:270:1)
				node[midway,left]{$[v_1^\vee,v_2^\vee]$};
		\draw[ultra thick,orange] 
			(270:1) arc (270:360:1)
				node[midway,anchor=north west]
					{$[v_1,v_2^\vee]$};
	\end{tikzpicture}
    \] 
    \end{minipage}%
    \begin{minipage}[t]{0.5\textwidth}\vspace{2em}
	    \[\includegraphics[scale=1]{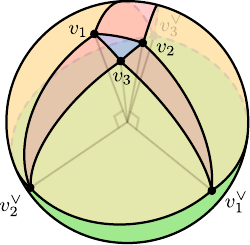}\]
    \end{minipage}%
    \caption{The cases of $n = 2$ and $n = 3$ of \cref{sphere_decomposition}.}
\end{figure}

\begin{proof}
    We first establish that for every subset $S \subseteq \{v_1,\ldots,v_n\}$, the set of vectors $S \amalg (S^c)^\vee$ is linearly independent in $V$. 
    Using \cref{orthogonal_lemma}, suffices to project $S^c$ to the orthogonal complement of $S$, and show that its dual inside the vector space $V \ominus \langle S \rangle$ is linearly independent. 
    But the dual of any basis (in the above sense) is always a basis again, for if it were not linearly independent, it would lie in a proper subspace $U$, and then the double dual would again lie in $U$, which is a contradiction.

    For the next part of the argument, we think of $[S \amalg (S^c)^\vee]$ as a conical region in the vector space $V$, and show that as $S$ varies these regions cover $V$. 
    The first of these regions we denote
    \[ 
    	C = [v_1,\ldots,v_n] 
		  = 
		\left\{\ 
			\sum_{i=1}^n a_i v_i \ \Big\vert \ a_i \geq 0 \ 
		\right\}. 
	\]
    % which corresponds to the case where $S$ is the entire basis $\{v_1,\ldots,v_n\}$.
    For any other vector $x \in V$, let $\bar x \in C$ be the closest point in $C$ to $x$, which exists because $C$ is a closed set. %(It suffices to consider those points in $C$ whose distance to the origin is no more than $|x|$.)
    We write
    \[ 
    	\bar x = \sum_{i=1}^n a_i v_i 
	\]
    and let $S$ be the set of indices for which $a_i > 0$. 
    Since $S \amalg (S^c)^\vee$ is a basis and the vectors in $(S^c)^\vee$ are all orthogonal to those in $S$, we can express the original vector $x$ as
    \begin{equation}\label{sphere_decomposition_key_step}
        x = \bar x + (x - \bar x) 
          = \sum_{v_i \in S} a_i v_i 
          + \sum_{v_i \in S^c} b_i v_i^\vee, 
        \qquad a_i > 0, \ \ b_i \in \R.
    \end{equation}
    For each $j \in S^c$, the line segment connecting $\bar x$ to $v_j$ lies in $C$. 
    Therefore the angle formed by $x$, $\bar x$, and $v_j$ must be greater than or equal to $\pi/2$, for otherwise moving $\bar x$ in the direction of $v_j$ would yield a point in $C$ that is closer to $x$. 
    Therefore
    $(x - \bar x) \cdot v_j \leq 0$, so
    \[ 
    	b_j(v_j^\vee \cdot v_j) 
			= \Big(
				\sum_{v_i \in S^c} b_i v_i^\vee
			  \Big) 
			  \cdot v_j 
			= (x - \bar x) \cdot v_j \leq 0. 
	\]
    % (Note that the $v_i$ are orthogonal to the $v_j^\vee$, though within each of these two sets the vectors do not have to be orthogonal.)
    Since by assumption $(v_j^\vee \cdot v_j) < 0$, we conclude that $b_j \geq 0$. 
    Therefore $x$ is a nonnegative linear combination of the vectors in $S \amalg (S^c)^\vee$, so it lies in the conical region $[S \amalg (S^c)^\vee]$.

    This shows that the regions $[S \amalg (S^c)^\vee]$ cover the vector space $V$, and it remains to show that their interiors are disjoint. 
    If $x$ is in the interior of $[S \amalg (S^c)^\vee]$, then it can be expressed as in \cref{sphere_decomposition_key_step}, with all $a_i > 0$ and all $b_i > 0$. 
    Since the vectors in $(S^c)^\vee$ are all orthogonal to those in $S$, the first sum $\sum_{i \in S} a_iv_i$ must be the closest point $\bar x$ in $C$, which is uniquely determined by $x$. 
    Therefore $S$ is determined to be the set of all nonzero coordinates for $\bar x$ in the basis $\{v_1,\ldots,v_n\}$, so there is only one possible choice for $S$, finishing the proof.
\end{proof}

\begin{proposition}\label{lem:polytope_functor_hopf}
    The above product, coproduct and involution make the polytope functor $\tPt(V) = \Pt(V)/\Sigma$ into a Hopf algebra in the symmetric monoidal category $\Fun(\Dip, \Ab)$.
\end{proposition}

\begin{proof}
    Let $\boxtimes$ denote the symmetric monoidal structure on $\Fun(\Dip,\Ab)$. To show that the product and coproduct give a bi-algebra structure, we need to check that the diagram
    \begin{equation}\label{eq:P_is_bialg}
	    \begin{tikzcd}
    			& 
			\tPt \boxtimes \tPt 
				\arrow[rr, "\mu"] 
				\arrow[dl, "\delta \boxtimes \delta"'] 
				& 
				& 
			\tPt
				\arrow[dr, "\delta"] 
				& 
				\\
	        \tPt\boxtimes \tPt \boxtimes \tPt \boxtimes \tPt
    	    	\arrow[rr, "1 \boxtimes \tau \boxtimes 1"] 
				& 
				&  
			\tPt\boxtimes \tPt \boxtimes \tPt \boxtimes \tPt 
				\arrow[rr, "\mu \boxtimes \mu"] 
				& 
				& 
			\tPt \boxtimes \tPt
	    \end{tikzcd}
	\end{equation} 
	commutes, where $\tau$ is the symmetry isomorphism. 
	Starting with $[v_1,\dots, v_n]$ in $\tPt(V)$ and $[u_1,\dots, u_m]$ in $\tPt(U)$, first applying the product and then the coproduct gives
	\[
		\sum_{R\subseteq\{v_1,\dots, v_n,u_1,\dots, u_m\}} 
			[R] 
			\otimes 
			[\pr_{(V\oplus U) \ominus \langle R\rangle} R^c].
	\]
	Here we denote by $[R]$ the element of $\tPt(V)$ given by the vectors in $R$, by $\langle R\rangle $ the linear subspace of $V\oplus U$ spanned by the vectors in $R$, and by $R^c$ the complement of $R$ in $\{v_1,\dots, v_n,u_1,\dots, u_m\}$.
 
	On the other hand, applying the coproduct and $\tau$ first gives the element 
	\[
		\sum_{S\subseteq \{v_1,\dots, v_n\}} 
			\sum_{T\subseteq\{u_1,\dots, u_m\}} 
				[S] \otimes [T] 
					\otimes 
				[\pr_{V \ominus \langle S\rangle} S^c ]
					\otimes 
				[\pr_{U\ominus \langle T\rangle} T^c].
	\]
	Now we observe that sets $R$ correspond one-to-one with pairs of a set $S$ and a set $T$, and the map $\mu\boxtimes \mu$ turns 
	\(
		[S] \otimes [T] 
		    \otimes [\pr_{V \ominus \langle S\rangle} S^c ]
			\otimes [\pr_{U\ominus \langle T\rangle} T^c]
	\) 
	into 
	\(
		[R] 
			\otimes 
		[\pr_{(V\oplus U) \ominus \langle R\rangle} R^c].
	\) 
	This shows that \cref{eq:P_is_bialg} commutes, as desired.

	In order to show the involution makes this bialgebra into a Hopf algebra, we need to verify that the diagram 
	\begin{center}
		\begin{tikzcd}
    			& 
			\tPt \boxtimes \tPt 
				\arrow{rr}{\alpha \otimes \id} 
				&  
				&
			\tPt \boxtimes \tPt 
				\arrow[dr, "\mu"]
				& 
				\\
		    \tPt 
		    	\arrow[rr, "\eta"] 
				\arrow[ru, "\delta"] 
				\arrow[rd, "\delta"'] 
				&  
				& 
			\mathbbm{1} 
				\arrow[rr, "\varepsilon"] 
				& 
				& 
			\tPt 
				\\
	    		& 
			\tPt \boxtimes \tPt 
				\arrow{rr}{\id \otimes \alpha} 
				&  
				&
			\tPt \boxtimes \tPt 
				\arrow[ur, "\mu"']
		\end{tikzcd}   
	\end{center}
	commutes, where $\mathbbm{1}$ is the diagram on $\Dip$ that such that $\mathbbm{1}(0)=\mathbb{Z}$ and $\mathbbm{1}(V) = 0$ for $V\neq 0$; the maps $\eta$ and $\varepsilon$ are induced by the identity on $\mathbb{Z}$. 
	We check this on a single element $[v_1,\ldots,v_n]$. 
	If $n = 0$ then all maps are the identity of $\Z$, so there is nothing to check. 
	If $n > 0$, then we compute
	\begin{align*}
        \mu \circ (\id \otimes \alpha) \circ \delta([v_1,\ldots,v_n])
        	&= 
		\mu\left( 
			\sum_{S \subseteq \{v_1,\ldots,v_n\}} 
				[S] 
				\otimes 
				\alpha[\pr_{V \ominus \langle S\rangle} S^c] 
		\right) 
			\\
        	&= 
		\mu\left( 
			\sum_{S \subseteq \{v_1,\ldots,v_n\}} 
				[S] \otimes [(S^c)^\vee] 
		\right) 
			\\
        	&= 
		\sum_{S \subseteq \{v_1,\ldots,v_n\}} 
			[S \amalg (S^c)^\vee] 
			\\
        	&= 
		0,
	\end{align*}
	where in the second line we use \cref{orthogonal_lemma} to conclude that 
	\(
		\alpha[\pr_{\langle S \rangle^\perp} S^c] 
		= [(S^c)^\vee],
	\) 
	as in the proof of \cite[Chapter 6, Proposition 3.7]{sah_79}, and in the last line we use \cref{sphere_decomposition}. 
	Going the other way,
	\begin{align*}
        \mu \circ (\alpha \otimes \id) \circ \delta([v_1,\ldots,v_n])
        	&= 
		\mu\left( 
			\sum_{S \subseteq \{v_1,\ldots,v_n\}} 
				\alpha[S] \otimes 
				[\pr_{V \ominus \langle S\rangle} S^c] 
		\right) 
			\\
        \alpha \circ \mu \circ (\alpha \otimes \id) \circ c([v_1,\ldots,v_n])
        	&= 
		\mu\left( 
			\sum_{S \subseteq \{v_1,\ldots,v_n\}} 
				[S] \otimes 
				a[\pr_{V \ominus \langle S\rangle} S^c] 
		\right) 
			\\
	        &= 
		\sum_{S \subseteq \{v_1,\ldots,v_n\}} 
			[S \amalg (S^c)^\vee] 
			\\
        	&= 
		0
	\end{align*}
	from which we deduce that $\mu \circ (\alpha \otimes \id) \circ \delta([v_1,\ldots,v_n]) = 0$ as well.
	The equality in the second line follows from the fact that $\alpha \circ \mu = \mu \circ (\alpha \otimes \alpha)$, which can be checked on generators $[v_1, ..., v_n], [u_1, ..., u_m]$. 
	The crucial property that enables this is that every $v_i$ is perpendicular to every $u_j$; note that this is also true for the vectors in the terms $[S]$ and $[\pr_{V \ominus \langle S\rangle} S^c]$ that appear in the second line. 
\end{proof}

Taking the colimit along $\Dip$, we recover Sah's original result:

\begin{corollary}[{\cite[Chapter 6, Proposition 3.7]{sah_79}}]\label{lem:polytope_colimit_hopf}
    The above product, coproduct and involution make the colimit of the polytope functor
    \[ 
    	\bigoplus_{n \geq 0} \tPt(\R^n)_{O(n)} 
	\]
    into a commutative Hopf algebra in the category of abelian groups.
\end{corollary}

\begin{proof}
    This follows from the fact that the colimit functor $\textup{Fun}(\Dip, \Ab) \to \Ab$ is strong symmetric monoidal, and therefore preserves commutative Hopf algebra objects.
\end{proof}

%%% LEE-SZCZARBA HOPF ALGEBRA
\subsection{The Lee-Szczarba Hopf algebra}

If we are willing to pick orientations for each vector space $V$, the above operations also make the groups $\Ls(V)$ into a Hopf algebra object in $\textup{Fun}(\Dip, \Ab)$. 
We spell this out since we will need the formulas in \cref{sec:hop_functor}.

Assume we have fixed an orientation for each finite-dimensional inner product space $V$. 
Along the isomorphism of \cref{pt_is_ls}, the product is the map $m \colon \Ls \boxtimes \Ls \to \Ls $ that is on $V$ in $\Dip$ given by the maps 
\[
	\begin{tikzcd}[row sep=0]
    	\Ls(U)\otimes \Ls(V\ominus U)
			\ar[r]
			& 
		\Ls(V) 
			\\
		{(u_1,\dots, u_n)\otimes (v_1,\dots, v_m)}
			\ar[r,mapsto]
			& 
		{(u_1,\dots, u_n,v_1,\dots, v_m),}
    \end{tikzcd}
\]
together with a sign if the chosen orientation of $V$ does not match the one on $U \oplus (V \ominus U)$. 
These assemble into a map 
\[
	\bigoplus_{U\subseteq V}\Ls(U)\otimes \Ls(V\ominus U) 
		\to \Ls(V)
\]
 
The coproduct  $\delta \colon \Ls\to \Ls \boxtimes \Ls$ is given on $V$ by
\[
	\begin{tikzcd}[row sep=0]
    	\textup{Dehn}\colon \Ls(V)
			\ar[r] 
    		&
		\displaystyle \bigoplus_{0 \subseteq U \subseteq V} \Ls(U) \otimes \Ls(V \ominus U) 
			\\
    	{(v_1,\ldots,v_n)}
			\ar[r,mapsto] 
			&
		\displaystyle 
			\sum_{\varnothing 
					\subseteq S 
					\subseteq \{v_1,\ldots,v_n\}
				 } 
			(\sgn(S)) 
			(S) \otimes (\pr_{\langle S \rangle ^\perp} S^c).
	\end{tikzcd}
\]
Here $\sgn(S) \in \{\pm 1\}$ is the sign of the bijection $\{v_1,\ldots,v_n\} \to S \amalg S^c$, written as a permutation of $\{1,\ldots,n\}$. 
Here's why the sign appears: suppose we order $\{v_1,\ldots,v_n\}$ to match the chosen orientation of $V$, and order $S$ and $S^c$ to match the orientations of $U$ and $V \ominus U$. 
Since the sets $\{v_1,\ldots,v_n\}$ and $S \amalg S^c$ have the same elements, these choices differ by a permutation, and we add the sign of that permutation. 
This guarantees that the map respects the orientations of our vector spaces in the correct way.

The antipode is the map
\[
	\begin{tikzcd}[row sep=0]
		\alpha \colon \Ls(V) 
			\ar[r]
			&
		\Ls(V)
			\\
		{(v_1,\ldots,v_n)} 
			\ar[r, mapsto]
			&
		{(-1)^n(v_1^\vee,\ldots,v_n^\vee),}
	\end{tikzcd}
\]
where 
\(
	v_i^\vee 
		= 
	\pr_{V\ominus
		 \langle
		 	\{v_1,\ldots,\hat{v_i},\ldots,v_n\}
		 \rangle
		} 
	(-v_i)
\)
as before. 
Note that the tuples $(v_1^\vee,\ldots,v_n^\vee)$ and $(v_1,\ldots,v_n)$ give orientations differing by a factor of $(-1)^n$, hence the sign.

\begin{lemma}
    This coproduct and antipode make the Lee Szczarba functor into a Hopf algebra in $\textup{Fun}(\Dip, \textup{Ab})$.
\end{lemma}

The proof is the observation that the Lee-Szczarba functor is isomorphic to the polytope functor, and so this follows from \cref{lem:polytope_functor_hopf}. 
Writing everything out explicitly in this case tends to be difficult, as the signs are harder to keep track of.

Taking the colimit, we conclude the result of Cathelineau \cite[Theorem 7.1.1]{cathelineau_sc} that the groups
\[ 
	\bigoplus_{n \geq 0} (\Ls(V) \otimes \det)_{O(n)}  
		= 
	\bigoplus_{n \geq 0} H_0(O(n);\Ls(V) \otimes \det) 
\]
form a Hopf algebra. 
In fact, Cathelineau's construction extends to the higher homology of $O(n)$ as well, although it does not obtain a Hopf algebra structure before taking coinvariants, which we will need to get the Hopf algebra structure on the underlying spectra.

%%%%%%%%%%%%%%%%%%%%%%%%%%%%%%%%%%%%
%%% SPECTRAL SAH FUNCTOR IS HOPF %%%
%%%%%%%%%%%%%%%%%%%%%%%%%%%%%%%%%%%% 
\section{The spectral Sah algebra is Hopf}\label{sec:hop_functor}

In this section we lift the combinatorial Hopf algebra structure of the previous section to the level of spectra. 
Our strategy will be to check that the spectral Sah functor is a Hopf algebra object in $\Fun(\Dip,\Sp^O)$ before passing to the colimit, as in the previous section. 
First, we check that the product and coproduct agree along the apartment maps.

\begin{proposition}\label{product_pi0}
    The product fits into the commuting diagram
    \[ 
    	\begin{tikzcd}[arrows=rightarrow, column sep=4em]
        	\tPt(U) \otimes \tPt(V) 
				\ar[d,"\textup{Join (\cref{pt_hopf})}"'] 
				\ar{r}{[\apt] \otimes [\apt]} 
				& 
			H_0(\Sigma^{-U}\Sigma\!\ST(U)) 
			\otimes H_0(\Sigma^{-V}\Sigma\!\ST(V)) 
				\ar{d}{\textup{\cref{cons:product}}}
				\\
	        \tPt(U \oplus V) 
	        	\ar{r}{[\apt]} 
				& 
			H_0(\Sigma^{-(U\oplus V)}\Sigma\!\ST(U \oplus V)).
    	\end{tikzcd}
	\]
\end{proposition}

\begin{proof}
    Choose orientations for $U$ and $V$ and any two elements $[u_1,\ldots,u_m] \in \tPt(U)$, $[v_1,\ldots,v_n] \in \tPt(V)$.
    Under this choice of convention, it suffices to show that the product of the apartment classes goes to the apartment class of the join:
    \[ 
    	\begin{tikzcd}[arrows=rightarrow,
			/tikz/column 1/.append style={anchor=base east},
    		/tikz/column 2/.append style={anchor=base west}
			]
        	[\apt(u_1,\ldots,u_m)] 
			\otimes [\apt(v_1,\ldots,v_n)] 
				\ar[r,phantom,"\in" description]
				& 
			H_m(\Sigma\!\ST(U)) \otimes H_n(\Sigma\!\ST(V)) 
				\ar[d, start anchor=-90, end anchor =26,"\textup{\cref{cons:product}}"] 
				\\
        	{[\apt(u_1,\ldots,u_m,v_1,\ldots,v_n)]} 
				\ar[r, phantom, "\in" description] 
				& 
			H_{m+n}(\Sigma\!\ST(U \oplus V)).
	    \end{tikzcd}
	\]
    We make the identification of $\Delta$-complexes
    \[ 
    	I^m \times I^n \cong I^{m+n}, 
	\]
    which takes each pair of subsets
    \[ 
    	(\varnothing \subseteq S \subseteq \{u_1,\ldots,u_m\}),
		(\varnothing \subseteq T \subseteq \{v_1,\ldots,v_n\}) 
	\]
    to the subset
    \[ 
    	\varnothing 
			\subseteq S \amalg T 
			\subseteq \{u_1,\ldots,u_m,v_1,\ldots,v_n\}. 
	\]
    With this isomorphism, the following diagram commutes:
    \[ 
    	\begin{tikzcd}[arrows=rightarrow, column sep=12em]
        	\frac{I^m}{\partial I^m}\sma\frac{I^n}{\partial I^n} 
				\ar{d}{\cong} 
				\ar{r}{\apt(u_1,\ldots,u_m)\, 
						\sma \, \apt(v_1,\ldots,v_n)
					  } 
				& 
			\Sigma\!\ST(U) \sma \Sigma\!\ST(V) 
				\ar{d}{\textup{\cref{cons:product}}} 
				\\
	        \frac{I^{m+n}}{\partial I^{m+n}} 
	        	\ar{r}{\apt(u_1,\ldots,u_m,v_1,\ldots,v_n)} 
				& 
			\Sigma\!\ST(U \oplus V).
	    \end{tikzcd}
	\]
    Both routes arise from the same map of simplicial complexes 
    \(
    	I^m \times I^n 
		\to |\varnothing \subseteq U \subseteq U \oplus V|,
	\) 
	sending the vertex corresponding to $(S,T)$ to the vertex corresponding to the vector space
    \[ 
    	\langle S \rangle \oplus \langle T \rangle 
		= \langle S \amalg T \rangle \subseteq U \oplus V. 
	\]
    This is enough to show that the product 
    \(
    	[\apt(u_1,\ldots,u_m)] \otimes [\apt(v_1,\ldots,v_n)]
	\) 
	gets mapped to $[\apt(u_1,\ldots,u_m,v_1,\ldots,v_n)]$ as required.
\end{proof}

\begin{proposition}\label{coproduct_pi0}
    The coproduct fits into the commuting diagram
    
\begin{adjustbox}{width=\textwidth}
\begin{tikzcd}[column sep=4em]
        \tPt(V) 
        	\ar[d,"\textup{Dehn (\cref{pt_hopf})}"'] 
			\ar{r}{[\apt]} 
			& 
		H_0(\Sigma^{-V}\Sigma\!\ST(V)) 
			\ar[d,"\textup{\cref{const:E_1_coalg}}"] 
			\\
        \displaystyle 
        	\bigoplus_{0 \subseteq U \subseteq V} 
				\tPt(U) \otimes 
				\tPt(V \ominus U) 
			\ar{r}{[\apt] \otimes [\apt]} 
			& 
		\displaystyle 
			\bigoplus_{0 \subseteq U \subseteq V} 
				H_0(\Sigma^{-U}\Sigma\!\ST(U)) \otimes 
				H_0(\Sigma^{-(V \ominus U)}
					\Sigma\!\ST(V \ominus U)).
    \end{tikzcd}
\end{adjustbox}
    
\end{proposition}

\begin{proof}
    Choose an orientation for $V$ and any element $[v_1,\ldots,v_n] \in \tPt(V)$ with the same orientation. 
    Give each span of each subset of $\{v_1,\ldots,v_n\}$ the orientation induced by its ordering inherited from $\{v_1,\ldots,v_n\}$. 
    Under this choice of convention, it suffices to show that the coproduct of the apartment class $[\apt(v_1,\ldots,v_n)]$ has the formula:
    
    \begin{adjustbox}{width=\textwidth}
    \begin{tikzcd}[column sep=0em]
        [\apt(v_1,\ldots,v_n)] 
        	\ar[r, phantom, "\in\hspace*{5mm}" description] 
			& 
		H_{\dim(V)}(\Sigma\!\ST(V)) 
			\ar[d,"\textup{\cref{const:E_1_coalg}}"] 
			\\
        {\displaystyle 
        	\sum_{S \subseteq \{v_1,\ldots,v_n\}} 
				\sgn(S)[\apt(S)] \otimes 
				[\apt (\pr_{V \ominus \langle S \rangle} S^c)]
		} 
			\ar[r, phantom, "\in" description] 
			& 
		\displaystyle 
			\bigoplus_{0 \subseteq U \subseteq V} 
				H_{\dim U}(\Sigma\!\ST(U)) \otimes 
				H_{\dim (V \ominus U)}(\Sigma\!\ST(V \ominus U)).
    \end{tikzcd}
    \end{adjustbox}

    We select the point in $D_1(2)$ corresponding to the point that cuts the interval exactly in half: the pair of intervals $[0,1/2]$ and $[1/2,1]$ in $[0,1]$. 
    The corresponding pair of embeddings is
    \[ 
    	e_1(x) = \frac12 x \quad e_2(x) = \frac12 + \frac12 x. 
	\]
    Then we describe the action of the coproduct for this point on the apartment class $I^n/\partial I^n \to \Sigma\!\ST(V)$ given by $(v_1,\ldots,v_n)$.

    The map $I^n \to |0 \subseteq \bullet \subseteq V|$ is a homeomorphism to a subcomplex consisting of $n!$ top-dimensional simplices and their faces, described as all flags of subspaces that are obtained from the vectors $v_1,\ldots,v_n$ by taking spans of their subsets. 
    We can represent each point in this cube using a step function, as in \cref{prop:set_model}. 
    However, instead of labeling sub-intervals $[s_j,s_{j+1}]$ with subspaces, we can label the cut points $s_j$ by vectors $v_i$, where the $v_i$ can appear in any order. 
    The corresponding step function then sends the interval $[s_j,s_{j+1}]$ to the subspace spanned by the vectors belonging to cut points $s_1, \dots, s_j$.

\begin{figure}[h]
\begin{center}
\begin{tikzpicture}[scale=2]
	% unit square
	\draw[thick] 
		(0,0) node[anchor=north east]{$\langle v_1,v_2 \rangle$}
			-- 
		(1,0) node[anchor=north west]{$\langle v_1 \rangle$} 
			-- 
		(1,1) node[anchor=south west]{$\langle\varnothing\rangle$} 
			-- 
		(0,1) node[anchor=south east]{$\langle v_2 \rangle$} 
			-- 
		cycle;
        \draw[thick] (0,0) -- (1,1);
        \draw[thick,fill,OI2] 
        	(1/3,2/3) circle (0.02) node[above]{$a$};
        \draw[thick,fill,OI1] 
        	(1/3,1/3) circle (0.02) node[above]{$b$};
        \draw[thick,fill,OI3] 
        	(4/5,1/5) circle (0.02) node[above]{$c$};
	
		% c
        \begin{scope}[xshift=2cm,scale=1.5]
        	\draw[ultra thick,OI3]
				(0,0) 
					edge[|-|] 
				(1/5,0) 
					edge[-|] 
				(4/5,0) 
					edge[-|] 
				(1,0);
			\node[scale=0.75] at (1/5,0.1) {$v_1$};
			\node[scale=0.75] at (4/5,0.1) {$v_2$};
			\node[OI3] at (-0.15,0) {$c = $};
			\draw[thick,|->] (1.1,0) -- (1.4,0);
			\draw[ultra thick,OI3,]
				(1.5,0) 
					edge[|-|] 
				(1.7,0)
					edge[-|] 
				(2.3,0)
					edge[-|]
				(2.5,0);
			\node[scale=0.75] at (1.6,0.1) {$0$};
			\node[scale=0.75] at (2,0.1) {$\langle v_1 \rangle$};
			\node[scale=0.75] at (2.4,0.1) 
				{$\langle v_1,v_2 \rangle$};
        \end{scope}

		% b
        \begin{scope}[xshift=2cm,scale=1.5,yshift=0.33cm]
        	\draw[ultra thick,OI1]
				(0,0) 
					edge[|-|] 
				(1/3,0) 
					edge[-|]  
				(1,0);
			\node[scale=0.75] at (1/3,0.1) {$v_1,v_2$};
			\node[OI1] at (-0.15,0) {$b = $};
			\draw[thick,|->] (1.1,0) -- (1.4,0);
			\draw[ultra thick,OI1]
				(1.5,0) 
					edge[|-|] 
				({1.5+1/3},0)
					edge[-|] 
				(2.5,0);
			\node[scale=0.75] at (5/3,0.1) {$0$};
			\node[scale=0.75] at (13/6,0.1) {$\langle v_1,v_2 \rangle$};
        \end{scope} 
        
        % a
        \begin{scope}[xshift=2cm,scale=1.5,yshift=0.66cm]
        	\draw[ultra thick,OI2]
				(0,0) 
					edge[|-|] 
				(1/3,0) 
					edge[-|] 
				(2/3,0) 
					edge[-|] 
				(1,0);
			\node[scale=0.75] at (1/3,0.1) {$v_1$};
			\node[scale=0.75] at (2/3,0.1) {$v_2$};
			\node[OI2] at (-0.15,0) {$a = $};
			\draw[thick,|->] (1.1,0) -- (1.4,0);
			\draw[ultra thick,OI2,]
				(1.5,0) 
					edge[|-|] 
				({1.5+1/3},0)
					edge[-|] 
				({1.5+2/3},0)
					edge[-|]
				(2.5,0);
			\node[scale=0.75] at (5/3,0.1) {$0$};
			\node[scale=0.75] at (2,0.1) {$\langle v_1 \rangle$};
			\node[scale=0.75] at (7/3,0.1) 
				{$\langle v_1,v_2 \rangle$};
        \end{scope} 
    \end{tikzpicture}
\end{center}
\caption{The map from a 2-cube to flags of subspaces from the proof of \cref{coproduct_pi0}.}
\end{figure}
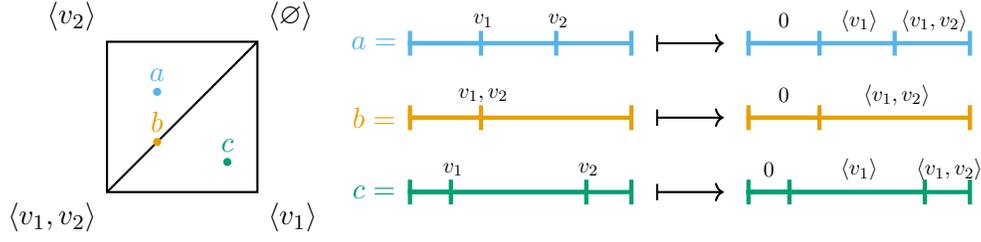

    Under this description, the cut points give the coordinates in the cube $I^n$, with the somewhat bizarre property that the vertex $(1,1,\ldots,1)$ is the initial vertex corresponding to the zero subspace, while the vertex $(0,0,\ldots,0)$ is the terminal vertex corresponding to all of $V$. 
    Each of the other vertices of the cube corresponds to a span of some subset $S \subseteq \{v_1,\ldots,v_n\}$, and its $i$th coordinate is 0 when $v_i \in S$ and $1$ when $v_i \not\in S$.
    
    Examining the formula for the coproduct (\cref{const:E_1_coalg}) for this particular choice of embeddings $(e_1,e_2) \in D_1(2)$, the coproduct on each point of $I^n$ vanishes when one of the coordinates (cut points) is equal to 1/2, and otherwise takes as $U$ the value of the function at $1/2$, takes the configuration the left of 1/2 and regards it as a configuration of subspaces of $U$, and takes the configuration to the right and projects the spaces to $V \ominus U$. 
    In particular, the image of this apartment lands in a finite sub-wedge of the wedge product $\bigvee_{0 \subseteq U \subseteq V} \Sigma\!\ST(U) \sma \Sigma\!\ST(V \ominus U)$, in which the only values of $U$ that can appear are those that are spans of the vectors $v_1,\ldots,v_n$.

    The map cuts the cube into $2^n$ sub-cubes with disjoint interiors, one for each subset $S \subseteq \{v_1,\ldots,v_n\}$, where the sub-cube is those points whose $i$-th coordinate is strictly smaller than $1/2$ when $v_i \in S$ and strictly greater than $1/2$ when $v_i \not\in S$. 
    The boundary of each of these sub-cubes is collapsed to the basepoint, and then within each one, we re-scale the coordinates to fit in $[0,1]$ and then apply the product of apartment maps for $S$ and for $\pr_{V \ominus \langle S \rangle}$ of $S^c$.
    This gives the commuting diagram
    \[ 
    \begin{tikzcd}[arrows=rightarrow, column sep=8em]
        \frac{I^n}{\partial I^n} 
        	\ar{d}{\textup{collapse}} 
			\ar{r}{\apt(v_1,\ldots,v_n)} 
			& 
		\Sigma\!\ST(V) 
			\ar{d}{\textup{\cref{const:E_1_coalg}}} 
			\\
        \bigvee_{\varnothing \subseteq 
        		 S \subseteq 
		 		 \{v_1,\ldots,v_n\}
				} 
        \frac{I^S}{\partial I^S} \sma 
		\frac{I^{S^c}}{\partial I^{S^c}} 
			\ar{r}{\apt(S) \sma 
				   \apt(\pr_{V \ominus \langle S \rangle} S^c)
				  } 
			& 
		\bigvee_{0 \subseteq U \subseteq V} 
			\Sigma\!\ST(U) \sma \Sigma\!\ST(V \ominus U),
    \end{tikzcd}
    \]
    where the vertical map on the left-hand side collapses the cube $I^n$ down to each of these $2^n$ sub-cubes. 
    The identification of each of these sub-cubes with the product $I^S \times I^{S^c}$ requires a permutation to be applied to the coordinates, and so this map on homology is a sum of maps with degrees $\sgn(S)$. 
    All together, this shows that the coproduct takes the homology class $[\apt(v_1,\ldots,v_n)]$ to the sum over all $S \subseteq \{v_1,\ldots,v_n\}$ of the classes $\sgn(S)([\apt (S)] \otimes [\apt (\pr_{V \ominus \langle S \rangle} S^c)])$, as required.
\end{proof}

\begin{remark}
	Although it is not necessary for our main theorem, it is also enlightening to write down a topological version of the antipode from \cref{pt_hopf}: we send each simplex corresponding to the flag of subspaces $0 \subsetneq U_0 \subsetneq \cdots \subsetneq U_k \subsetneq V$ to the simplex corresponding to the flag of complements $0 \subsetneq V \ominus U_k \subsetneq \cdots \subsetneq V \ominus U_0 \subsetneq V$, reversing the order of the vertices in the process. 
	It is straightforward to check that this takes the apartment of the basis $(v_1,\ldots,v_n)$ to the apartment of the basis $(v_1^\vee,\ldots,v_n^\vee)$, together with an orientation change of $(-1)^n$ because every coordinate direction of the cube is reversed. 
	Therefore this is a topological version of the involution.
	However, it is not clear how to lift the algebraic proof from \cref{lem:polytope_functor_hopf} to a topological proof that the diagram encoding the Hopf condition from \cite[Lemma 5.7]{KKMMW0} commutes with this choice of involution.
\end{remark}

\begin{proposition}\label{Hopfupourlife}
    The spectral Sah functor $\funSah$ is a $(\Comm,E_1)$-bialgebra in $\textup{Fun}(\Dip, \Sp^O)$ whose shear map is an isomorphism in the homotopy category. 
\end{proposition}

This proposition says that the spectral Sah functor $\funSah$ satisfies the most natural definition of a Hopf algebra for a model category, in which the shear map is not literally an isomorphism, only an isomorphism in the homotopy category \cite[Definition 5.10]{KKMMW0}. 
In particular, $\funSah$ becomes a Hopf algebra in the homotopy category of $\textup{Fun}(\Dip, \Sp^O)$.

\begin{proof}
    By \cref{cor:sah_functor_bialg}, the spectral Sah functor $\funSah$ is a $(\Comm,E_1)$-bialgebra in the category $\textup{Fun}(\Dip, \Sp^O)$. We will show that the shear map
    \[
    	\sh \colon 
    	\funSah\boxwedge \funSah 
			\xrightarrow{\id \, \boxwedge \, \Delta} 
		\funSah\boxwedge\funSah\boxwedge\funSah 
			\xrightarrow{\mu \,\boxwedge \, \id} 
		\funSah\boxwedge \funSah
	\]
    is an equivalence. 
    On $V \in \Dip$ this evaluates to a morphism
    \begin{equation}\label{eq:shear}
       	\sh_V \colon 
       	\bigvee_{U\subseteq V} 
			\Sigma^{-V}\Sigma ST(U) 
				\wedge\Sigma^{-U} \Sigma ST(V) 
			\to 
		\bigvee_{U\subseteq V} 
			\Sigma^{-V}\Sigma ST(U) 
				\wedge\Sigma^{-U} \Sigma ST(V). 
    \end{equation}
    The source and target of this map are wedges of spheres, so it suffices to prove that this is an isomorphism after applying $H_0$. 
    Moreover, we have $H_0\circ \funSah \cong \tPt$ via the apartment maps. 
    By \cref{lem:polytope_functor_hopf}, $\tPt$ is a Hopf algebra, so its shear map is an isomorphism. 
    Therefore \cref{eq:shear} is an equivalence.
\end{proof}

\subsection{Building the Hopf algebra spectrum}

It remains now to show that when we take the homotopy colimit along $\Dip$, we get a Hopf algebra in spectra. 
For this last step, it is necessary for us to work in $\infty$-categories, since there are very few coalgebras in the usual symmetric monoidal model categories of spectra. 
In fact, \cite[Theorem 1.1]{PS_coalgebras} shows that any strict coalgebra in one of these categories must be cocommutative. 
By contrast, we show in \cref{sec:no_cocomm} below that the classical Sah algebra (and therefore also the spectral Sah algebra) is not cocommutative.

Recall that every model category $M$ or simplicial model category $M_\bullet$ has an underlying $\infty$-category $\mathbfcal{M}$. 
This is formed by taking the nerve $N(M)$, or the simplicial nerve $N^s(M_\bullet)$, and applying a Dwyer-Kan localization to invert the weak equivalences, see e.g.~\cite[Section 4]{KKMMW0}. 
In the case of the model category of (orthogonal) spectra, the underlying $\infty$-category is the $\infty$-category of spectra, denoted in bold as $\Spp$ to distinguish it from any model category of spectra. 

In the companion paper \cite{KKMMW0} we recall how every symmetric monoidal simplicial model category $M_\bullet$ gives a symmetric monoidal $\infty$-category $\mathbfcal M$, and how every commutative bialgebra object in $M_\bullet$ gives a commutative bialgebra in $\mathbfcal M$, provided the underlying object of the algebra is ``cofibrant enough''.
We will use this result to pass the spectral Sah algebra to the underlying $\infty$-category of $\Fun(\Dip,\Sp^O)$, and \emph{then} take the colimit, in order to ensure that the colimit is a Hopf algebra.
Recall that $\Sp^O_+$ denotes the category of orthogonal spectra with the positive model structure introduced in \cref{sec-dip-diagrams}. 

\begin{lemma}\label{lem:HopfinInfty}
	The spectral Sah functor $\funSah$ defines a commutative Hopf algebra in the underlying symmetric monoidal $\infty$-category of $\Fun(\Dip,\Sp^O_+)$ or $\Fun(\Dip,\Sp^O)$.
\end{lemma}

\begin{proof}
	Recall from \cite[Definition 5.1]{KKMMW0} that every symmetric monoidal simplicial category $C_\bullet$ has an associated category of $(\Comm,D_1)$-bialgebras, in which the objects are objects of $C_\bullet$ with compatible structures of a commutative algebra and of a $D_1$-coalgebra.
	\cref{cor:sah_functor_bialg} proves that the spectral Sah functor $\funSah$ is a $(\Comm,D_1)$-bialgebra in the  symmetric monoidal simplicial category $C_\bullet = \Fun(\Dip,\Sp^O)^{\operatorname{pc}}$.
	Here the ``$\operatorname{pc}$'' refers to pointwise cofibrant diagrams of spectra. 

	By \cite[Theorem 5.7]{KKMMW0} it therefore gives a commutative bialgebra in the underlying $\infty$-category of $\Fun(\Dip,\Sp^O)$. 
	Furthermore, \cref{Hopfupourlife} shows that the shear map becomes an equivalence after we pass from $C_\bullet$ to the Dwyer-Kan localization $C_\bullet[W^{-1}]$ in which the pointwise stable equivalences have been inverted. 
	Therefore, again by \cite[Theorem 5.7]{KKMMW0}, the shear map of this bialgebra in the $\infty$-category is an equivalence, and so by definition this is a commutative Hopf algebra.
\end{proof}

\begin{remark}
Note that by a forthcoming result of Keenan--P\'eroux \cite{keenan_peroux}, the underlying symmetric monoidal $\infty$-category of $\Fun(\Dip,\Sp^O)$ is the symmetric monoidal $\infty$-category of diagrams $\Fun(N(\Dip),\Spp)$ under Day convolution. 
While this is nice to know, we can skip past this without getting into the technicalities of Day convolution, and show that the colimit is a Hopf algebra, as desired. 
Simply observe that the colimit defines a functor of symmetric monoidal 1-categories
\[ 
	\Fun(\Dip,\Sp^O) \to \Sp^O, 
\]
as in e.g. \cite{ben-moshe-schlank}, and we apply this to the projective cofibrant (not pointwise cofibrant!) diagrams, take the nerve, and invert the weak equivalences to get a functor of symmetric monoidal $\infty$-categories
\[ 
	N(\Fun(\Dip,\Sp^O)^{c}_0)[(W^{c})^{-1}]) 
		\to 
	N((\Sp^O)^c_0)[W^{-1}] \simeq \Spp. 
\]
Here we have implicitly used \cite[Proposition 4.7]{KKMMW0} to identify the underlying symmetric monoidal $\infty$-category of these symmetric monoidal model categories. 
By Lemma 5.6 of loc.\ cit.\ this yields an element of $\CHopf(\Spp)$.
\end{remark}

\begin{theorem} \label{cor:spectral_Sah_alg_is_Hopf}
	The resulting object of $\CHopf(\Spp)$ is a commutative Hopf algebra structure on the spectral Sah algebra $\Sah$ from \cref{spectral_sah_algebra}.
\end{theorem}

\begin{proof}
    Most of the work was already done. 
    We just need to observe that if we forget the Hopf algebra structure, then the above recipe replaces $\funSah$ by a cofibrant diagram and then takes its colimit. 
    This is equivalent to the homotopy colimit of $\funSah$, which is $\Sah$ by \cref{prop-hocolimS}.
\end{proof}

\begin{corollary}\label{cor:hopf_alg_struct_agree}
    As Hopf algebras, 
    \(
    	\pi_0(\Sah) \cong 
		\bigoplus_{n\geq 0} \widetilde{\mathcal{P}}(S^{n-1})
	\) 
	with Sah's original product and coproduct from \cref{pt_hopf}.
\end{corollary}

\begin{proof}
    In \cref{sec:hop_functor} we showed that $H_0 \circ \funSah \cong \widetilde{\Pt}$ as Hopf algebras in $\textup{Fun}(\Dip,\Ab)$. 
    After taking the colimit over $\Dip$, by definition this is the classical Sah algebra. 
    Therefore we also have that 
    \[
    	\pi_0(\Sah) 
			= \pi_0(\hocolim \funSah) 
			\cong \colim (\pi_0 \circ \funSah)
			\cong \colim \widetilde{\Pt}
	\]
is the classical Hopf algebra. For the last isomorphism above, we have used the fact that the functor $\funSah \colon \Dip \to \Sp^O$ takes values in connective spectra, so $\pi_0 \circ \funSah \cong H_0 \circ \funSah$. 
\end{proof}

	Concretely, this version of $\Sah$ is the spectrum $\underset{\Dip}\colim\, Q\funSah$, where $Q$ denotes cofibrant replacement in $\Fun(\Dip,\Sp^O)$. 
	The multiplication map agrees in the homotopy category with the composite
	\begin{align*}
    	\underset{\Dip}\colim\,Q\funSah 
	    	\wedge 
		\underset{\Dip}\colim\, Q\funSah
    		&\cong 
		\underset{\Dip}\colim\,
		\left( Q\funSah \boxwedge Q\funSah \right) 
			\\
	    	&\simeq 
		\underset{\Dip}\colim\,
		Q\left( \funSah \boxwedge \funSah \right) 
			\\
    		&\to 
		\underset{\Dip}\colim\, Q\funSah
	\end{align*}
	where the second map is the canonical identification in the homotopy category 
	\(
		Q\funSah \boxwedge Q\funSah 
			\simeq 
		Q(\funSah \boxwedge \funSah),
	\) 
	and the last map is the product 
	\(
		\funSah \boxwedge \funSah \to \funSah
	\) 
	from \cref{cons:product}. 
	Similarly, the coproduct is the composite
	\begin{align*}
    	\underset{\Dip}\colim\, Q\funSah
	    	&\to 
		\underset{\Dip}\colim\, 
		Q\left( \funSah \boxwedge \funSah \right) 
			\\
	    	&\simeq 
		\underset{\Dip}\colim\, 
		\left( Q\funSah \boxwedge Q\funSah \right) 
			\\
    		&\cong 
		\underset{\Dip}\colim\, Q\funSah \wedge 
		\underset{\Dip}\colim\, Q\funSah
	\end{align*}
	where the first map is the coproduct 
	\(
		\funSah \to \funSah \boxwedge \funSah
	\)
	from  \cref{const:E_1_coalg}. 
	The work we do in \cite{KKMMW0} may be summarized as a careful check that the point-set coherences between the product and coproduct that we proved in \cref{sec:Sah_functor_is_bialg} ascend to an elaborate infinite system of homotopy coherences between the above product and coproduct in the $\infty$-category of spectra.

%%%%%%%%%%%%%%%%%%%%%%%%%%%%%%%%%%
%%% FAILURE OF COCOMMUTATIVITY %%%
%%%%%%%%%%%%%%%%%%%%%%%%%%%%%%%%%% 
\section{The failure of cocommutativity}
\label{sec:no_cocomm}

It is natural to ask whether the $(E_\infty,E_1)$-Hopf algebra structure from the previous section can be improved. 
For instance, is it possible to refine the comultiplication to one that is $E_2$? 
In this section we answer this question in the negative, by showing that the classical Sah algebra fails to be cocommutative. 
This in turn prevents the spectral Sah algebra's comultiplication from being more than $E_1$.

Our counterexample occurs in the first degree where cocommutativity can possibly fail:

\begin{theorem}\label{not_cocomm}
    The Dehn invariant map 
    \(
    	\tcP(S^3) \to \tcP(S^1) \otimes \tcP(S^1)
	\) 
	is not co-commutative.
\end{theorem}

We prove this by considering a regular tetrahedron $T_a \subseteq S^3$ of side length $a$, with $0 < a < \arccos(-\frac{1}{3})$. 
Let $D$ be the dihedral angle, that is, the angle formed between any two of the two-dimensional faces of $T_a$. 
The Dehn invariant of $T_a$ is
\[ 
	6(a \otimes D) \in (\R/\pi\Q) \otimes (\R/\pi\Q). 
\]
It suffices to show that for at least one value of $a$, this element is not preserved by the symmetry isomorphism: $a \otimes D \neq D \otimes a$.

To show this we first give an explicit formula for $D$ in terms of $a$. 
As is typical in spherical geometry, we conflate the length $a$ with the angle it subtends at the center of the sphere; working on a sphere of radius $1$, this angle (in radians) is exactly the length $a$.

\begin{lemma}\label{a_to_d}
    The dihedral angle $D$ of the spherical regular tetrahedron of side length $a$ satisfies
    \[ 
    	\cos(D) = \frac{\cos (a)}{1 + 2\cos (a)}. 
	\]
\end{lemma}

\begin{proof}
    This is an exercise in spherical trigonometry. 
    Each face of the tetrahedron is an spherical equilateral triangle of side length $a$. Let $A$ be the angle at each vertex of this triangle, as in \Cref{trig_figure}. 
    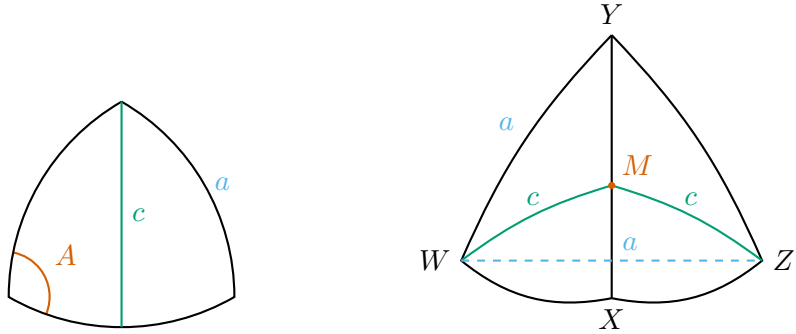
\begin{figure}[h]     
    \begin{center}
    \begin{tikzpicture}[scale=2]
        \draw[OI6,thick] 	
        	(-0.5,-0.41) arc (-22:79:0.3) 
				node[midway,anchor=south west]{$A$}; 
        \draw[thick] 
        	(90:1) 
				arc (120:180:1.5) 
				arc (-120:-60:1.5) 
				arc (0:60:1.5) 
					node[midway,anchor=west,color=OI2]{$a$};
        \draw[thick,OI3] 
        	(90:1) -- node[right]{$c$} (0,-0.5);
    \end{tikzpicture}
    \hspace*{2cm}
    \begin{tikzpicture}[scale=2]
        \draw[thick] 
        	(0,-0.25) node[below]{$X$} 
				-- 
			(0,1.5) node[above]{$Y$};
        \draw[thick] 
        	(0,-0.25) 
				to[bend right=25] 
			(1,0) node[right]{$Z$} 
				to[bend right=10] 
			(0,1.5);
        \draw[thick] 
        	(0,-0.25) 
				to[bend left=25] 
			(-1,0) node[left]{$W$} 
				to[bend left=10] 
				node[anchor = south east,color=OI2]{$a$} 
			(0,1.5);
        \draw[thick,OI3] 
        	(-1,0) 
				to[bend left=10] 
				node[above]{$c$} 
			(0,0.5);
        \draw[thick,OI3] 
        	(0,0.5) 
				to[bend left=10] 
				node[above]{$c$} 
			(1,0);
         \draw[fill,OI6] 
         	(0,0.5) 
				circle (0.02) 
				node[anchor = south west]{$M$};
         \draw[thick,dashed,OI2] 
         	(-1,0) 
				to[bend right=0] 
				node[anchor = south west]{$a$} (1,0);
%         \draw[thick] (-0.35,-0.278) -- (0.35,-0.278);
    \end{tikzpicture}
    \end{center}
    \caption{An illustration of the angles and lengths used to prove \cref{a_to_d}.}
    \label{trig_figure}
    \end{figure} 
    
    By the spherical law of cosines, we have 
    \(
    	\cos(a) = \cos^2(a) + \sin^2(a) \cos(A).
    \)
    Solving for $\cos(A)$ yields
    \[
        \cos (A) 
        	= \frac{\cos (a) - \cos^2 (a)}{\sin^2 (a)}
	        = \frac{\cos (a) - \cos^2 (a)}{1 - \cos^2 (a)}
	        = \frac{\cos (a)}{1+\cos (a)}. 
	\]
	Therefore, 
	\[
		\sin^2 (A) 
			= 1 - \left(\frac{\cos (a)}{1+\cos (a)}\right)^2
        	= \frac{1+2\cos (a)+\cos^2 (a)}{(1+\cos (a))^2} 
			- \frac{\cos^2 (a)}{(1+\cos (a))^2}
        	= \frac{1+2\cos (a)}{(1+\cos (a))^2}.
	\]
    Let $c$ be the altitude of this equilateral triangle.  By the spherical law of sines, we have:
    \[
    	\frac{\sin(c)}{\sin(A)} = \frac{\sin(a)}{1}. 
    \]
    Solving for $\sin(c)$ and combining with the equation for $\sin^2(A)$ above, we see that
	\[
        \sin^2 (c) 
        	= \sin^2 (a) \sin^2 (A)
	        = (1-\cos^2 (a))
	       	  \left(\frac{1+2\cos (a)}{(1+\cos (a))^2}\right)
    	    = \frac{(1-\cos (a))(1+2\cos (a))}{1+\cos (a)}.
	\]
    Now take two adjacent faces $WXY$ and $XYZ$ of the tetrahedron $T_a$ that meet at the common edge $XY$, and drop altitudes from $W$ to $XY$ and $Z$ to $XY$, meeting at the midpoint $M$ of $XY$. 
    This is illustrated in \Cref{trig_figure}. 
    The triangle $WZM$ has side lengths $a$, $c$, and $c$, and the dihedral angle $D$ of the tetrahedron is the angle between the sides of length $c$. 
    By the spherical law of cosines applied to $WZM$ we have:
    \begin{align*}
        \cos (a) 
        	&= 
		\cos^2 (c) + \sin^2 (c) \cos (D) 
			\\
       		&= 
		1 + \sin^2 (c)(\cos (D) - 1) 
			\\
        	&= 
		1 + \left(
				\frac{(1-\cos (a))(1+2\cos (a))}{1+\cos (a)}
			\right)
		(\cos (D) - 1).
    \end{align*}
	Subtracting $1$ from both sides, we see that
	\[
        -(1 - \cos (a)) 
        	= 
		\frac{(1-\cos (a))(1+2\cos (a))}%
			 {1+\cos (a)}(\cos (D) - 1).
	\]
	Cancelling out the factor $1 - \cos(a)$ from both sides and solving for $\cos(D)$ then yields
    \begin{align*} 
        \cos (D) 
        	&= 
		1 + \frac{-(1+\cos (a))}{(1+2\cos (a))}
			\\
        	&= 
		\frac{1+2\cos (a)}{1+2\cos (a)} 
		- \frac{1+\cos (a)}{1+2\cos (a)}
			\\[.5em]
         	&= 
		\frac{\cos (a)}{1+2\cos (a)} \qedhere
    \end{align*}
\end{proof}

As a sanity check, the limiting case of $a \to 0$ gives $A = \frac{\pi}{3}$ and $D = \arccos(\frac{1}{3})$, which agrees with the dihedral angle of a Euclidean regular tetrahedron. 
In the case of a right-angled spherical tetrahedron, $a = A = D = \frac{\pi}{2}$. 
Finally, the largest possible regular tetrahedron has all faces lying along an equatorial $S^2$, giving $a = \arccos(-\frac{1}{3})$, $A = \frac{2\pi}{3}$, and $D = \pi$.

We next establish an elementary analytic lemma.

\begin{lemma}
    If $f\colon (t_0,t_1) \to \R$ is a non-constant analytic function, then the preimage of any one point is countable.
\end{lemma}

\begin{proof}
    The zeros of $f'$ are isolated, because otherwise by analytic continuation $f'$ would have to be identically zero.
    Away from those zeros, $f$ consists of countably many order-preserving bijections from intervals to intervals. 
    Any one point in the codomain has at most one preimage in each of these intervals, and therefore has countably many preimages in total.
\end{proof}

\begin{corollary}
\label{preimageofQ}
     If $f\colon (t_0,t_1) \to \R$ is a non-constant analytic function, then $f^{-1}(\Q)$ is countable.
\end{corollary}

\begin{proposition}\label{cocomm_counterexample}
    For all but countably many values of $a$, 
    \(
    	a \otimes D \neq D \otimes a.
	\)
\end{proposition}

\begin{proof}
	Regarding $D$ as a function of $a$ on the interval $a \in \big(0,\arccos(-\frac13)\big)$, the following three functions are all analytic and non-constant by \cref{a_to_d}:
	\[ 
		\frac{a}{\pi}, \ \frac{D}{\pi}, \ \frac{D}{a}. 
	\]
	By \cref{preimageofQ}, there are only countably many values of $a$ for which at least one of them is rational. 
	For the remaining values of $a$, all three quantities are irrational. 
	It follows that $a$ and $D$ are both nonzero elements of the rational vector space $\R/\pi\Q$, and that they are linearly independent. 
	Since they are linearly independent, there is a basis for $(\R/\pi\Q) \otimes (\R/\pi\Q)$ in which $a \otimes D$ and $D \otimes a$ are distinct basis vectors, and therefore they are non-equal.
\end{proof}

This finishes the proof of \cref{not_cocomm}.

%%%%%%%%%%%%%%%%%%%%%%%%%%%%%%%%%%%%%%%%%%%%%%%%%%%%%%%%%%%
%%% APPLICATION: NONTRIVIAL ELEMENTS IN HIGHER K-GROUPS %%%
%%%%%%%%%%%%%%%%%%%%%%%%%%%%%%%%%%%%%%%%%%%%%%%%%%%%%%%%%%%
\section{Application: non-trivial elements in higher K-groups}\label{sec:application}

Finally, we turn our attention to the rational homotopy groups of the spectral Sah algebra, which form a bigraded Hopf algebra
\begin{equation}\label{rational_sah}
	\bigoplus_{m \geq 0} \pi_m(\Sah)_\Q 
		\cong 
	\bigoplus_{m,n \geq 0} \ 
		\widetilde K_m(\mathcal{P}^{S^{n-1}}_{O(n)}) \otimes \Q,
\end{equation}
%which can equivalently be described as $\bigoplus_{m,n \geq 0} H_m\big( O(n) ; \St(\R^n)^t) \otimes \Q$; see \cite{scissors_thom}. Building on previous calculations, we construct a large nonzero subalgebra.
In a few cases, we can understand the integral homotopy groups of the spectral Sah algebra.  

\begin{example}
The degree $n$ piece of the spectral Sah algebra is 
\(
	\widetilde{K}(\mathcal{P}^{S^{n-1}}_{O(n)}) 
		\simeq 
	\ST(\R^n)_{hO(n)}^{1 - \R^n}.
\) 
When $n = 1$, this is equivalent to
\[
	\Sph^{1-\sigma}_{hO(1)} 
		\simeq 
	\Sigma^{1-\sigma}_+ \mathbb{RP}^\infty,
\] 
where $\sigma$ denotes the sign representation, which corresponds to the canonical bundle over $\mathbb{RP}^\infty$.
This spectrum vanishes rationally, but not integrally. 
On integral homology, it has a $\Z/2$ in every even degree and a 0 in every odd degree.
\end{example} 

The typical calculation, however, requires inverting some elements.

\begin{lemma}\label{odd_vanishing}
    When $n$ is odd, the reduced spherical scissors congruence $K$-theory spectrum $\widetilde K(\Pol{S^{n-1}}{O(n)})$ has trivial homology after 2 is inverted.
\end{lemma}

In particular, its rationalization is trivial when $n$ is odd. So the bigraded Hopf algebra \cref{rational_sah} is concentrated in degrees where $n$ is even.

\begin{proof}
    This follows from the formula
    \[ 
    	H_m\big(\widetilde K(\mathcal{P}^{S^{n-1}}_{O(n)});\Z[1/2]\big) 
			\cong 
		H_m\big( O(n) ; \St(\R^n)^t \otimes \Z[1/2]\big), 
	\]
    which is the reduced version of the formula in \cite[Theorem 1.5]{scissors_thom} and is an immediate consequence of \cite[Theorem 1.9]{scissors_thom}. 
    The $(-)^t$ decoration adds a minus sign for any element that reverses orientation.
    
    The element $-\id \in O(n)$ is central and acts trivially on the Steinberg module $\St(\R^n)$. 
    Since $n$ is odd, it also reverses orientation, so in total it acts by multiplication by $-1$ on the coefficients in this homology group. 
    By the ``center-kills'' lemma from \cite[Lemma 5.4]{dupont_book}, therefore multiplication by $2$ acts trivially on the homology. 
    Note that a version of this calculation (with the same proof) also appears in \cite[1.23]{cz-hilbert}.
\end{proof}

\begin{lemma}\label{rational_r0_and_r2}
    For $n = 0$ and $n = 2$, the rational reduced spherical scissors congruence groups are given by
    \[
        \widetilde K_m(\mathcal{P}^{S^{-1}}_{O(0)}) \otimes \Q 
        	\cong 
			\begin{cases} 
				\Q & m = 0 \\ 
				0 & m > 0, 
			\end{cases} 
		\hspace*{1cm}
        \widetilde K_m(\mathcal{P}^{S^{1}}_{O(2)}) \otimes \Q 
        	\cong 
			\begin{cases} 
				\Lambda^{m+1}(\R/\Q) & \textup{$m$ even} \\ 
				0 & \textup{$m$ odd}. 
			\end{cases}
    \]
\end{lemma}

\begin{proof}
    The $n = 0$ case follows from \cite[Example 7.1]{scissors_thom}, which says that before rationalization $\widetilde K(\mathcal{P}^{S^{-1}}_{O(0)}) \simeq \Sph$. The $n = 2$ case follows from \cite[Example 7.7]{scissors_thom}, which says that $\widetilde{K}(\mathcal{P}^{S^{1}}_{O(2)})$ is the homotopy fiber of a map of the form
\[ 
	\Sigma^\infty_+ B(\Z/2 \times \Z/2) 
		\to 
	\Sigma^\infty_+ BO(2),
\]
where $O(2)$ has the discrete topology and the Klein-four subgroup is the subgroup that stabilizes a single line through the origin.\end{proof}

We can use these calculations and the Hopf algebra structure of $\Sah$ to construct a nonzero subalgebra. Recall that an element of a Hopf algebra $x \in A$ is \emph{primitive} if its coproduct satisfies
\[ 
	\Delta(x) = 1 \otimes x + x \otimes 1. 
\]
Note that sums and scalar multiples of primitive elements are also primitive. 

\begin{proposition}
    If $A$ is any graded commutative Hopf algebra over the rationals, $V \subseteq A$ is any sub-vector space consisting entirely of primitive elements, and $\Lambda(V)$ is the free graded-commutative algebra on $V$, then the natural map $\Lambda(V) \to A$ is injective.
\end{proposition}

\begin{proof}
    Let $A'$ be the subalgebra of $A$ generated by $V$. 
    Then by \cite[Chapter 3, Corollary 3.2]{may_hopf}, which is a corollary of \cite[Theorem 5.18]{milnor_moore}, $A'$ is isomorphic to the free (graded-)commutative algebra on its primitive elements. 
    These contain $V$, but we argue that they must coincide with $V$, for otherwise $A'$ is the free commutative algebra on some larger vector space $V \oplus W$. 
    But if $W \neq 0$ then $A' \cong \Lambda(V \oplus W)$ is not actually generated by $V$ as an algebra, which is a contradiction.
\end{proof}

In the spectral Sah algebra, the classes in $\widetilde K_{2k}(S^1) \cong \Lambda^{2k+1}(\R/\Q)$ are all primitive for degree reasons. 
Therefore the free commutative algebra on these classes is a subalgebra of the spectral Sah algebra.  
In this case, we only have primitive elements of even degree, so the free graded-commutative algebra $\Lambda(-)$ becomes the free commutative algebra or polynomial algebra $P(-)$.
The resulting classes in this subalgebra are given geometrically by iterated joins of arcs lying in perpendicular copies of $S^1$ inside $S^{2d-1}$, and higher scissors congruence classes that come from applying interval exchange transformations to these arcs.
    
\begin{corollary}\label{cor:nonzero_subalg}
    The rational homotopy groups of the spectral Sah algebra \cref{rational_sah} contain as a subalgebra the free commutative algebra
    \[ 
    	P\left( 
			\bigoplus_{\textup{$m$ even}} 
				\Lambda^{m+1}(\R/\Q)[m,2] 
		 \right), 
	\]
    where the $[m,2]$ means that the term appears in bidegree $(m,n) = (m,2)$.
\end{corollary}

Note that this subalgebra is both commutative and cocommutative, so the class of the regular tetrahedron $T_a$ that we considered in \cref{sec:no_cocomm} is a class in bidegree $(0,4)$ that for almost every value of $a$ does not lie in this subalgebra.

\begin{figure}[h]
\begin{tikzpicture}[xscale=2.5,yscale=1]
	% grid
	\begin{scope}[xshift=-0.7cm, yshift=-0.7cm]
%		\draw[gray] (-0.25,-0.25) grid (5.25,7.25); 
		\draw[->,thick] (-0.25,0) to (5.25,0) node[right]{$m$}; 
		\draw[->,thick] (0,-0.5) to (0,7.25) node[above]{$n$};

        % axes labels 	
    	\foreach \x in {0,...,4}{
	   	\draw[thick] (\x+0.7,0.1) --     (\x+0.7,-0.1) node[below]{$\x$};
		}
    	\foreach \y in {0,...,6}{
	   	\draw[thick] (0.04,\y+0.7) -- (-0.04,\y+0.7) node[left]{$\y$};
		}
	\end{scope}	
     
	% row 0
	\node at (0,0) {$\Q$};
	\node at (0,2) {\Large$\sfrac{\R}{\Q}$};
	\node at (0,4) {$\Sym^2\!\left(\sfrac{\R}{\Q}\right)$};
	\node at (0,6) {$\Sym^3\!\left(\sfrac{\R}{\Q}\right)$};
	
	% row 2
	\node at (2,2) {$\Lambda^3\!\left(\sfrac{\R}{\Q}\right)$};
	\node at (2,4) {$\sfrac{\R}{\Q}\! \otimes\! \Lambda^3\!\left(\sfrac{\R}{\Q}\right)$};
	\node[fill=white] at (2,6) {$\Sym^2\!\left(\sfrac{\R}{\Q}\right) \otimes \Lambda^3\!\left(\sfrac{\R}{\Q}\right)$};
	
	% row 4
	\node at (4,2) {$\Lambda^5\!\left(\sfrac{\R}{\Q}\right)$};

	\node at (4,4) {$\oplus$};	
	\node at (4,4.4) {$\Sym^2\!\Lambda^3\!\left(\sfrac{\R}{\Q}\right)$};
	\node at (4,3.5) {$\sfrac{\R}{\Q}\! \otimes\! \Lambda^5\!\left(\sfrac{\R}{\Q}\right)$};

	\node at (4,6.25) {$\oplus$};		
	\node at (4,6.65) { $\sfrac{\R}{\Q} \otimes \Sym^2\!\Lambda^3\!\left(\sfrac{\R}{\Q}\right)$};

	\node at (4,5.75) { $\Sym^2\!\left(\sfrac{\R}{\Q}\right) \otimes \Lambda^5\!\left(\sfrac{\R}{\Q}\right)$};
\end{tikzpicture}
\caption{The subalgebra of \cref{cor:nonzero_subalg} up to bidegree $(4,6)$. Here, $\Sym^n$ stands for the $n$-th symmetric power. Unfilled bidegrees are zero.}
\end{figure}
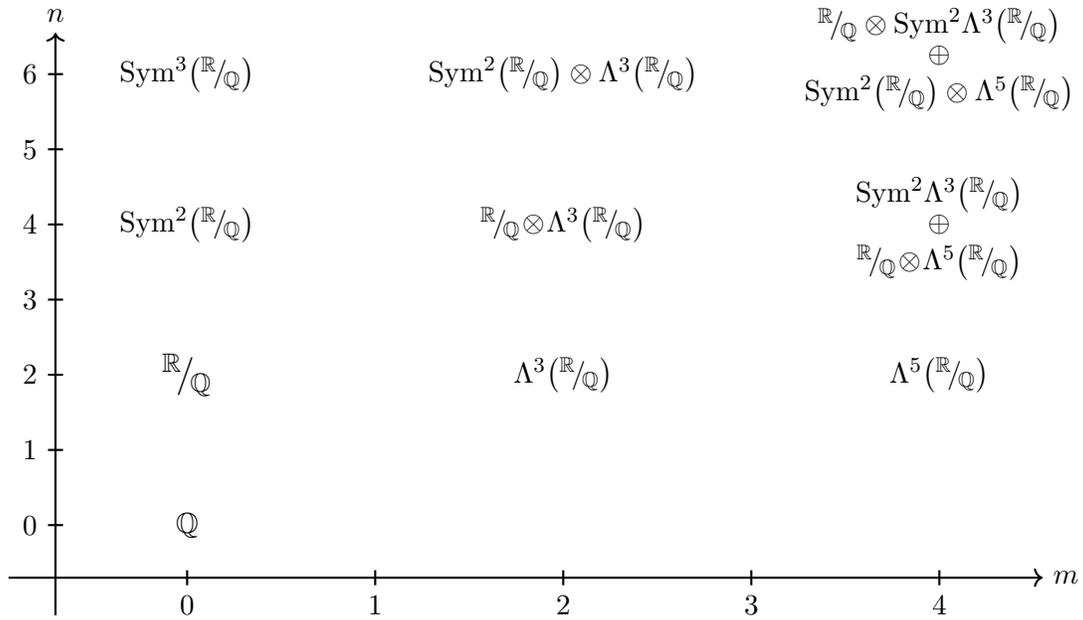

%%%%%%%%%%%%%%%%%%%%
%%% BIBLIOGRAPHY %%%
%%%%%%%%%%%%%%%%%%%%

\newpage

\bibliographystyle{alpha}
\bibliography{references}

\end{document}